\documentclass[11pt]{article}

\usepackage{natbib}

\usepackage{bbm}
\usepackage{epsfig,lscape}
\usepackage{amssymb,amsfonts,amsmath}

\usepackage{enumitem}
\usepackage{amsthm}
\usepackage{tikz}
\usetikzlibrary{arrows}

\usepackage{booktabs}
\usepackage{rotating}
\usepackage[referable]{threeparttablex}

\usepackage{color}

\newcommand{\tzq}{\textcolor{black}}

\def\me{\mathrm e}

\def\dif{\mathrm d}
\def\diag{\mathrm{diag}}

\def\T{ {\mathrm{\scriptscriptstyle T}} }

\def\bbR{\mathbb R}
\def\dom{\mathrm{dom}}
\def\argmax{\mathrm{argmax}}
\def\argmin{\mathrm{argmin}}

\def\one{\mathbbm{1}}

\def\dprime{ {\prime\prime} }

\newenvironment{prf}{\par\noindent{\bf Proof.}}{\hfill$\blacksquare$\\[2mm]}

\newtheorem{lem}{Lemma}
\newtheorem{pro}{Proposition}
\newtheorem{cor}{Corollary}

\theoremstyle{definition}

\theoremstyle{definition}
\newtheorem{rem}{Remark}

\begin{document}

\def\spacingset#1{\renewcommand{\baselinestretch}%
{#1}\small\normalsize} \spacingset{1.5}

\begin{titlepage}

\begin{center}
{\Large On Loss Functions and Regret Bounds for Multi-category Classification}

\vspace{.1in} Zhiqiang Tan\footnotemark[1] and Xinwei Zhang\footnotemark[1]

\vspace{.1in}
\today
\end{center}

\footnotetext[1]{Department of Statistics, Rutgers University. Address: 110 Frelinghuysen Road,
Piscataway, NJ 08854. E-mail: ztan@stat.rutgers.edu, xinwei.zhang@rutgers.edu.}

\paragraph{Abstract.}
We develop new approaches in multi-class settings for constructing proper scoring rules and hinge-like losses and establishing corresponding regret bounds
with respect to the zero-one or cost-weighted classification loss.
Our construction of losses involves deriving new inverse mappings from a concave generalized entropy to a loss through the use of a convex dissimilarity function
related to the multi-distribution $f$-divergence. Moreover, we identify new classes of multi-class proper scoring rules, which also recover and reveal interesting relationships between
various composite losses currently in use. We establish new classification regret bounds in general
for multi-class proper scoring rules by exploiting the Bregman divergences of the associated generalized entropies, and, as applications, provide simple meaningful
regret bounds for two specific classes of proper scoring rules.
Finally, we derive new hinge-like convex losses, which are tighter convex extensions than related hinge-like losses and geometrically
simpler with fewer non-differentiable edges, while achieving similar regret bounds. \tzq{We also establish a general classification regret bound for all losses which induce
the same generalized entropy as the zero-one loss.}

\vspace{.1in}

\end{titlepage}

\section{Introduction}\label{sect:intro}

Multi-category classification has been extensively studied in machine learning and statistics. For concreteness, let $\{ (X_i,Y_i): i=1,\ldots,n\}$ be training data generated from a certain
probability distribution on $(X,Y)$, where $X$ is a covariate or feature vector and $Y$ is a class label, with possible values from 1 to $m \, (\ge 2)$.
Various learning methods are developed in the form of minimizing an empirical risk function,
\begin{align}
\hat R_L (\alpha) = \frac{1}{n} \sum_{i=1}^n L (Y_i, \alpha(X_i)), \label{eq:emp-risk}
\end{align}
where $L(y, \alpha(x))$ is a loss function, and $\alpha(x)$ is a vector-valued function of covariates, taken from a potentially rich family of functions, for example, reproducing kernel Hilbert spaces
or neural networks. For convenience, $\alpha(x)$ is called an action function, following the terminology of decision theory \citep{degroot1962,grunwald2004}.
The performance of $\alpha(x)$ is typically evaluated by the zero-one risk on test data,
\begin{align}
E \big\{ L^{\text{zo}} (Y_0, \tilde \alpha (X_0)) \big\}, \label{eq:test-risk}
\end{align}
where $(X_0,Y_0)$ is a new observation, independent of training data, and $\tilde\alpha(x)$ is either $\alpha(x)$ or an $m$-dimensional vector converted from $\alpha(x)$,
and $L^{\text{zo}}(y, \tilde\alpha(x))$ is the zero-one loss, defined as 0 if the $y$th component of $\tilde\alpha(x)$ is a maximum  and 1 otherwise.
Due to discontinuity, using $L^{\text{zo}}$ directly as $L$ in (\ref{eq:emp-risk}) is computationally intractable. Hence
the loss $L$ used in (\ref{eq:emp-risk}) is also referred to as a surrogate loss for $L^{\text{zo}}$.

It is helpful to distinguish two types of loss functions $L$ commonly used for training in (\ref{eq:emp-risk}).
One type of losses, called scoring rules, involves an action defined as a probability vector $q(x): \mathcal X \to \Delta_m$, where
$\mathcal X$ is the covariate space and $\Delta_m$ is the probability simplex with $m$ categories \citep{savage1971,buja2005}.
The elements of $q(x)$ can be interpreted as class probabilities.
Typically, the probability vector $q(x)$ is parameterized in terms of a real vector $h(x)$ as $q^h(x)$, via an invertible link such as the multinomial logistic (or softmax) link.
The resulting loss $L(y,q^h(x))$ is then called a composite loss, with $h(x)$ as an action function \citep{williamson2016}.
It is often desired to combine a {\it proper} scoring rule, which ensures infinite-sample consistency of probability estimation (see Section~\ref{sec:scoring}),
with a link function $q^h$ such that $L(y,q^h(x))$ is convex in $h$.
In multi-class settings, composite losses satisfying these properties include the standard likelihood loss and two variants of exponential losses related to boosting \citep{zou2008,mukherjee2013},
all combined with the multinomial logistic link.

Another type of losses involves an action defined as a real vector, $\alpha(x): \mathcal X \to \bbR^m$, where
the elements of $\alpha(x)$, loosely called margins, can be interpreted as relative measures of association of $x$ with the $m$ classes.
Although a composite loss $L(y,q^h(x))$ based on a scoring rule can be considered with $h(x)$ as a margin vector,
it is mainly of interest to include in this type hinge-like losses, where the margins are designed not to be directly mapped to probability vectors.
The hinge loss is originally related to support vector machines in two-class settings.
This loss, $L^{\text{hin}}(y, \tau(x))$, is known to be convex in its action $\tau$, and achieve classification calibration (or infinite-sample classification consistency),
which means that a minimizer of the hinge loss in the population version
leads to a Bayes rule minimizing the zero-one risk (\ref{eq:test-risk}) \citep{lin2002,zhang2004b,bartlett2006}.
There are various extensions of the hinge loss to multi-class settings.
The hinge-like losses in \cite{lee2004} and \cite{duchi2018} are shown to achieve classification calibration,
whereas those in \cite{weston1998} and \cite{crammer2001} fail to achieve such a property \citep{zhang2004a,tewari2007}.

Classification calibration is also called Fisher consistency, although it is appropriate to distinguish two types of Fisher consistency,
in parallel to the two types of losses above:
Fisher probability consistency as satisfied by a proper scoring rule or
Fisher classification consistency as achieved by a hinge-like loss.
In general, Fisher probability consistency (or properness) implies classification consistency, but not vice versa \citep{williamson2016}.
On the other hand, there are interesting results indicating that only hinge-like losses are (classification) consistent with respect to the zero-one loss
when both an action function and a data quantizer are estimated \citep{nguyen2009,duchi2018}.

The purpose of this article is broadly two-fold: first constructing new multi-class losses while studying existing ones, and second
establishing corresponding classification regret bounds.
Our development in both directions is facilitated by the concept of a generalized entropy, defined as the minimum Bayes risk for a loss \citep{grunwald2004}.
In particular, our construction of losses (including proper scoring rules
and hinge-like losses) involves deriving inverse mappings from
a (concave) generalized entropy to a loss through a (convex) dissimilarity function $f$, which is motivated by the multi-distribution $f$-divergence
\citep{gyorfi1978, garcia2012}.
Moreover, our regret bounds for proper scoring rules are directly determined by the Bregman divergence of the associated generalized entropy.
\tzq{Finally, we provide a general characterization of losses with the same generalized
entropy as the zero-one loss and establish a general classification regret bound for all such losses, beyond the hinge-like losses
specifically constructed.}

\vspace{.1in}\textbf{Summary of main results.} We develop a new general representation of proper scoring rules based on a dissimilarity function $f$
directly related to the multi-distribution $f$-divergence (Propositions~\ref{pro:f-loss2}--\ref{pro:f-loss3}).
The general representation is then employed to derive two new classes of proper scoring rules:
multi-class pairwise losses corresponding an additive dissimilarity function $f$ and
multi-class simultaneous losses with non-additive dissimilarity functions corresponding to generalized entropies defined as $L_\beta$ norms (Section~\ref{sec:scoring-examples}).
These two classes of losses not only reveal interesting relationships between
the likelihood and exponential losses mentioned above, but also lead to specific new losses including
a pairwise likelihood loss distinct from the standard likelihood loss.

We propose a novel approach for constructing hinge-like, convex losses in multi-class settings: we first derive a new loss with actions restricted to the probability simplex $\Delta_m$
and its generalized entropy identical to that of the zero-one or cost-weighted classification loss (Proposition~\ref{pro:f-loss}; Lemma~\ref{lem:cw-conj2}),
and then we find a convex extension of the loss such that its actions are defined on $\bbR^m$ and its generalized entropy remains unchanged.
As a result, we derive two new hinge-like losses (Propositions~\ref{pro:cw-conj3}--\ref{pro:cw-conj4}), related to \cite{lee2004} and \cite{duchi2018} respectively.
In each case, our new loss and the existing one admit the same generalized entropy and coincide with each other for actions restricted to $\Delta_m$,
but our loss is uniformly lower (hence a tighter extension outside the probability simplex) and geometrically simpler with fewer non-differentiable edges.

We establish classification regret bounds not only for multi-class proper scoring rules in general (Section~\ref{sec:regret-bd}),
but also for our new hinge-like losses (Proposition~\ref{pro:hinge-regret-bd}) and \tzq{more broadly {\it all} losses with the same generalized entropy as the zero-one loss (Proposition~\ref{pro:general-hinge-bd}).}
In each case, a regret bound compares the regret of the loss studied with that of the zero-one or cost-weighted classification loss and implies that classification calibration is achieved
with a quantitative guarantee. As applications in the first case, we derive simple meaningful regret bounds
for two specific classes of proper scoring rules including the standard likelihood loss and the pairwise likelihood and exponential losses (Proposition~\ref{pro:bregman-bd}).

\vspace{.1in}\textbf{Related work.} There is an extensive literature on multi-category classification including and beyond the special case of binary classification.
We discuss directly related work to ours, in addition to
those mentioned above. An inverse mapping from a generalized entropy to a proper scoring rule can be seen in the canonical representation of proper scoring rules \citep{savage1971,gneiting2007}.
\tzq{This and related representations are extensively used in the design and study of composite binary losses \citep{buja2005, reid2011}
and composite multi-class losses \citep{williamson2016}.}

Recently, an inverse mapping is constructed by \cite{duchi2018} from a generalized entropy to a convex loss with actions in $\bbR^m$,
hence different from the canonical representation of proper scoring rules.
\tzq{Our construction of losses is in a similar direction as \cite{duchi2018}, but operates explicitly through a dissimilarity function $f$.
Our approach is applicable to handling both proper scoring rules and hinge-like losses and leads to interesting new findings.}
In particular, using an additive function $f$ provides a convenient, general extension of two-class proper scoring rules to
multi-class settings. By comparison, using an additive generalized entropy does not seem to achieve a similar effect.
Moreover, our inverse mapping in terms of $f$ are applied to discover new hinge-like losses, by first identifying a hinge-like loss on the probability simplex and
then constructing a convex extension. The hinge-like losses in \cite{lee2004} and \cite{duchi2018} are also such convex extensions.
This point of view enriches our understanding of multi-class hinge-like losses.

\tzq{Our new regret bounds for proper scoring rules generalize two-class results in \cite{reid2011}, Section 7.1, to
multi-class settings, by carefully exploiting the Bregman representation for the regret of a proper scoring rule
together with a novel bound on the regret of the zero-one or cost-weighted classification loss (Lemmas~\ref{lem:misclass}--\ref{lem:misclass2}).
Such a generalization seems to be previously unnoticed \citep[cf.][]{williamson2016}.}
These classification regret bounds provide a quantitative guarantee on classification calibration, a qualitative property studied in
\cite{zhang2004a}, \cite{steinwart2007}, and \cite{tewari2007} among others.
Although two-class regret bounds can be obtained more generally for all margin-based losses similarly as for proper scoring rules \citep{zhang2004b,bartlett2006,scott2012},
such results seem to rely on simplification due to two classes.

\tzq{Our regret bounds for the new hinge-like losses are similar to those in \cite{duchi2018}. However, we also establish in general that
all losses with the same generalized entropy as the zero-one loss achieve a regret bound which ensures classification calibration.
Hence our result provides a more concrete sufficient condition for achieving classification calibration than in \cite{tewari2007}.
Recently, necessary and sufficient conditions are obtained for classification calibration while allowing prediction labels to differ from the original class labels,
for example, using a refrain option \citep{ramaswamy2016}. It is interesting to study possible extensions in that direction.}

\vspace{.1in}\textbf{Notation.}
Denote $\overline \bbR = \bbR \cup \{\infty\}$, $\bbR_+=\{b \in \bbR: b \ge 0\}$, and $\overline\bbR_+=\{b \in \overline\bbR: b \ge 0\}$.
%Define $0/0=0$ and $0\cdot \infty=0$.
For $m \ge 2$, denote $[m]$ as the set $ \{1,\ldots,m \}$, $1_m$ as the $m\times 1$ vector of all ones, $I_m$ as the $m\times m$ identity matrix,
and $\Delta_m$ as the probability simplex $\{q \in \bbR_+^m: 1_m^\T q =1\}$.
For $j\in[m]$, a basis vector $e_j \in \Delta_m$ is defined such that its $j$th element is 1 and the remaining elements are 0.
The indicator function $1\{\cdot\}$ is defined as 1 if the argument is true or 0 otherwise.

\section{Background}

We provide a selective review of basic concepts and results which are instrumental to our subsequent development. See \cite{grunwald2004}, \cite{buja2005},
\cite{gneiting2007}, and \cite{duchi2018} among others for more information.

\subsection{Losses, risks and entropies} \label{sec:loss}

Consider the population version of the multi-category classification problem. Let $X \in \mathcal X$ be a vector of observed covariates or features, but $Y \in [m]$ an unobserved class label,
where $(X,Y)$ are generated from some joint probability distribution which can be assumed to be known unless otherwise noted.
It is of interest to predict the value of $Y$ based on $X$ (i.e., assign $X$ to one of the $m$ classes).
The prediction can be performed using an action function $\alpha(x): \mathcal X \to \mathcal A$, and evaluated through a loss function $L(y,\alpha(x))$ when the true label of $x$ is $y$.
Typically, an action in the space $\mathcal A$ is a vector whose components, as probabilities or margins, measure the strengths of association with the $m$ classes.
The risk or expected loss of the action function $\alpha(x)$ is
\begin{align}
R_L ( \alpha ) = E( L(Y, \alpha(X))  = E \left\{ \sum_{j=1}^m \pi_j(X) L(j, \alpha(X)) \right\}, \label{eq:risk}
\end{align}
where $\pi_j(x) = P(Y=j| X=x)$, the conditional probability of class $j$ given covariates $x$,
and the second expectation is taken over the marginal distribution of $X$ only.

From another perspective, the preceding problem can also be formulated as a Bayesian experiment with $m$ probability distributions $(P_1,\ldots,P_m)$ on $\mathcal X$,
corresponding to the within-class distributions of covariates \citep{degroot1962}.
Denote by $p_j(x)$ the density function of $P_j$ with respect to a baseline measure $\mu$.
Given a label $Y=j$ (regarded as an $m$-valued parameter), the random variable $X$ is drawn from the distribution $P_j$.
Let $\pi^0 = (\pi^0_1,\ldots, \pi^0_m)^\T \in \Delta_m$ be the prior probabilities of $Y$, corresponding to the marginal class probabilities.
Then the posterior probabilities of $Y$ given $X=x$ are
$\pi_j (x) = \pi^0_j p_j(x) / \{\sum_{k=1}^m \pi^0_k p_k(x)$\}, the same as the conditional class probabilities given covariates mentioned above.
In this context, $R_L (\alpha)$ is also called the Bayes risk of $\alpha(x)$. By standard Bayes theory \citep[][Eq. (5)]{garcia2012},
the minimum Bayes risk, or even shortened as the Bayes risk, can be obtained as
\begin{align}
\inf_{\alpha(x) : \mathcal X \to \mathcal A} R_L (\alpha) = E \{ H_L (\pi(X)) \}, \label{eq:min-bayes}
\end{align}
where $\alpha(x) : \mathcal X \to \mathcal A$ can be any measurable function, $\pi(x) = (\pi_1(x), \ldots,\pi_m(x))^\T$,
and $H_L$ is a function defined on $\Delta_m$ such that for $\eta=(\eta_1,\ldots,\eta_m)^\T \in \Delta_m$,
\begin{align}
H_L (\eta) = \inf_{\gamma \in \mathcal A} \left\{ \sum_{j=1}^m \eta_j L(j,\gamma) \right\}, \label{eq:LtoH}
\end{align}
The function $H_L$, which is concave on $\Delta_m$, is called an uncertainty function \citep{degroot1962} or a generalized entropy associated with the loss $L$ \citep{grunwald2004}.

A subtle point is that minimization in (\ref{eq:min-bayes}) is over all measurable functions $\alpha(x): \mathcal X \to \mathcal A$, whereas
minimization in (\ref{eq:LtoH}) is over all elements $\gamma \in \mathcal A$.
The generalized entropy $H_L$ is merely a function on $\Delta_m$, induced by the loss $L(j,\gamma)$ on $[m]\times \mathcal A$, where the covariate vector $X$ is conditioned on (or lifted out).
Similarly, the risk of an action $\gamma \in \mathcal A$ is defined as $R_L (\eta, \gamma) = \sum_{j=1}^m \eta_j L(j , \gamma)$,
and the regret of the action is defined as
\begin{align}
B_L( \eta, \gamma) = R_L(\eta,\gamma) - H_L(\eta), \label{eq:regret}
\end{align}
where $H_L (\eta ) = \inf_{\gamma^\prime\in\mathcal A} R_L(\eta,\gamma^\prime)$ by (\ref{eq:LtoH}).
This simplification where the covariate vector $X$ is lifted out is often useful when studying losses and regrets.

\subsection{Scoring rules} \label{sec:scoring}

A scoring rule is a particular type of loss $L(j,q)$, where its action $q$ is a probability vector in $\Delta_m$, interpreted as the predicted class probabilities \citep{grunwald2004}.
Sometimes, the expected loss, $R_L (\eta,q) = \sum_{j=1}^m \eta_j L(j,q)$, is also referred to as a scoring rule, for measuring the discrepancy between
underlying and predicted probability vectors, $\eta$ and $q$ \citep{gneiting2007}.

A scoring rule $L(j,q)$ is said to be proper if $H_L (\eta) = R_L(\eta,\eta)$, i.e.,
\begin{align*}
R_L ( \eta, \eta) \le R_L(\eta, q) , \quad \eta,q \in \Delta_m.
\end{align*}
The rule is strictly proper if the inequality is strict for $q\not=\eta$.
Hence for a proper scoring rule, the expected loss $R_L(\eta, q)$ is minimized over $q\in \Delta_m$ when $q= \eta$, the predicted probability vector coincides with the underlying probability vector.
This condition is typically required for establishing large-sample consistency of (conditional) probability estimators \citep[e.g.,][]{zhang2004b,buja2005}.

As shown in \cite{savage1971} and \cite{gneiting2007}, a proper scoring rule $L(j,q)$ in general admits the following representation:
\begin{align}
R_L (\eta, q) = H_L (q) - (q-\eta)^\T \partial H_L (q), \quad \eta,q \in\Delta_m, \label{eq:savage}
\end{align}
where $-\partial H_L$ is a sub-gradient of the convex function $-H_L$ on $\bbR^m$. Note that the generalized entropy $H_L$ is evaluated at $q$, not $\eta$, in (\ref{eq:savage}).
Then the regret in (\ref{eq:regret}) becomes
\begin{align}
B_L (\eta, q) =  H_L (q) - H_L(\eta) - (q-\eta)^\T \partial H_L (q), \label{eq:bregman}
\end{align}
which is the Bregman divergence from $q$ to $\eta$ associated with the convex function $-H_L$.

An important example of proper scoring rules is the logarithmic scoring rule \citep{good1952},
$L (j, q) = - \log q_j$.
The corresponding expected loss is $R_L(\eta,q)=-\sum_{j=1}^m \eta_j \log q_j$, which is, up to scaling, the negative expected log-likelihood
of the predicted probability vector $q$ with the underlying probability vector $\eta$ for multinomial data.
The generalized entropy is $H_L(\eta)=-\sum_{j=1}^m \eta_j \log \eta_j$, the negative Shannon entropy.
The regret is $B_L(\eta,q) = \sum_{j=1}^m \eta_j \log (\eta_j / q_j)$, the Kullback--Liebler divergence.

\subsection{Classification losses} \label{sec:misclass}

Consider the zero-one loss, formally defined as
\begin{align*}
L^{\text{zo}} (j, \gamma) = 1 \{j \not= \argmax_{k\in [m]}\gamma_k \}, \quad j\in[m], \gamma\in\bbR^m,
\end{align*}
where, if not unique, $\argmax_{k\in [m]}\gamma_k$ can be fixed as the index of any maximum component of $\gamma$.
As mentioned below (\ref{eq:test-risk}), $L^{\text{zo}}$ is typically used to evaluate performance,
but not as the loss $L$ for training, and the action $\gamma$ can be transformed from the action of $L$.
Nevertheless, the generalized entropy defined by (\ref{eq:LtoH}) with $L=L^{\text{zo}}$ is
\begin{align}
H^{\text{zo}} ( \eta ) = 1- \max_{k\in[m]} \eta_k, \quad \eta\in\Delta_m .  \label{eq:zo-H}
\end{align}
This function is concave and continuous, but not everywhere differentiable.

In practice, there can be different costs of misclassification, depending on which classes are involved. For example,
the cost of classifying a cancerous tumor as benign can be greater than in the other direction.
Let $C = (c_{jk})_{j,k \in [m]}$ be a cost matrix, where
$c_{jk} \ge 0$ indicates the cost of classifying class $j$ as class $k$. For each $j \in [m]$, assume that $c_{jj}=0$ and $c_{jk}>0$ for some $k\not=j$.
Consider the cost-weighted classification loss
\begin{align*}
L^{\text{cw}} (j, \gamma) = c_{jk} \mbox{ if } k = \argmax_{l\in [m]} \gamma_l, \quad j\in[m], \gamma\in\bbR^m.
\end{align*}
As shown in \cite{duchi2018}, the generalized entropy defined by (\ref{eq:LtoH}) with $L=L^{\text{cw}}$ is
\begin{align}
H ^{\text{cw}} ( \eta ) = \min_{k\in[m]} \eta^\T C_k, \quad \eta\in\Delta_m, \label{eq:cw-H}
\end{align}
where $C=(C_1, \ldots, C_m)$ is the column representation of $C$.
The standard zero-one loss corresponds to the special choice $C= 1_m 1_m^\T - I_m$.

An intermediate case is the class weighted classification loss,
\begin{align*}
L^{\text{cw0}} (j, \gamma) = c_{j0} 1\{ j \not= \argmax_{k\in [m]}\gamma_k \}, \quad j\in[m], \gamma\in\bbR^m,
\end{align*}
where $c_{j0} >0$ is the cost associated with misclassification of class $j$.
This loss is more general than the standard zero-one loss $L^{\text{zo}}$,
although a special case of the cost-weighted loss $L^{\text{cw}}$ with $C = C_0 1_m^\T - \diag(C_0)$, where $C_0 = (c_{10}, \ldots, c_{m0})^\T$.
The generalized entropy associated with $L^{\text{cw0}} $ is $H^{\text{cw0}} (\eta) =\eta^\T C_0 - \max_{k\in[m]} \,\eta_k c_{k0}$.

\subsection{Entropies and divergences} \label{sec:div}

In \citeauthor{degroot1962}'s \citeyearpar{degroot1962} theory, any concave function $H$ on $\Delta_m$ can be used as an uncertainty function.
The information of $X$ about label (``parameter") $Y$ is defined as the reduction of uncertainty (or entropy) from the prior to the posterior:
\begin{align*}
I_H (X; \pi^0) = H( \pi^0) - E \{ H (\pi(X)) \} ,
\end{align*}
which is nonnegative by the concavity of $H$. The information $I_H(X;\pi^0)$ is closely related to the
$f$-divergence between the multiple distributions $(P_1,\ldots, P_m)$, which is a generalization of the $f$-divergence between two distributions \citep{ali1966,csiszar1967}
Heuristically, the more dissimilar $(P_1,\ldots,P_m)$ are from each other, the more information about $Y$ is obtained after observing $X$.

For a convex function $f$ on $\overline\bbR_+^{m-1}$ with $f(1_{m-1})=0$, the $f$-divergence between $(P_1,\ldots, P_{m-1})$ and $ P_m$ with  densities
$(p_1,\ldots,p_{m-1})$ and $p_m$ is
\begin{align*}
D_f ( P_{1:(m-1)} \| P_m ) = \frac{1}{m} \int f\left( \frac{p_1(x)}{p_m(x)}, \ldots, \frac{p_{m-1}(x)}{p_m(x)} \right) p_m(x) \,\dif \mu(x),
\end{align*}
which is nonnegative by the convexity of $f$.
Compared with the standard definition of multi-way $f$-divergences \citep{gyorfi1978,duchi2018}, our definition above involves a rescaling factor $m^{-1}$,
for notational simplicity in the later discussion; otherwise, for example, rescaling would be needed in Eqs.~(\ref{eq:ftoH}) and (\ref{eq:Htof}).

There is a one-to-one correspondence between the statistical information $I_H (X; \pi^0) $ and multi-way $f$-divergences, as discussed in \cite{garcia2012}.
For any prior probability $\pi^0 \in \Delta_m$ and probability distributions $(P_1,\ldots,P_m)$,
if a convex function $f$ on $\overline\bbR_+^{m-1}$ with $f(1_{m-1})=0$ and a concave function $H$ on $\Delta_m$ are related such that for $\eta=(\eta_1,\ldots,\eta_m)^\T \in \Delta_n$,
\begin{align}
H( \eta ) = -\frac{\eta_m}{m \pi^0_m} f\left( \frac{\pi^0_m}{\pi^0_1} \frac{\eta_1}{\eta_m}, \ldots, \frac{\pi^0_m}{\pi^0_{m-1}} \frac{\eta_{m-1}}{\eta_m} \right), \label{eq:f-H}
\end{align}
then $I_H( X; \pi^0) = D_f ( P_{1:(m-1)} \| P_m )$ or, because $H(\pi^0) = -f(1_{m-1}) =0$ here,
\begin{align}
-E \{ H (\pi(X)) \} = D_f ( P_{1:(m-1)} \| P_m ) , \label{eq:f-H2}
\end{align}
where the expectation is taken over $X \sim \sum_{j=1}^m \pi^0_j P_j$.

\section{Construction of losses} \label{sec:construct-loss}

In practice, a learning method for classification involves minimization of (\ref{eq:emp-risk}), an empirical version of the risk (\ref{eq:risk}) based on training data,
with specific choices of a loss function $L(y, \alpha)$ and a potentially rich family of action functions $\alpha(x)$.
As suggested in Section~\ref{sec:loss}, we study construction of the loss $L(y,\alpha)$ as a function of a label $y$ and
a freely-varying action $\alpha$, with the dependency on covariates (or features) lifted out.
As a result, we not only derive new general classes of losses, but also improve understanding of various existing losses as shown in Sections~\ref{sec:scoring-examples}--\ref{sec:hinge-like}.
Nevertheless, the interplay between losses and function classes remains important, but challenging to study, for further research.

Equation (\ref{eq:LtoH}) is a mapping from a loss $L$ to a generalized entropy $H_L$, which is in general many-to-one (i.e., different losses can lead to the same generalized entropy).
\cite{duchi2018} constructed an inverse mapping  from a generalized entropy to a convex loss.
For a closed, concave function $H$ on $\Delta_m$, define a loss with action space $\mathcal A = \bbR^m$ such that for $\gamma=(\gamma_1,\ldots,\gamma_m)^\T \in \bbR^m$,
\begin{align}
L_H (j, \gamma) = -\gamma_j + (-H)^*(\gamma), \label{eq:HtoL}
\end{align}
where $(-H)^*(\gamma) = \sup_{\eta \in \Delta_m} \{\gamma^\T \eta + H(\eta) \}$, the conjugate of $-H$.
Then $L_H(j,\gamma)$ is convex in $\gamma$ and (\ref{eq:LtoH}) is satisfied with \tzq{$H_{L_H} = H$}, by \cite{duchi2018}, Proposition 3.
Hence a convex loss is obtained for a concave function on $\Delta_m$ to be the generalized entropy.
Note that the loss $L_H$ is over-parameterized, because
$ (-H)^*(\gamma - b 1_m) = -b + (-H)^*(\gamma)$
and hence $L_H (j, \gamma-b 1_m)= L_H (j, \gamma)$ for any constant $b \in \bbR$.

We derive a new mapping from generalized entropies to convex losses, by working with perspective-like functions related to multi-distribution $f$-divergences.
First, there exists a one-to-one correspondence between concave functions $H$ on $\Delta_m$ and convex functions $f$ on $\overline\bbR_+^{m-1}$.
%Without loss of generality, assume $H(1_m/m)=0$ and $f(1_{m-1}) =0$; otherwise, redefine $H$ by subtracting $H(1_m/m)$ or $f$ by subtracting $f(1_{m-1})$.
For a convex function $f$ on $\overline\bbR_+^{m-1}$, define a function on $\Delta_m$:
\begin{align}
H_f (\eta) = - \eta_m f \left(\frac{\eta_1}{\eta_m}, \ldots, \frac{\eta_{m-1}}{\eta_m} \right) . \label{eq:ftoH}
\end{align}
Conversely, for a concave function $H$ on $\Delta_m$,  define a function on $\overline\bbR_+^{m-1}$:
\begin{align}
f_H (t) = - t_\bullet \,H\left( \frac{t_1}{t_\bullet} , \ldots, \frac{t_{m-1}}{t_\bullet}, \frac{1}{t_\bullet} \right), \label{eq:Htof}
\end{align}
where $t= (t_1,\ldots,t_{m-1})^\T$ and $t_\bullet=1+\sum_{j=1}^{m-1} t_j$. The mappings $H_f$ and $f_H$ are of a similar form to perspective functions associated with $f$ and $H$ respectively,
although neither fits the standard definition of perspective functions \citep{boyd2004}.

\begin{lem}[\citealt*{garcia2012}] \label{lem:f-H}
For a convex function $f$ on $\overline\bbR_+^{m-1}$, the function $H_f$ defined by (\ref{eq:ftoH}) is concave on $\Delta_m$ such that (\ref{eq:Htof}) is satisfied with \tzq{$f_{H_f} = f$}.
Conversely, for a concave function $H$ on $\Delta_m$, the function $f_H$ defined by (\ref{eq:Htof}) is convex on $\overline\bbR_+^{m-1}$ such that (\ref{eq:ftoH}) is satisfied with \tzq{$H_{f_H} = H$}.
Moreover, %the restrictions $H(1_m/m)=0$ and $f(1_{m-1}) =0$ are preserved under (\ref{eq:ftoH}) and (\ref{eq:Htof}).
it is preserved that $H(1_m/m)= - m^{-1} f(1_{m-1})$ under (\ref{eq:ftoH}) and (\ref{eq:Htof}).
\end{lem}

\begin{rem} \label{rem:f-H}
Equations (\ref{eq:ftoH}) and (\ref{eq:Htof}) can be obtained as a special case of (\ref{eq:f-H}) with $\pi^0 = 1_m /m$, from \cite{garcia2012}.
As mentioned in Section~\ref{sec:div}, (\ref{eq:f-H}) is originally determined such that identity (\ref{eq:f-H2}) holds for linking the expected entropy and multi-way $f$-divergences,
which are, by definition, concerned with the covariates and within-class distributions.
Nevertheless, our subsequent development is technically independent of this connection, because covariates are lifted out in our study.
In other words, we merely use (\ref{eq:ftoH}) and (\ref{eq:Htof}) as convenient mappings between $H$ and $f$.
The usual restriction $f(1_{m-1})=0$ used in $f$-divergences does not need to be imposed.
\end{rem}

Our first main result shows
a mapping from a convex function $f$ to a convex loss $L$ such that the concave function $H_f$ is the generalized entropy associated with $L$.

\begin{pro} \label{pro:f-loss}
For a closed, convex function $f$ on $\overline\bbR_+^{m-1}$, define a loss with action space $\mathcal A = \dom(f^*)$ such that for $s  = (s_1,\ldots,s_{m-1})^\T \in \dom(f^*)$,
\begin{align}
L_f(j, s) = \left\{ \begin{array}{cl}
- s_j, & j \in [m-1],\\
f^*(s), & j=m,
\end{array} \right. \label{eq:ftoL}
\end{align}
where $f^*(s) = \sup_{t \in \overline\bbR_+^{m-1}} \{ s^\T t - f(t)\}$ and $\dom(f^*)=\{s \in \bbR_+^{m-1}: f^*(s)<\infty\}$.
Then $L_f(j,s)$ is convex in $s$. Moreover, the concave function $H_f$ defined by (\ref{eq:ftoH})
is the generalized entropy associated with $L_f$, that is, (\ref{eq:LtoH}) is satisfied with \tzq{$H_{L_f} = H_f$}.
\end{pro}

\begin{prf}
For $\eta\in\Delta_m$ and $s\in\dom(f^*)$, the definition of $L_f$ implies that $\sum_{j=1}^m \eta_j L_f(j,s) = - \sum_{j=1}^{m-1} \eta_j s_j + \eta_m f^*(s)$.
Hence
\begin{align*}
& \inf_{s \in \mathcal A}\left\{\sum_{j=1}^m \eta_j L_f(j,s) \right\} = - \sup_{s \in \mathcal A} \left\{ \sum_{j=1}^{m-1} \eta_j s_j - \eta_m f^*(s) \right\} \\
& =- \eta_m \sup_{s \in \dom(f^*)} \left\{ \sum_{j=1}^{m-1} \frac{\eta_j}{\eta_m} s_j - f^*(s) \right\}
 = -\eta_m f \left( \frac{\eta_1}{\eta_m },\ldots, \frac{\eta_{m-1}}{\eta_m} \right) =H_f(\eta),
\end{align*}
where the second last equality holds by Fenchel's conjugacy relationship.
\end{prf}

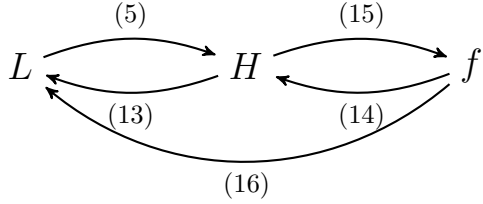
\begin{figure} %[!h]
\centering
\begin{tikzpicture}[->,>=stealth',shorten >=1pt,auto,node distance=3cm,
                    thick,main node/.style={font=\sffamily\Large\bfseries}]\label{sm1}

  \node[main node] (1) {$H$};
  \node[main node] (2) [left of=1] {$L$};
  \node[main node] (3) [right of=1] {$f$};
 \path[every node/.style={font=\sffamily\small}]
  (1) edge [bend left=20] node  {$(13)$} (2)
      edge [bend left=20] node  {$(15)$} (3)  
  (2) edge [bend left=20] node  {$(5)$} (1)   
  (3) edge [bend left=20] node  {$(14)$} (1) 
  (3) edge [bend left = 40] node {$(16)$} (2);
\end{tikzpicture} \vspace{-.2in}
\caption{Relations between loss $L$, generalized entropy $H$, and dissimilarity function $f$.}\label{fig:relation}
\end{figure}

Compared with (\ref{eq:HtoL}), Eq.~(\ref{eq:ftoL}) together with (\ref{eq:Htof}) presents an alternative approach for determining a convex loss $L$
from a generalized entropy $H$ through a dissimilarity function $f$.
See Figure~\ref{fig:relation} which illustrates various relationships discussed.
For ease of interpretation, a convex function $f$ on $\overline\bbR_+^{m-1}$ can be called a dissimilarity function, similarly as a concave function $H$ on $\Delta_m$ can be a generalized entropy.

In spite of the one-to-one correspondence between entropy and dissimilarity functions $H$ and $f$
by (\ref{eq:ftoH}) and (\ref{eq:Htof}), we stress that the new loss (\ref{eq:ftoL})
is in general distinct from (\ref{eq:HtoL}). An immediate difference, which is further discussed in Section~\ref{sec:hinge-like}, is that the action space for loss (\ref{eq:ftoL}),
$\dom(f^*)$, can be a strict subset of $\bbR^{m-1}$, whereas the action space for loss (\ref{eq:HtoL}) is either $\bbR^{m}$ with over-parametrization as noted above or $\bbR^{m-1}$
with, for example, $\gamma_m=0$ fixed to remove over-parametrization.
Moreover, loss (\ref{eq:ftoL}) can also be used to derive a new class of closed-form losses based on arbitrary convex functions $f$ as shown in the following result
and, as discussed in Section~\ref{sec:hinge-like}, to find novel multi-class, hinge-like losses related to the zero-one or cost-weighted classification loss.

\begin{pro} \label{pro:f-loss2}
For a closed, convex function $f$ on $\overline\bbR_+^{m-1}$, define a loss with action space $\mathcal A = \overline\bbR_+^{m-1}$ such that for $u  = (u_1,\ldots,u_{m-1})^\T \in \overline\bbR_+^{m-1}$,
\begin{align}
L_{f2}(j, u) = \left\{ \begin{array}{cl}
- \partial_j f(u), & j \in [m-1],\\
u^\T \partial f(u) - f(u), & j=m,
\end{array} \right. \label{eq:ftoL2}
\end{align}
where $\partial f = (\partial_1 f, \ldots, \partial_{m-1} f)^\T$ is a sub-gradient of $f$, arbitrarily fixed (if needed).
Then the concave function $H_f$ defined by (\ref{eq:ftoH})
is the generalized entropy associated with $L_{f2}$, that is, (\ref{eq:LtoH}) is satisfied with \tzq{$H_{L_{f2}} = H_f$}.
\end{pro}

\noindent\textbf{Proof (outline).}\;
A basic idea is to use the parametrization $s = \partial f(u)$ and Fenchel's conjugacy property $ f^*(s) = u^\T s -f(u)$,
and then obtain the loss $L_{f2}$ from $L_f$ in Proposition~\ref{pro:f-loss}.
This argument gives a one-sided inequality for the desired equality (\ref{eq:LtoH}). A complete proof is provided in the Supplement.
\hfill$\blacksquare$\\[2mm]

Compared with loss (\ref{eq:ftoL}), the preceding loss (\ref{eq:ftoL2}) is of a closed form without involving the conjugate $f^*$,
which can be nontrivial to calculate. On the other hand, loss (\ref{eq:ftoL2}) may not be convex in its action $u$.
Nevertheless, it is often possible to choose a link function, for example, $u^h= (u^h_1, \ldots, u^h_{m-1})^\T$ with $u^h_j = \exp(h_j)$ such that
$L_{f2}(j, u^h)$ becomes convex in $(h_1,\ldots,h_{m-1})$.
This link can be easily identified as the multinomial logistic link after reparameterizing $(u_1, \ldots, u_{m-1})$ as probability ratios below.

The following result shows that a reparametrization of loss (\ref{eq:ftoL2}) with actions defined as probability vectors in $\Delta_m$
automatically yields a proper scoring rule. See Section~\ref{sec:scoring} for the related background on scoring rules.
Together with the relationship between $f$ and $H$ by (\ref{eq:ftoH}) and (\ref{eq:Htof}), our construction gives a mapping from
a dissimilarity function $f$ or equivalently a generalized entropy $H$ to a proper scoring rule.

\begin{pro} \label{pro:f-loss3}
For a closed, convex function $f$ on $\overline\bbR_+^{m-1}$, define a loss with action space $\mathcal A = \Delta_m$ such that for $q  = (q_1,\ldots,q_m)^\T \in \Delta_m$,
\begin{align}
 L_{f3}(j, q) = \left\{ \begin{array}{cl}
- \partial_j f(u^q), & j \in [m-1],\\
u^{q\T} \partial f(u^q) - f(u^q), & j=m,
\end{array} \right. \label{eq:ftoL3}
\end{align}
where $u^q = (q_1/q_m, \ldots, q_{m-1}/ q_m)^\T$.
Then $L_{f3}$ is a proper scoring rule, with $H_f$ defined  by (\ref{eq:ftoH}) as  the generalized entropy, satisfying
\begin{align*}
\inf_{q \in \Delta_m} \left\{ \sum_{j=1}^m \eta_j L_{f3}(j,q) \right\} = H_f(\eta) = \sum_{j=1}^m \eta_j L_{f3}(j,\eta)  .
\end{align*}
\end{pro}

\begin{prf}
The generalized entropy from $L_{f3}$ is $H_f$, due to Proposition~\ref{pro:f-loss2} and the one-to-one mapping $u^q = (q_1/q_m, \ldots, q_{m-1}/ q_m)^\T$. Then direct calculation
shows that
\begin{align*}
& \sum_{j=1}^m \eta_j L_{f3}(j,\eta) = -\sum_{j=1}^{m-1} \eta_j \partial_j f( u^\eta)  +  \eta_m \left\{ \sum_{j=1}^{m-1} \frac{\eta_j}{\eta_m} \partial_j f(u^\eta) -f(u^\eta) \right\} \\
& = -\eta_m f(u^\eta) =  H_f(\eta) = \inf_{q \in \Delta_m} \left\{ \sum_{j=1}^m \eta_j L_{f3}(j,q) \right\}.
\end{align*}
Hence $L_{f3}$ is a proper scoring rule.
\end{prf}

For completeness, the expected loss associated with $L_{f_3}$ can be shown to satisfy the canonical representation (\ref{eq:savage}) with $H_f$ defined  by (\ref{eq:ftoH}):
\begin{align}
& \sum_{j=1}^m \eta_j L_{f3}(j,q) = -\sum_{j=1}^{m-1} \eta_j \partial_j f( u^q)  +  \eta_m \left\{ \sum_{j=1}^{m-1} \frac{q_j}{q_m} \partial_j f(u^q) -f(u^q) \right\} \nonumber \\
& = H_f( q) - \sum_{j=1}^m  (q_j - \eta_j) \partial_j H_f( q) ,   \label{eq:savage2}
\end{align}
where $-\partial H_f = (-\partial_1 H_f, \ldots, -\partial_m H_f)^\T$ is the sub-gradient of $-H_f$. See the Supplement for a proof.
Conversely, the loss $L_{f3}$ can also be obtained by calculating the canonical representation (\ref{eq:savage}) for the concave function $H_f$ in (\ref{eq:ftoH})
and then taking $\eta$ to be a basis vector, $e_1,\ldots,e_m$, one by one in the resulting expression, which is on the left of the second equality in (\ref{eq:savage2}).
%Similarly, taking $\eta$ to be unit vectors in the canonical representation (\ref{eq:savage}) yields a proper scoring rule directly from a generalized entropy $H$ as
%\begin{align}
%L_{H2} (j, q) = H (q) - q^\T \partial H(q) + \partial_j H(q), \quad j \in [m], q \in \Delta_m . \label{eq:savage3}
%\end{align}
%Note that the term $H (q) - q^\T \partial H(q)$ cannot in general be dropped, even though independent of $j$.
Moreover, by the necessity of the representation (\ref{eq:savage}), we see that $L_{f_3}$ in (\ref{eq:ftoL3}) is the only proper scoring rule with the generalized entropy $H_f$.

\begin{cor} \label{cor:f-loss3}
For a closed, convex function $f$ on $\overline\bbR_+^{m-1}$,
any proper scoring rule with $H_f$ in (\ref{eq:ftoH}) as the generalized entropy can be expressed as $L_{f3}$ in (\ref{eq:ftoL3}), up to possible choices of sub-gradients of $f$, $\{\partial_j f: j\in [m-1]\}$.
\end{cor}

While the preceding use of the canonical representation (\ref{eq:savage}) seems straightforward, our development from Propositions~\ref{pro:f-loss} to \ref{pro:f-loss3} remains worthwhile.
The proper scoring rule $L_{f3}$ in (\ref{eq:ftoL3}) is of simple form, depending explicitly on a dissimilarity function $f$.
Moreover, as shown in Section~\ref{sec:hinge-like},
Proposition~\ref{pro:f-loss} can be further exploited to derive new convex losses which are related to classification losses but are not proper scoring rules.

\section{Examples of proper scoring rules} \label{sec:scoring-examples}

We examine various examples of multi-class proper scoring rules, obtained from Proposition~\ref{pro:f-loss3}.
In particular, it is of interest to study how commonly used two-class losses can be extended to multi-class ones.
These examples lead to new multi-class proper scoring rules and shed new light on existing ones.
See Section~\ref{sec:hinge-like} for a discussion of multi-class hinge-like losses related to zero-one classification losses, derived using Proposition~\ref{pro:f-loss}.

\vspace{.1in}\textbf{Two-class losses.}
For two-class classification ($m=2$) and a univariate convex function $f_0$ on $\overline \bbR_+$,
the proper scoring rule (\ref{eq:ftoL3}) in Proposition~\ref{pro:f-loss3} reduces to
\begin{align}
L_{f_0} (j, q) = - \one_1(j) \partial f_0 (u^q ) + \one_2(j) \big\{ u^q \partial f_0 (u^q) - f(u^q) \big\}, \quad j=1,2, \label{eq:2class-f-loss}
\end{align}
where $u^q = q_1 / q_2$ for $q = (q_1,q_2)^\T \in \Delta_2$, $\partial f_0$ denotes a sub-gradient of $f_0$, and $\one_k(j)$ is an indicator defined as 1 if $j=k$ or 0 otherwise.
%In fact, the representation (\ref{eq:2class-f-loss}) is broad enough to encompass all (two-class) proper scoring rules.
For a twice-differentiable function $f_0$, the gradient of the loss (\ref{eq:2class-f-loss}) can be directly calculated as %(Schervish 1989; Buja et al.~2005),
\begin{align}
\frac{\dif}{\dif q_1} L_{f_0}(j, q) = -\{\one_1(j) - q_1\} w(q_1), \label{eq:2class-grad}
\end{align}
where $\dif /\dif q_1$ denotes a derivative taken with respect to $q_1$ with $q_2=1-q_1$, and
the weight function $w(q_1) =f_0^\dprime(u^q) / q_2^3$ with $f_0^\dprime$ the second derivative of $f_0$. %equivalently $w(q_1) = -\frac{\dif^2} {\dif q_1^2} H_{f_0} (q)$
\tzq{From (\ref{eq:2class-grad}), $L_{f_0}(j,q)$ can be put into an integral representation in terms of $w(\cdot)$ and the cost-weighted binary classification loss \citep{schervish1989, buja2005, reid2011}.
The formula (\ref{eq:2class-f-loss}) in terms of $f_0$ differs from the integral representation or the canonical representation (\ref{eq:savage}),
even though they can be transformed into each other.}

For concreteness, consider the following examples of two-class losses:
\begin{itemize}
\item {\it Likelihood loss}:
$L_\ell(j, q) = -\log q_j$ with $f_0 = t\log t - (1+t) \log (1+t)$,

\item {\it Exponential loss}: $L_e (j, q) = \one_1(j) \sqrt{q_2/q_1} + \one_2(j) \sqrt{q_1/q_2}$ with  $f_0 = (\sqrt t -1)^2$,

\item {\it Calibration loss}:  $L_c (j, q) =\{ \one_1(j) (q_2/q_1) + \one_2(j) \log(q_1/q_2)\}/2 $ with $f_0 = - (\log t)/2$,
\end{itemize}
where all the expressions for $L(j,q)$ are up to additive constants in $q$. See Supplement Table~S1 for further information.
While the likelihood loss is tied to maximum likelihood estimation, the exponential loss is associated with boosting algorithms \citep{friedman2000,schapire2012}.
The calibration loss is studied in \cite{tan2020} for logistic regression, where the fitted probabilities are  used for inverse probability weighting.
\tzq{See the Supplement for a discussion on convexity of these losses with a logistic link.}

\begin{rem} \label{rem:f-gan}
The loss (\ref{eq:2class-f-loss}) was also derived in \cite{tan2019} for training a discriminator in generative adversarial learning \citep{goodfellow2014,nowozin2016}.
In that context, the loss for training a generator is, in a nonparametric limit, the negative Bayes risk from discrimination or the $f_0$-divergence by relationship (\ref{eq:f-H2}) with $\pi^0=1_m/m$,
\begin{align*}
D_{f_0} (P_1 \| P_2 ) = - \inf_{q(x): \mathcal X \to \Delta_2} E \big\{ L_{f_0}(Y, q(X)) \big\},
\end{align*}
where $P_1$ is the data distribution represented by training data and $P_2$ is the model distribution represented by simulated data from the generator.
Hence the generator can be trained to minimize various $f_0$-divergences, including forward and reverse Kullback--Liebler and Hellinger divergences. See Supplement Table S1 in \cite{tan2019}.
\end{rem}

\vspace{.1in}\textbf{Multi-class pairwise losses.}
There can be numerous choices for extending a two-class loss (\ref{eq:2class-f-loss}) to multi-class ones, just as a univariate convex function $f_0$
can be extended in multiple ways to multivariate ones.
A simple approach is to use an additive extension,
$f(u_1, \ldots, u_{m-1} ) = \sum_{k=1}^{m-1} f_0 ( u_k)$. The corresponding loss (\ref{eq:ftoL3}) is then
\begin{align}
 L_{f_0}^{\text{pw,a}}(j, q) =\sum_{k=1}^{m-1} \left[- \one_k(j) \partial  f_0( \frac{q_k}{q_m} ) +
 \one_m(j) \left\{ \frac{q_k}{q_m} \partial f_0( \frac{q_k}{q_m}) - f_0( \frac{q_k}{q_m}) \right\} \right],   \label{eq:pairwise-a}
\end{align}
Equivalently, the loss (\ref{eq:pairwise-a}) can be obtained by applying the two-class loss (\ref{eq:2class-f-loss}) to a pair of classes, $k$ and $m$, and
summing up such pairwise losses for $k \in [m-1]$. In this sense, the loss (\ref{eq:pairwise-a}) can be interpreted as performing multi-class classification
via pairwise comparison of each class $k \in[m-1]$ with class $m$.

The preceding loss (\ref{eq:pairwise-a}) is asymmetric with class $m$ compared with the remaining classes $k\in[m-1]$. A symmetrized version can be obtained
by varying the choice of a base class and summing up the resulting losses as
\begin{align}
 & L_{f_0}^{\text{pw,s}}(j, q) = \sum_{l,k \in [m], k\not= l} \left[- \one_k(j) \partial  f_0( \frac{q_k}{q_l} ) +
 \one_l(j) \left\{ \frac{q_k}{q_l} \partial f_0( \frac{q_k}{q_l}) - f_0( \frac{q_k}{q_l}) \right\} \right] \nonumber \\
 & = \sum_{k\in [m], k\not=j} \left\{ - \partial  f_0( \frac{q_j}{q_k} ) + \frac{q_k}{q_j} \partial f_0( \frac{q_k}{q_j}) - f_0( \frac{q_k}{q_j}) \right\}.  \label{eq:pairwise-s}
\end{align}
See the Supplement for a proof.
The symmetrized loss (\ref{eq:pairwise-s}) can also be deduced from (\ref{eq:ftoL3}) with the choice
$f (u_1, \ldots, u_{m-1}) =  \sum_{l,k \in [m], k\not= l} u_l f_0( \frac{u_k}{u_l} )$, where $u_m \equiv 1$.
\tzq{In spite of the interpretation via pairwise comparison, our approach involves optimizing the loss (\ref{eq:pairwise-a}) or (\ref{eq:pairwise-s}) {\it jointly} over $q \in \Delta_m$ using all $m$ labels,
and hence differs from the usual one-against-all or all-pairs approach, which performs binary classification with 2 reduced labels {\it separately} for multiple times.
Further comparison of these approaches can be studied in future work.}

Consider a multinomial logistic link $q^h = (q_1^h, \ldots, q_m^h)^\T$, where $h=(h_1,\ldots,h_m)^\T$ and
\begin{align}
q_j^h = \frac{\exp( h_j) }{\sum_{k=1}^m \exp( h_k)}, \quad j \in [m]. \label{eq:multinomial-link}
\end{align}
The link is a natural extension of the logistic link, because log ratios between $(q_1,\ldots,q_m)$ are related to contrasts between $(h_1,\ldots,h_m)$.
To remove over-parametrization, a restriction is often imposed such as $h_m \equiv 0$ or $\sum_{k=1}^m h_k \equiv 0$.
By the additive construction, the composite losses obtained from (\ref{eq:pairwise-a}) and (\ref{eq:pairwise-s}) can be easily shown to be convex in $h$
whenever the two-class loss (\ref{eq:2class-f-loss}) with a logistic link  $q_1^{h_0}/q_2^{h_0} = \exp(h_0)$ is convex in $h_0$.

For the two-class likelihood, exponential, and calibration losses above, the pairwise extensions (\ref{eq:pairwise-s}) can be calculated as follows:
\begin{itemize}
\item {\it Pairwise likelihood loss}:
$L_\ell^{\text{pw,s}}(j, q) = 2 \sum_{k\in[m], k\not=j} \log(1+ \frac{q_k}{q_j})$,

\item {\it Pairwise exponential loss}: $L_e^{\text{pw,s}} (j, q) = 2 \sum_{k\in[m], k\not=j} \sqrt{\frac{q_k}{q_j}}$,

\item {\it Multi-class calibration loss}:  $L_c^{\text{pw,s}} (j, q) = \sum_{k\in[m], k\not=j} \{\log(\frac{q_k}{q_j}) + \frac{q_k}{q_j} \} /2$,
\end{itemize}
where additive constants in $q$ are dropped for simplicity. See Supplement Table~S1 for the expressions of the corresponding $f$, $H$, and gradients.
By convexity of the associated two-class composite losses \citep{buja2005}, we see that
with the multinomial logistic link (\ref{eq:multinomial-link}), the three composite losses, $L_\ell^{\text{pw,s}}(j, q^h)$,
$L_e^{\text{pw,s}}(j, q^h)$, and $L_c^{\text{pw,s}}(j, q^h) $, are all convex in $h$.
In particular, the  pairwise exponential composite loss is
\begin{align*}
L_e^{\text{pw,s}} (j, q^h) = 2 \sum_{k\in[m], k\not=j} \me^{ (h_k-h_j)/2 } ,
\end{align*}
which is associated with multi-class boosting algorithms AdaBoost.M2 \citep{freund1997} or AdaBoost.MR \citep{schapire1999}.
See \cite{mukherjee2013} for further study.
The pairwise likelihood and calibration losses appear to be new.
The pairwise likelihood loss, with $m \ge 3$, differs from the standard likelihood loss based on multinomial data, which will be discussed later.
The multi-class calibration loss can be useful for inverse probability weighting with multi-valued treatments.

\vspace{.1in}\textbf{Multi-class simultaneous losses.}
Apparently, there exist various multi-class proper scoring rules, which cannot be expressed as pairwise losses (\ref{eq:pairwise-a}) or (\ref{eq:pairwise-s})
and hence will be referred to as simultaneous losses.
A notable example as mentioned above is the standard likelihood loss (or the logarithmic scoring rule) for multinomial data,
$L(j,q) = -\log q_j$. In fact, a large class of multi-class simultaneous losses can be defined with the generalized entropy in the form
\begin{align*}
H_\beta (q) = \left\{
\begin{array}{cl}
\| q \|_\beta, & \mbox{if } \beta \in (0,1), \\
-  \|q \|_\beta , & \mbox{if } \beta \in (1,\infty),
\end{array} \right.
\end{align*}
where $ \| q \|_\beta = \{ \sum_{j=1}^m q_j^\beta \}^{1/\beta} $ is the $L_\beta$ norm.
The corresponding dissimilarity function is $f_\beta(t) = -\| \tilde t \|_\beta$ if $\beta \in (0,1)$ or $\|\tilde t\|_\beta$ if $\beta \in (1,\infty)$,
where $\tilde t= (t_1,\ldots, t_{m-1}, 1)^\T$.
The resulting scoring rule can be calculated by (\ref{eq:ftoL3}) as %for $j \in [m]$,
\begin{align}
L_\beta (j,q) = \left\{
\begin{array}{cl}
( q_j/ \|q\|_\beta )^{\beta-1}, & \mbox{if } \beta \in (0,1), \\
- ( q_j/ \|q\|_\beta )^{\beta-1} , & \mbox{if } \beta \in (1,\infty).
\end{array} \right.  \label{eq:simul-loss}
\end{align}
The case $\beta >1$ is called a pseudo-spherical score \citep{good1971,gneiting2007}.
The limiting case $\beta \to 1$ is also known to yield the logarithmic score, $L(j,q) = -\log q_j$, after suitable rescaling.
The case $\beta\in(0,1)$ seems previously unstudied. There are also two additional limiting cases as $\beta \to 0+$ or $\infty$. See Supplement Table~S1 for further details.

\begin{pro} \label{pro:limit}
Define a rescaled version of $H_\beta$ as
\begin{align}
H^{\text{r}}_\beta(q) = \frac{\|q\|_\beta - 1}{m^{1/\beta-1}-1}, \label{eq:H-beta2}
\end{align}
if $\beta \in (0,1)\cup (1,\infty)$, and
$H^{\text{r}}_\beta(q) =\lim_{\beta^\prime \to \beta}  H^{\text{r}}_{\beta^\prime}(q)$, if $ \beta=0, 1, \infty$.
Then the following proper scoring rules are obtained.
\begin{itemize}
\item[(i)] Simultaneous exponential loss ($\beta=0$): $L^{\text{r}}_0 (j,q) = ( \prod_{k=1}^m \frac{q_k}{q_j})^{1/m}$ corresponding to $H^{\text{r}}_0(q) =   m(\prod_{j=1}^m q_j)^{1/m}$.

\item[(ii)] Pairwise exponential loss ($\beta=1/2$): $L^{\text{r}}_{1/2} (j,q) = (m-1)^{-1}\sum_{k\in[m],k\not=j} \sqrt{\frac{q_k}{q_j}}$ corresponding to
$H^{\text{r}}_{\frac{1}{2}} (q) = (m-1)^{-1}(\|q\|_{\frac{1}{2}}-1)$.

\item[(iii)] Multinomial  likelihood loss ($\beta=1$): $L^{\text{r}}_1 (j,q) = - (\log{m})^{-1}\log q_j $ corresponding to $H^{\text{r}}_1(q) =  - (\log{m})^{-1}\sum_{j=1}^m q_j\log q_j$.

\item[(iv)] Multi-class zero-one loss ($\beta=\infty$): $L^{\text{r}}_\infty (j,q) = (1-m^{-1})^{-1} 1\{j\neq \argmax_{k\in[m]}q_k\}$ corresponding to $H^{\text{r}}_\infty (q) =  (1-m^{-1})^{-1}(1-\max_{ j \in [m]}q_j)$.
\end{itemize}
Moreover, with a multinomial logistic link (\ref{eq:multinomial-link}), the composite loss $L^{\text{r}}_\beta (j,q^h)$ is convex in $h$ if $\beta \in [0,1]$, but non-convex in $h$ if $\beta >1$.
\end{pro}

There are several interesting features. First, with a multinomial logistic link (\ref{eq:multinomial-link}), the scoring rule $L^{\text{r}}_0 (j,q)$ leads to a composite loss
\begin{align*}
L^{\text{r}}_0 (j,q^h) = \me^{ \frac{1}{m} \sum_{k=1}^m (h_k - h_j) },
\end{align*}
which coincides with the exponential loss in \cite{zou2008}.
For this reason, $L^{\text{r}}_0 (j,q)$ is called the simultaneous exponential loss.
Moreover, the scoring rule $L^{\text{r}}_{1/2} (j,q)$ yields, up to a multiplicative factor, the pairwise exponential loss $L^{\text{pw,s}}_e (j,q)$,
which is connected with the boosting algorithms in \cite{freund1997} and \cite{schapire1999} as mentioned earlier.
The logarithmic rule $L^{\text{r}}_1 (j,q)$ corresponds to the standard likelihood loss based on multinomial data.
%With a multinomial logistic link (\ref{eq:multinomial-link}), the resulting composite losses, $L^{\text{r}}_{1/2} (j,q^h)$ and $L^{\text{r}}_1 (j,q^h)$, are known to be convex in $h$.
Finally, the loss $L^{\text{r}}_\infty (j,q)$ obtained as $\beta\to\infty$ recovers the zero-one loss, which is a proper scoring rule (although not strictly proper).
Further research is desired to study relative merits of these losses.

\section{Hinge-like losses} \label{sec:hinge-like}

The purpose of this section is three-fold. We derive novel hinge-like, convex losses which induce the same generalized entropy as the zero-one, or more generally,
cost-weighted classification loss in multi-class settings. Our hinge-like losses are uniformly lower (after suitable alignment) and geometrically simpler (with fewer non-differentiable ridges)
than related hinge-like losses in \cite{lee2004} and \cite{duchi2018}.
Moreover, we show that similar classification regret bounds are achieved by our hinge-like losses and those in \cite{lee2004} and \cite{duchi2018}.
These regret bounds give a quantitative guarantee on classification calibration
as studied in \cite{zhang2004a} and \cite{tewari2007} among others.
\tzq{Finally, we provide a general characterization of losses with the same generalized entropy as the zero-one loss
and establish a general classification regret bound for all such losses, beyond the hinge-like losses specifically constructed.}

\subsection{Construction of hinge-like losses} \label{sec:hinge-construction}

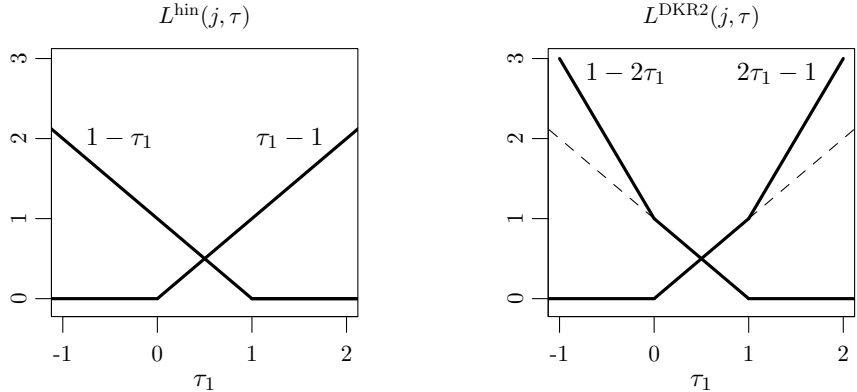
\begin{figure} %[!h]
\centering
% \resizebox{.9\textwidth}{!}{
% Created by tikzDevice version 0.12.3 on 2020-03-31 13:32:42
% !TEX encoding = UTF-8 Unicode
\begin{tikzpicture}[x=1pt,y=1pt]
\definecolor{fillColor}{RGB}{255,255,255}
\path[use as bounding box,fill=fillColor,fill opacity=0.00] (0,0) rectangle (375.80,173.45);
\begin{scope}
\path[clip] (  0.00,  0.00) rectangle (187.90,173.45);
\definecolor{drawColor}{RGB}{0,0,0}

\node[text=drawColor,anchor=base,inner sep=0pt, outer sep=0pt, scale=  0.90] at ( 93.95, 21.60) {$\tau_1$};
\end{scope}
\begin{scope}
\path[clip] ( 36.00, 48.00) rectangle (151.90,149.45);
\definecolor{drawColor}{RGB}{0,0,0}

\path[draw=drawColor,line width= 1.2pt,line join=round,line cap=round] (  0.00,149.52) -- (111.84, 54.79);

\path[draw=drawColor,line width= 1.2pt,line join=round,line cap=round] (111.84, 54.79) -- (219.15, 54.79);

\path[draw=drawColor,line width= 1.2pt,line join=round,line cap=round] (183.38,145.69) -- ( 76.06, 54.79);

\path[draw=drawColor,line width= 1.2pt,line join=round,line cap=round] ( 76.06, 54.79) -- (  0.00, 54.79);

\path[draw=drawColor,line width= 1.2pt,line join=round,line cap=round] (111.84, 54.79) -- (219.15, 54.79);

\node[text=drawColor,anchor=base,inner sep=0pt, outer sep=0pt, scale=  0.90] at ( 61.76,113.14) {$1-\tau_1$};

\node[text=drawColor,anchor=base,inner sep=0pt, outer sep=0pt, scale=  0.90] at (126.15,113.14) {$\tau_1-1$};
\end{scope}
\begin{scope}
\path[clip] (  0.00,  0.00) rectangle (375.80,173.45);
\definecolor{drawColor}{RGB}{0,0,0}

\path[draw=drawColor,line width= 0.4pt,line join=round,line cap=round] ( 40.29, 48.00) -- (151.90, 48.00);

\path[draw=drawColor,line width= 0.4pt,line join=round,line cap=round] ( 40.29, 48.00) -- ( 40.29, 42.00);

\path[draw=drawColor,line width= 0.4pt,line join=round,line cap=round] ( 76.06, 48.00) -- ( 76.06, 42.00);

\path[draw=drawColor,line width= 0.4pt,line join=round,line cap=round] (111.84, 48.00) -- (111.84, 42.00);

\path[draw=drawColor,line width= 0.4pt,line join=round,line cap=round] (147.61, 48.00) -- (147.61, 42.00);

\node[text=drawColor,anchor=base,inner sep=0pt, outer sep=0pt, scale=  0.80] at ( 40.29, 31.20) {-1};

\node[text=drawColor,anchor=base,inner sep=0pt, outer sep=0pt, scale=  0.80] at ( 76.06, 31.20) {0};

\node[text=drawColor,anchor=base,inner sep=0pt, outer sep=0pt, scale=  0.80] at (111.84, 31.20) {1};

\node[text=drawColor,anchor=base,inner sep=0pt, outer sep=0pt, scale=  0.80] at (147.61, 31.20) {2};

\path[draw=drawColor,line width= 0.4pt,line join=round,line cap=round] ( 36.00, 48.00) -- ( 36.00,145.69);

\path[draw=drawColor,line width= 0.4pt,line join=round,line cap=round] ( 36.00, 54.79) -- ( 30.00, 54.79);

\path[draw=drawColor,line width= 0.4pt,line join=round,line cap=round] ( 36.00, 85.09) -- ( 30.00, 85.09);

\path[draw=drawColor,line width= 0.4pt,line join=round,line cap=round] ( 36.00,115.39) -- ( 30.00,115.39);

\path[draw=drawColor,line width= 0.4pt,line join=round,line cap=round] ( 36.00,145.69) -- ( 30.00,145.69);

\node[text=drawColor,rotate= 90.00,anchor=base,inner sep=0pt, outer sep=0pt, scale=  0.80] at ( 26.40, 54.79) {0};

\node[text=drawColor,rotate= 90.00,anchor=base,inner sep=0pt, outer sep=0pt, scale=  0.80] at ( 26.40, 85.09) {1};

\node[text=drawColor,rotate= 90.00,anchor=base,inner sep=0pt, outer sep=0pt, scale=  0.80] at ( 26.40,115.39) {2};

\node[text=drawColor,rotate= 90.00,anchor=base,inner sep=0pt, outer sep=0pt, scale=  0.80] at ( 26.40,145.69) {3};

\path[draw=drawColor,line width= 0.4pt,line join=round,line cap=round] ( 36.00, 48.00) --
	(151.90, 48.00) --
	(151.90,149.45) --
	( 36.00,149.45) --
	( 36.00, 48.00);
\end{scope}
\begin{scope}
\path[clip] (  0.00,  0.00) rectangle (187.90,173.45);
\definecolor{drawColor}{RGB}{0,0,0}

\node[text=drawColor,anchor=base,inner sep=0pt, outer sep=0pt, scale=  0.80] at ( 93.95,159.05) { $L^{\text{hin}}(j,\tau)$};
\end{scope}
\begin{scope}
\path[clip] (187.90,  0.00) rectangle (375.80,173.45);
\definecolor{drawColor}{RGB}{0,0,0}

\node[text=drawColor,anchor=base,inner sep=0pt, outer sep=0pt, scale=  0.90] at (281.85, 21.60) {$\tau_1$};
\end{scope}
\begin{scope}
\path[clip] (223.90, 48.00) rectangle (339.80,149.45);
\definecolor{drawColor}{RGB}{0,0,0}

\path[draw=drawColor,line width= 1.2pt,line join=round,line cap=round] (156.65, 54.79) -- (263.97, 54.79);

\path[draw=drawColor,line width= 1.2pt,line join=round,line cap=round] (263.97, 54.79) -- (299.74, 85.09);

\path[draw=drawColor,line width= 1.2pt,line join=round,line cap=round] (299.74, 85.09) -- (335.51,145.69);

\path[draw=drawColor,line width= 1.2pt,line join=round,line cap=round] (299.74, 54.79) -- (371.28, 54.79);

\path[draw=drawColor,line width= 1.2pt,line join=round,line cap=round] (263.97, 85.09) -- (299.74, 54.79);

\path[draw=drawColor,line width= 1.2pt,line join=round,line cap=round] (263.97, 85.09) -- (228.19,145.69);

\path[draw=drawColor,line width= 0.4pt,dash pattern=on 4pt off 4pt ,line join=round,line cap=round] (299.74, 85.09) -- (371.28,145.69);

\path[draw=drawColor,line width= 0.4pt,dash pattern=on 4pt off 4pt ,line join=round,line cap=round] (263.97, 85.09) -- (192.42,145.69);

\node[text=drawColor,anchor=base,inner sep=0pt, outer sep=0pt, scale=  0.90] at (253.24,137.38) {$1-2\tau_1$};

\node[text=drawColor,anchor=base,inner sep=0pt, outer sep=0pt, scale=  0.90] at (310.47,137.38) {$2\tau_1-1$};
\end{scope}
\begin{scope}
\path[clip] (187.90,  0.00) rectangle (375.80,173.45);
\definecolor{drawColor}{RGB}{0,0,0}

\node[text=drawColor,anchor=base,inner sep=0pt, outer sep=0pt, scale=  0.80] at (281.85,159.05) { $L^{ \text{DKR2}}(j,\tau)$};
\end{scope}
\begin{scope}
\path[clip] (  0.00,  0.00) rectangle (375.80,173.45);
\definecolor{drawColor}{RGB}{0,0,0}

\path[draw=drawColor,line width= 0.4pt,line join=round,line cap=round] (228.19, 48.00) -- (339.80, 48.00);

\path[draw=drawColor,line width= 0.4pt,line join=round,line cap=round] (228.19, 48.00) -- (228.19, 42.00);

\path[draw=drawColor,line width= 0.4pt,line join=round,line cap=round] (263.97, 48.00) -- (263.97, 42.00);

\path[draw=drawColor,line width= 0.4pt,line join=round,line cap=round] (299.74, 48.00) -- (299.74, 42.00);

\path[draw=drawColor,line width= 0.4pt,line join=round,line cap=round] (335.51, 48.00) -- (335.51, 42.00);

\node[text=drawColor,anchor=base,inner sep=0pt, outer sep=0pt, scale=  0.80] at (228.19, 31.20) {-1};

\node[text=drawColor,anchor=base,inner sep=0pt, outer sep=0pt, scale=  0.80] at (263.97, 31.20) {0};

\node[text=drawColor,anchor=base,inner sep=0pt, outer sep=0pt, scale=  0.80] at (299.74, 31.20) {1};

\node[text=drawColor,anchor=base,inner sep=0pt, outer sep=0pt, scale=  0.80] at (335.51, 31.20) {2};

\path[draw=drawColor,line width= 0.4pt,line join=round,line cap=round] (223.90, 48.00) -- (223.90,145.69);

\path[draw=drawColor,line width= 0.4pt,line join=round,line cap=round] (223.90, 54.79) -- (217.90, 54.79);

\path[draw=drawColor,line width= 0.4pt,line join=round,line cap=round] (223.90, 85.09) -- (217.90, 85.09);

\path[draw=drawColor,line width= 0.4pt,line join=round,line cap=round] (223.90,115.39) -- (217.90,115.39);

\path[draw=drawColor,line width= 0.4pt,line join=round,line cap=round] (223.90,145.69) -- (217.90,145.69);

\node[text=drawColor,rotate= 90.00,anchor=base,inner sep=0pt, outer sep=0pt, scale=  0.80] at (214.30, 54.79) {0};

\node[text=drawColor,rotate= 90.00,anchor=base,inner sep=0pt, outer sep=0pt, scale=  0.80] at (214.30, 85.09) {1};

\node[text=drawColor,rotate= 90.00,anchor=base,inner sep=0pt, outer sep=0pt, scale=  0.80] at (214.30,115.39) {2};

\node[text=drawColor,rotate= 90.00,anchor=base,inner sep=0pt, outer sep=0pt, scale=  0.80] at (214.30,145.69) {3};

\path[draw=drawColor,line width= 0.4pt,line join=round,line cap=round] (223.90, 48.00) --
	(339.80, 48.00) --
	(339.80,149.45) --
	(223.90,149.45) --
	(223.90, 48.00);
\end{scope}
\end{tikzpicture}
% }
\vspace{-.4in}
\caption{Two-class hinge loss (left) and hinge-like loss in \cite{duchi2018}.}\label{fig:hinge-m2}
\end{figure}

While the proper scoring rules discussed in Section~\ref{sec:scoring-examples} are based on Proposition~\ref{pro:f-loss3}, our construction of hinge-like losses
relies on Proposition~\ref{pro:f-loss}. To see the difference, it is helpful to consider the two-class setting.
The generalized entropy for the zero-one loss is $H^{\text{zo}}(\eta) = \min(\eta_1,\eta_2)$ and the dissimilarity function is $f^{\text{zo}}(t_1) = - \min(1,t_1)$.
For this choice of $f$, it remains valid to apply Proposition~\ref{pro:f-loss3}. With $\partial f^{\text{zo}}(t_1) = -1\{t_1\le 1\}$, the resulting loss can be shown to be
$L_{f3} (1, q) = 1\{ q_1 \le q_2\}$ and $L_{f3} (2, q) = 1\{q_2 > q_1\}$, which is just the zero-one loss with action $q\in \Delta_2$. Such a discontinuous loss is computationally intractable for training,
even though it is a proper but not strictly proper scoring rule (consistently with Proposition~\ref{pro:f-loss3}).
In contrast, as illustrated in Figure~\ref{fig:hinge-m2}, the popular hinge loss can be defined, for notational consistency with our later results, such that for $\tau\in\bbR$,
\begin{align}
L^{\text{hin}}(1, \tau) = \max(0, 1-\tau), \quad L^{\text{hin}}(2,\tau) = \max(0, \tau) . \label{eq:hinge}
\end{align}
which is continuous and convex  in $\tau$ and known to yield the same generalized entropy $H^{\text{zo}}$ as the zero-one loss (Nguyen et al.~2009).
In the following, we show that Proposition~\ref{pro:f-loss} can be leveraged to develop convex, hinge-like losses in multi-class settings.

\vspace{.1in}\textbf{Application of Proposition~\ref{pro:f-loss}.}
Our application of Proposition~\ref{pro:f-loss} is facilitated by the following lemma, which gives the conjugate function of the dissimilarity function $f^{\text{cw}}$,
corresponding to the generalized entropy $H^{\text{cw}}$ in (\ref{eq:cw-H}) for the cost-weighted classification loss $L^{\text{cw}}$.
By definition (\ref{eq:Htof}), the dissimilarity function $f^{\text{cw}}$ can be calculated as
\begin{align*}
f^{\text{cw}}(t) = - \min_{k\in[m]} \, C_k^\T \tilde t ,
\end{align*}
where $\tilde t =(t^\T,1)^\T= (t_1, \ldots, t_{m-1},1)^\T$ and, as before, $C=(C_1, \ldots, C_m)$ is a column representation of the cost matrix
for the cost-weighted classification loss $L^{\text{cw}}$.

\begin{lem} \label{lem:cw-conj}
The conjugate of the convex function $f^{\text{cw}}$ is
\begin{align*}
f^{\text{cw} \,*} (s) =
\min_{ \{\lambda\in\Delta_m: s_j \le - (C\lambda)_j, j\in [m-1]\} } (C\lambda)_m  ,
\end{align*}
where $(C\lambda)_j$ denotes the $j$th component of $C\lambda$ for $j \in [m]$, and the minimum over an empty set is defined as $\infty$.
\end{lem}

From Lemma~\ref{lem:cw-conj}, the domain of $f^{\text{cw} *}$ is a strict subset of $ \bbR^{m-1}$, a phenomenon mentioned earlier in the discussion of Proposition~\ref{pro:f-loss}:
$$
\dom(f^{\text{cw} *} ) = \{ s \in \bbR ^{m-1}:  s_j \le - (C\lambda)_j, j\in [m-1] \mbox{ for some } \lambda\in\Delta_m \}.
$$
The following loss can be obtained from Proposition~\ref{pro:f-loss} with the convex function $f= f^{\text{cw}}$
and further simplification with a reparametrization $s_j = -(C\lambda)_j$ for $j \in [m-1]$.

\begin{lem} \label{lem:cw-conj2}
Define a loss with action space $\mathcal A = \Delta_m$ such that for $\lambda \in \Delta_m$,
\begin{align}
 L^{\text{cw2}} (j, \lambda) = \left\{ \begin{array}{cl}
 (C\lambda)_j, & j \in [m-1],\\
 (C \lambda)_m, & j=m,
\end{array} \right. \label{eq:cw-loss2}
\end{align}
Then the loss $L^{\text{cw2}}$ induces the same generalized entropy $H^{\text{cw}}$ in (\ref{eq:cw-H}) as does the cost-weighted classification loss $L^{\text{cw}}$.
\end{lem}

It is interesting that the loss $L^{\text{cw2}}$ is defined with actions restricted to the probability simplex $\Delta_m$. But $L^{\text{cw2}}$ is not a proper scoring rule,
because in general
\begin{align*}
&\inf_{\lambda \in \Delta_m} \left\{ \sum_{j=1}^m \eta_j L^{\text{cw2}} (j, \lambda) \right\} = \inf_{\lambda\in \Delta_m} \eta^\T C\lambda = \min_{k\in[m]} \eta^\T C_k \\
& \not= \eta^\T C \eta = \sum_{j=1}^m \eta_j L^{\text{cw2}} (j, \eta).
\end{align*}
by Lemma~\ref{lem:cw-conj2}. In fact, the minimum risk in the first line is achieved by $\lambda$ equal to a basis vector $e_l\in \Delta_m$ such that $\eta^\T C_l = \min_{k\in[m]} \eta^\T C_k$.

\vspace{.1in}\textbf{Extension beyond the probability simplex.}
The loss $L^{\text{cw2}}(j,\lambda)$ is convex (more precisely, linear!) in its action
$\lambda$ when restricted to $ \Delta_m$. To handle this restriction, there are several possible approaches.
One is to introduce a link function such as the multinomial logistic link $\lambda^h =(\lambda^h_1,\ldots, \lambda^h_m)^\T$, where $\lambda^h_j = \exp(h_j)/\sum_{k=1}^m \exp(h_k)$
with $h=(h_1,\ldots,h_{m-1})^\T \in \bbR^{m-1}$ unrestricted  and $h_m=0$ fixed.
But the resulting loss $L^{\text{cw2}}(j, \lambda^h)$ would be non-convex in $h$.
Another approach is to define a trivial extension of $L^{\text{cw2}}$ such that  $L^{\text{cw2}}(j,\lambda) = \infty$ whenever $\lambda$ lies outside the restricted set $\Delta_m$.
But for numerical implementation with this extension, either a link function such as the multinomial logistic link would still be needed,
or the predicted action for a new observation is likely to lie outside the probability simplex $\Delta_m$, which then requires additional treatment.
By comparison, our approach is to carefully construct an extension of $L^{\text{cw2}}$
which remains convex in its action and induces the same generalized entropy $H^{\text{cw}}$, while avoiding the infinity value outside the restricted set $\Delta_m$.

The version of $L^{\text{cw2}}$ in (\ref{eq:cw-loss2}) with $C = 1_m 1_m^\T - I_m$ as in the zero-one loss is
\begin{align}
 L^{\text{zo2}} (j, \lambda) = 1- \lambda_j, \quad j \in [m], \lambda \in \Delta_m . \label{eq:zo-loss2}
\end{align}
In the two-class setting, the hinge loss $L^{\text{hin}}$ can be shown to be a desired convex extension of the loss $L^{\text{zo2}}$,
considered a function of $j$ and $\lambda_1$:
\begin{align*}
 L^{\text{zo2}} (1, \lambda) = 1-\lambda_1, \quad  L^{\text{zo2}} (2, \lambda) = \lambda_1, \quad \lambda_1 \in [0,1].
\end{align*}
See Figure~\ref{fig:hinge-m2} for an illustration. In multi-class settings, our first extension of the loss $L^{\text{cw2}}$ is as follows, related to
the multi-class hinge-like loss in \cite{lee2004}.

\begin{pro} \label{pro:cw-conj3}
Define a loss with action space $\mathcal A = \bbR^{m-1}$ such that for $\tau \in \bbR^{m-1}$,
\begin{align}
L^{\text{cw3}} (j, \tau) =   \left\{
\begin{array}{cl}
 c_{jm} (\tau^{(j)}_m) _+ + \sum_{k\in[m-1],k\not=j} c_{jk} \tau_{k+}, & \mbox{if } j \in [m-1],\\
 \sum_{k \in [m-1]} c_{mk} \tau_{k+} , & \mbox{if } j=m,
\end{array} \right. \label{eq:cw-loss3}
\end{align}
where $b_+ = \max(0,b)$ for $b \in \bbR$, and
\begin{align}
\textstyle{ \tau^{(j)}_m = 1- \tau_j - \sum_{k\in[m-1],k\not=j} \tau_{k+} }, \quad j \in [m-1]. \label{eq:tau-m}
\end{align}
Then $L^{\text{cw3}}(j,\tau)$ is convex in $\tau$, and coincides with $L^{\text{cw2}}(j,\tilde\tau)$ provided $\tilde\tau \in \Delta_m$,
where $\tilde \tau = (\tau_1,\ldots,\tau_{m-1}, 1-\sum_{k=1}^{m-1} \tau_k)^\T$.
Moreover, $L^{\text{cw3}}$ induces the same generalized entropy $H^{\text{cw}}$ in (\ref{eq:cw-H}) as does the cost-weighted classification loss $L^{\text{cw}}$.
\end{pro}

\begin{figure}[t] %[!h]
\centering
\resizebox{\textwidth}{!}{
\input{LLW2_zo3.tex}} \vspace{-.4in}
\caption{Three-class hinge-like losses $L^{\text{LLW2}}$ (top) and $L^{\text{zo3}}$ (bottom).
Regions separated by solid lines are associated with the function values indicated.}\label{fig:LLW2}
\end{figure}

A special case of the loss $L^{\text{cw3}}$ with the cost matrix $C= 1_m 1_m^\T - I_m$ as in the zero-one loss can be expressed such that
for $\tau=(\tau_1,\ldots,\tau_{m-1})^\T \in \bbR^{m-1}$,
\begin{align*}
 L^{\text{zo3}} (j,\tau)= \left\{
\begin{array}{cl}
 \max\left(1-\tau_j, \,\sum_{k\in[m-1],k\not=j} \tau_{k+} \right), & \mbox{if } j \in [m-1],\\
 \sum_{k \in [m-1]} \tau_{k+} , & \mbox{if } j=m,
\end{array} \right.
\end{align*}
where the summation over an empty set is defined as 0. Then
$L^{\text{zo3}}$ induces the same generalized entropy $H^{\text{zo}}$ in (\ref{eq:zo-H}) as does the zero-one loss $L^{\text{zo}}$.
In the two-class setting, the loss $L^{\text{zo3}}$ can be easily seen to coincide with the hinge loss (\ref{eq:hinge}).

We compare the new loss with the hinge-like loss in \cite{lee2004} corresponding to the zero-one loss with $C= 1_m 1_m^\T - I_m$, which is defined such that for $\gamma\in\bbR^m$,
\begin{align*}
 L^{\text{LLW}} (j,\gamma)=  \sum_{k\in[m], k\not=j} (1+\gamma_k)_+ , \quad j \in [m],
\end{align*}
subject to the restriction that $\sum_{k=1}^m \gamma_k = 0$. The general case of cost-weighted classification can be similarly discussed.
To facilitate comparison, a reparametrization of the loss $L^{\text{LLW}}$ can be obtained such that for $\tau\in\bbR^{m-1}$,
\begin{align*}
 L^{\text{LLW2}} (j, \tau)=  \left\{
\begin{array}{cl}
 \sum_{k\in[m-1],k\not=j} \tau_{k+} + \left( 1-\sum_{k\in[m-1]} \tau_k \right)_+ , & \mbox{if } j \in [m-1],\\
 \sum_{k \in [m-1]} \tau_{k+} , & \mbox{if } j=m,
\end{array} \right.
\end{align*}
In the Supplement, it is shown that $L^{\text{LLW2}} (j,\tau) = L^{\text{LLW}} (j,\gamma)/m$ for $j \in [m]$,
provided that $\tau_k = (1+\gamma_k)/m$ for $k\in [m-1]$.
Figure~\ref{fig:LLW2} illustrates the two losses $L^{\text{zo3}}$ and $L^{\text{LLW2}}$  in the three-class setting.
The loss $L^{\text{zo3}}$ is a tighter convex extension than $L^{\text{LLW2}}$ from $L^{\text{zo2}}$ in (\ref{eq:zo-loss2}),
and $L^{\text{zo3}}(j,\tau)$ is geometrically simpler with fewer non-differentiable ridges than $L^{\text{LLW2}}(j,\tau)$ for $j \in [m-1]$.
\tzq{See the Supplement for further discussion.}

There are various ways in which the loss $L^{\text{cw2}}$ can be extended from the probability simplex $\Delta_m$ to $\bbR^m$.
We describe another extension, related to the multi-class hinge-like loss in \cite{duchi2018} associated with the zero-one loss.
The general case of cost-weighted classification can be handled through the transformation (\ref{eq:transform-loss}) in Section~\ref{sec:regret-bd}, although
such a general construction is not discussed in \cite{duchi2018}.

\begin{pro} \label{pro:cw-conj4}
Define a loss with action space $\mathcal A = \bbR^{m-1}$ such that for $\tau \in \bbR^{m-1}$,
\begin{align}
L^{\text{zo4}} (j,\tau)=
 1-\tilde \tau_j + S_\tau^{(j)} , \quad j \in [m], \label{eq:zo-loss4}
\end{align}
where $(\tilde\tau_1,\ldots,\tilde\tau_{m-1}) = (\tau_1,\ldots, \tau_{m-1})$,  $\tilde\tau_m=1-\sum_{k=1}^{m-1} \tau_k$, and for $j\in[m]$,
\begin{align*}
S_\tau^{(j)} =\max  \left\{0, \tilde\tau_j-1, \frac{\tilde\tau_j+\tilde\tau_{j(1)} -1}{2},  \ldots,
 \frac{\tilde\tau_j+\tilde\tau_{j(1)}+\cdots+\tilde\tau_{j(m-2)}-1}{m-1}  \right\},
\end{align*}
with $\tilde\tau_{j(1)} \ge \ldots \ge \tilde\tau_{j(m-1)}$ the sorted components of $\tilde\tau=(\tilde\tau_1,\ldots,\tilde\tau_m)^\T$ excluding $\tilde\tau_j$.
Then $L^{\text{zo4}}(j,\tau)$ is convex in $\tau$, and coincides with $L^{\text{zo2}}(j,\tilde\tau)$ provided
$\tilde\tau\in \Delta_m$.
Moreover, $L^{\text{zo4}}$ induces the same generalized entropy $H^{\text{zo}}$ in (\ref{eq:zo-H}) as does the zero-one loss $L^{\text{zo}}$.
\end{pro}

\begin{figure}[t] %[!h]
\centering
\resizebox{\textwidth}{!}{
\input{Duchi_zo4.tex}} \vspace{-.4in}
\caption{Three-class hinge-like losses $L^{\text{DKR2}}$ (top) and $L^{\text{zo4}}$ (bottom).
Regions separated by solid lines are associated with the function values indicated.}\label{fig:DKR2}
\end{figure}

The hinge-like loss in \cite{duchi2018} is defined such that for $\gamma\in \bbR^m$,
\begin{align*}
 L^{\text{DKR}} (j,\gamma)=  1-\gamma_j + S_\gamma, \quad j \in [m],
\end{align*}
where
$S_\gamma=\max  \{\gamma_{(1)}-1, \frac{\gamma_{(1)} +\gamma_{(2)}-1}{2},  \ldots,
 \frac{\gamma_{(1)}+\cdots+\gamma_{(m)}-1}{m} \}$,
and $\gamma_{(1)} \ge \ldots \ge \gamma_{(m)}$ are the sorted components of $\gamma \in \bbR^m$.
This loss is invariant to any translation in $\gamma$, that is, $L^{\text{DKR}} (j,\gamma-b1_m) = L^{\text{DKR}} (j,\gamma)$ for any $b\in\bbR$.
It suffices to consider $L^{\text{DKR}}(j,\gamma)$ subject to the restriction that $\sum_{k=1}^m \gamma_k = 1$, or equivalently consider the loss
\begin{align*}
 L^{\text{DKR2}} (j,\tau)=  1-\tilde\tau_j + S_\tau, \quad j \in [m],
\end{align*}
where  $(\tilde\tau_1,\ldots,\tilde\tau_{m-1}) = (\tau_1,\ldots, \tau_{m-1})$,  $\tilde\tau_m=1-\sum_{k=1}^{m-1} \tau_k$, and
\begin{align*}
S_\tau=\max  \left\{0, \tilde\tau_{(1)}-1, \frac{\tilde\tau_{(1)} +\tilde\tau_{(2)}-1}{2},  \ldots,
 \frac{\tilde\tau_{(1)}+\cdots+\tilde\tau_{(m-1)}-1}{m-1} \right\} ,
\end{align*}
with $\tilde\tau_{(1)} \ge \ldots \ge \tilde\tau_{(m)}$ the sorted components of $\tilde\tau=(\tilde\tau_1,\ldots,\tilde\tau_m)^\T$.
There does not seem to be a direct transformation between the two losses $L^{\text{zo4}}$ and $L^{\text{DKR2}}$, in spite of their similar expressions.
An illustration is provided by Figures~\ref{fig:hinge-m2} and \ref{fig:DKR2} in two- and three-class settings.
The loss $L^{\text{zo4}}$ is a tighter convex extension than $L^{\text{DKR2}}$ from $L^{\text{zo2}}$ in (\ref{eq:zo-loss2}),
and $L^{\text{zo4}}(j,\tau)$ is geometrically simpler with fewer non-differentiable ridges than $L^{\text{LLW2}}(j,\tau)$ for $j \in [m]$.
\tzq{See the Supplement for further discussion.}

%\vspace{.1in}\textbf{Hinge-like regret bounds.}
\subsection{Regret bounds for hinge-like losses} \label{sec:hinge-regret-bd}

The preceding section mainly focuses on constructing multi-class hinge-like losses which induce the generalized entropy $L^{\text{zo}}$ or
$H^{\text{cw}}$ as does the zero-one or cost-weighted classification loss, while achieving certain desirable properties geometrically compared with
hinge-like losses in \cite{lee2004} and \cite{duchi2018}. Here we derive classification regret bounds, which compare
the regrets of our hinge-like losses with those of the zero-one and cost-weighted losses, where the actions are take from those of the hinge-like losses by a prediction mapping.
Such bounds provide a quantitative guarantee on classification calibration, a qualitative property which leads to infinite-sample classification consistency
under suitable technical conditions \citep{zhang2004a,tewari2007}.

\begin{pro} \label{pro:hinge-regret-bd}
The following regret bounds hold for the hinge-like losses
$L^{\text{cw3}}$ and $L^{\text{zo4}}$.
\begin{itemize}
\item[(i)] For $\eta\in\Delta_m$ and $\tau\in\bbR^{m-1}$, $m^{-1} B_{L^{\text{cw}}}  ( \eta, \tau^\dag) \le B_{L^{\text{cw3}} }  ( \eta, \tau)$, that is,
\begin{align}
\frac{1}{m} \left\{  \sum_{j=1}^m \eta_j L^{\text{cw}}  (j, \tau^\dag) - H^{\text{cw}}( \eta)  \right\} \le
\sum_{j=1}^m \eta_j L^{\text{cw3}}  (j, \tau) - H_{L^{\text{cw3}}}  ( \eta) , \label{eq:hinge-regret-bd}
\end{align}
where $ \tau^\dag = (\tau_1, \ldots, \tau_{m-1}, 1- \sum_{k=1}^{m-1} \tau_{k+} )^\T$.
\item[(ii)] For $\eta\in\Delta_m$ and $\tau\in\bbR^{m-1}$, $m^{-1} B_{L^{\text{zo}}}  ( \eta, \tilde\tau) \le B_{L^{\text{zo4}} }  ( \eta, \tau)$, that is,
\begin{align}
\frac{1}{m} \left\{  \sum_{j=1}^m \eta_j L^{\text{zo}}  (j, \tilde\tau) - H^{\text{zo}}( \eta)  \right\} \le
\sum_{j=1}^m \eta_j L^{\text{zo4}}  (j, \tau) - H_{L^{\text{zo4}}}  ( \eta) , \label{eq:hinge-regret-bd2}
\end{align}
where $ \tilde\tau = (\tau_1, \ldots, \tau_{m-1}, 1- \sum_{k=1}^{m-1} \tau_k )^\T$.
\end{itemize}
\end{pro}

The regret bounds (\ref{eq:hinge-regret-bd}) and (\ref{eq:hinge-regret-bd2}) directly lead to classification calibration,
which can be defined as follows, allowing a prediction mapping \citep{zhang2004a, tewari2007}.
\tzq{For a loss $L(j,\gamma)$ with action space $\mathcal A$, let $\sigma=(\sigma_1,\ldots,\sigma_m)^\T: \mathcal A\to \bbR^m$ be a  prediction mapping
which carries an action in $\mathcal A$ to a vector in $\bbR^m$, to be used as the corresponding action in the zero-one or cost-weighted classification loss.
The prediction mapping can be defined directly as the identity mapping, $\sigma(\gamma)=\gamma$, in the case of $\mathcal A \subset \bbR^m$,
but needs to convert an action $\gamma$ to a vector in $\bbR^m$ in the case of $\mathcal A \subset \bbR^{m-1}$.}
A loss $L(j,\gamma)$ with action space $\mathcal A$ and prediction mapping $\sigma(\cdot)$ is said to be classification calibrated for the zero-one loss if for any $\eta \in \Delta_m$ and $k\in[m]$ with $\eta_k < \max_{j\in[m]} \eta_j$,
\begin{align}
\inf_{\gamma \in \mathcal A}
\left\{ \sum_{j=1}^m \eta_j L(j,\gamma) - H_L(\eta): \sigma_k (\gamma) = \max_{j\in[m]} \sigma_j (\gamma) \right\} > 0 . \label{eq:class-calibration-zo}
\end{align}
For $L=L^{\text{zo4}}$ and $\sigma(\tau ) = \tilde\tau$, inequality (\ref{eq:hinge-regret-bd2}) implies that
the left-hand side of (\ref{eq:class-calibration-zo}) is no smaller than $(-\eta_k + \max_{j\in[m]} \eta_j)/m>0$. Hence
$L^{\text{zo4}}$ is classification calibrated for the zero-one loss.
Similarly, a loss $L(j,\gamma)$ with  action space $\mathcal A$ and prediction mapping $\sigma(\cdot)$
is said to be classification calibrated for cost-weighted classification with cost matrix $C$
if for any $\eta \in \Delta_m$ and $k\in[m]$ with $\eta^\T C_k > \max_{j\in[m]} \eta^\T C_j$,
\begin{align}
\inf_{\gamma \in \mathcal A}
\left\{ \sum_{j=1}^m \eta_j L(j,\gamma) - H_L(\eta): \sigma_k (\gamma) = \max_{j\in[m]} \sigma_j (\gamma) \right\} > 0. \label{eq:class-calibration-cw}
\end{align}
For $L=L^{\text{cw3}}$ and $\sigma(\tau ) = \tau^\dag$, inequality (\ref{eq:hinge-regret-bd}) implies that
the left-hand side of (\ref{eq:class-calibration-cw}) is no smaller than $(\eta^\T C_k - \min_{j\in[m]} \eta^\T C_j)/m>0$. Hence the loss
$L^{\text{cw3}}$ is classification calibrated with $\sigma(\tau ) = \tau^\dag$ for cost-weighted classification.

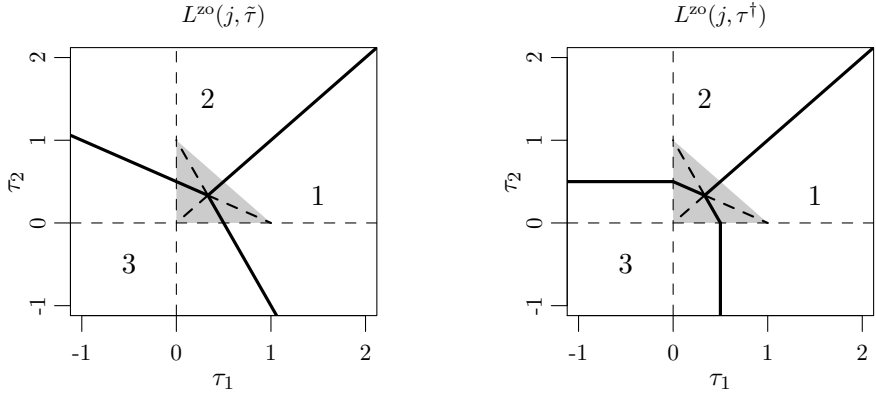
\begin{figure} %[!h]
\centering
% \resizebox{.8\textwidth}{!}{
% Created by tikzDevice version 0.12.3 on 2020-03-31 13:32:40
% !TEX encoding = UTF-8 Unicode
\begin{tikzpicture}[x=1pt,y=1pt]
\definecolor{fillColor}{RGB}{255,255,255}
\path[use as bounding box,fill=fillColor,fill opacity=0.00] (0,0) rectangle (375.80,173.45);
\begin{scope}
\path[clip] (  0.00,  0.00) rectangle (187.90,173.45);
\definecolor{drawColor}{RGB}{0,0,0}

\node[text=drawColor,anchor=base,inner sep=0pt, outer sep=0pt, scale=  0.90] at ( 93.95, 21.60) {$\tau_1$};

\node[text=drawColor,rotate= 90.00,anchor=base,inner sep=0pt, outer sep=0pt, scale=  0.90] at ( 16.80, 98.72) {$\tau_2$};
\end{scope}
\begin{scope}
\path[clip] ( 36.00, 48.00) rectangle (151.90,149.45);
\definecolor{fillColor}{gray}{0.80}

\path[fill=fillColor] ( 76.06, 83.07) --
	( 76.06,114.38) --
	(111.84, 83.07) --
	cycle;
\definecolor{drawColor}{RGB}{0,0,0}

\path[draw=drawColor,line width= 0.4pt,dash pattern=on 4pt off 4pt ,line join=round,line cap=round] ( 36.00, 83.07) -- (151.90, 83.07);

\path[draw=drawColor,line width= 0.4pt,dash pattern=on 4pt off 4pt ,line join=round,line cap=round] ( 76.06, 48.00) -- ( 76.06,149.45);

\path[draw=drawColor,line width= 1.2pt,line join=round,line cap=round] ( 87.99, 93.51) -- (141.40,  0.00);

\path[draw=drawColor,line width= 1.2pt,line join=round,line cap=round] ( 87.99, 93.51) -- (179.32,173.45);

\path[draw=drawColor,line width= 1.2pt,line join=round,line cap=round] ( 87.99, 93.51) -- (  0.00,132.01);

\path[draw=drawColor,line width= 0.8pt,dash pattern=on 4pt off 4pt ,line join=round,line cap=round] ( 87.99, 93.51) -- ( 76.06,114.38);

\path[draw=drawColor,line width= 0.8pt,dash pattern=on 4pt off 4pt ,line join=round,line cap=round] ( 87.99, 93.51) -- ( 76.06, 83.07);

\path[draw=drawColor,line width= 0.8pt,dash pattern=on 4pt off 4pt ,line join=round,line cap=round] ( 87.99, 93.51) -- (111.84, 83.07);

\node[text=drawColor,anchor=base,inner sep=0pt, outer sep=0pt, scale=  1.00] at (129.72, 90.03) {$1$};

\node[text=drawColor,anchor=base,inner sep=0pt, outer sep=0pt, scale=  1.00] at ( 87.99,126.56) {$2$};

\node[text=drawColor,anchor=base,inner sep=0pt, outer sep=0pt, scale=  1.00] at ( 58.18, 63.94) {$3$};
\end{scope}
\begin{scope}
\path[clip] (  0.00,  0.00) rectangle (375.80,173.45);
\definecolor{drawColor}{RGB}{0,0,0}

\path[draw=drawColor,line width= 0.4pt,line join=round,line cap=round] ( 40.29, 48.00) -- (147.61, 48.00);

\path[draw=drawColor,line width= 0.4pt,line join=round,line cap=round] ( 40.29, 48.00) -- ( 40.29, 42.00);

\path[draw=drawColor,line width= 0.4pt,line join=round,line cap=round] ( 76.06, 48.00) -- ( 76.06, 42.00);

\path[draw=drawColor,line width= 0.4pt,line join=round,line cap=round] (111.84, 48.00) -- (111.84, 42.00);

\path[draw=drawColor,line width= 0.4pt,line join=round,line cap=round] (147.61, 48.00) -- (147.61, 42.00);

\node[text=drawColor,anchor=base,inner sep=0pt, outer sep=0pt, scale=  0.80] at ( 40.29, 31.20) {-1};

\node[text=drawColor,anchor=base,inner sep=0pt, outer sep=0pt, scale=  0.80] at ( 76.06, 31.20) {0};

\node[text=drawColor,anchor=base,inner sep=0pt, outer sep=0pt, scale=  0.80] at (111.84, 31.20) {1};

\node[text=drawColor,anchor=base,inner sep=0pt, outer sep=0pt, scale=  0.80] at (147.61, 31.20) {2};

\path[draw=drawColor,line width= 0.4pt,line join=round,line cap=round] ( 36.00, 51.76) -- ( 36.00,145.69);

\path[draw=drawColor,line width= 0.4pt,line join=round,line cap=round] ( 36.00, 51.76) -- ( 30.00, 51.76);

\path[draw=drawColor,line width= 0.4pt,line join=round,line cap=round] ( 36.00, 83.07) -- ( 30.00, 83.07);

\path[draw=drawColor,line width= 0.4pt,line join=round,line cap=round] ( 36.00,114.38) -- ( 30.00,114.38);

\path[draw=drawColor,line width= 0.4pt,line join=round,line cap=round] ( 36.00,145.69) -- ( 30.00,145.69);

\node[text=drawColor,rotate= 90.00,anchor=base,inner sep=0pt, outer sep=0pt, scale=  0.80] at ( 26.40, 51.76) {-1};

\node[text=drawColor,rotate= 90.00,anchor=base,inner sep=0pt, outer sep=0pt, scale=  0.80] at ( 26.40, 83.07) {0};

\node[text=drawColor,rotate= 90.00,anchor=base,inner sep=0pt, outer sep=0pt, scale=  0.80] at ( 26.40,114.38) {1};

\node[text=drawColor,rotate= 90.00,anchor=base,inner sep=0pt, outer sep=0pt, scale=  0.80] at ( 26.40,145.69) {2};

\path[draw=drawColor,line width= 0.4pt,line join=round,line cap=round] ( 36.00, 48.00) --
	(151.90, 48.00) --
	(151.90,149.45) --
	( 36.00,149.45) --
	( 36.00, 48.00);
\end{scope}
\begin{scope}
\path[clip] (  0.00,  0.00) rectangle (187.90,173.45);
\definecolor{drawColor}{RGB}{0,0,0}

\node[text=drawColor,anchor=base,inner sep=0pt, outer sep=0pt, scale=  0.80] at ( 93.95,159.05) { $L^{\text{zo}}(j,\tilde{\tau})$};
\end{scope}
\begin{scope}
\path[clip] (187.90,  0.00) rectangle (375.80,173.45);
\definecolor{drawColor}{RGB}{0,0,0}

\node[text=drawColor,anchor=base,inner sep=0pt, outer sep=0pt, scale=  0.90] at (281.85, 21.60) {$\tau_1$};

\node[text=drawColor,rotate= 90.00,anchor=base,inner sep=0pt, outer sep=0pt, scale=  0.90] at (204.70, 98.72) {$\tau_2$};
\end{scope}
\begin{scope}
\path[clip] (223.90, 48.00) rectangle (339.80,149.45);
\definecolor{fillColor}{gray}{0.80}

\path[fill=fillColor] (263.97, 83.07) --
	(263.97,114.38) --
	(299.74, 83.07) --
	cycle;
\definecolor{drawColor}{RGB}{0,0,0}

\path[draw=drawColor,line width= 0.4pt,dash pattern=on 4pt off 4pt ,line join=round,line cap=round] (223.90, 83.07) -- (339.80, 83.07);

\path[draw=drawColor,line width= 0.4pt,dash pattern=on 4pt off 4pt ,line join=round,line cap=round] (263.97, 48.00) -- (263.97,149.45);

\path[draw=drawColor,line width= 1.2pt,line join=round,line cap=round] (275.89, 93.51) -- (367.22,173.45);

\path[draw=drawColor,line width= 1.2pt,line join=round,line cap=round] (275.89, 93.51) -- (281.85, 83.07);

\path[draw=drawColor,line width= 1.2pt,line join=round,line cap=round] (281.85, 83.07) -- (281.85,  0.00);

\path[draw=drawColor,line width= 1.2pt,line join=round,line cap=round] (275.89, 93.51) -- (263.97, 98.72);

\path[draw=drawColor,line width= 1.2pt,line join=round,line cap=round] (263.97, 98.72) -- ( 85.11, 98.72);

\path[draw=drawColor,line width= 0.8pt,dash pattern=on 4pt off 4pt ,line join=round,line cap=round] (275.89, 93.51) -- (263.97,114.38);

\path[draw=drawColor,line width= 0.8pt,dash pattern=on 4pt off 4pt ,line join=round,line cap=round] (275.89, 93.51) -- (263.97, 83.07);

\path[draw=drawColor,line width= 0.8pt,dash pattern=on 4pt off 4pt ,line join=round,line cap=round] (275.89, 93.51) -- (299.74, 83.07);

\node[text=drawColor,anchor=base,inner sep=0pt, outer sep=0pt, scale=  1.00] at (317.63, 90.03) {$1$};

\node[text=drawColor,anchor=base,inner sep=0pt, outer sep=0pt, scale=  1.00] at (275.89,126.56) {$2$};

\node[text=drawColor,anchor=base,inner sep=0pt, outer sep=0pt, scale=  1.00] at (246.08, 63.94) {$3$};
\end{scope}
\begin{scope}
\path[clip] (  0.00,  0.00) rectangle (375.80,173.45);
\definecolor{drawColor}{RGB}{0,0,0}

\path[draw=drawColor,line width= 0.4pt,line join=round,line cap=round] (228.19, 48.00) -- (335.51, 48.00);

\path[draw=drawColor,line width= 0.4pt,line join=round,line cap=round] (228.19, 48.00) -- (228.19, 42.00);

\path[draw=drawColor,line width= 0.4pt,line join=round,line cap=round] (263.97, 48.00) -- (263.97, 42.00);

\path[draw=drawColor,line width= 0.4pt,line join=round,line cap=round] (299.74, 48.00) -- (299.74, 42.00);

\path[draw=drawColor,line width= 0.4pt,line join=round,line cap=round] (335.51, 48.00) -- (335.51, 42.00);

\node[text=drawColor,anchor=base,inner sep=0pt, outer sep=0pt, scale=  0.80] at (228.19, 31.20) {-1};

\node[text=drawColor,anchor=base,inner sep=0pt, outer sep=0pt, scale=  0.80] at (263.97, 31.20) {0};

\node[text=drawColor,anchor=base,inner sep=0pt, outer sep=0pt, scale=  0.80] at (299.74, 31.20) {1};

\node[text=drawColor,anchor=base,inner sep=0pt, outer sep=0pt, scale=  0.80] at (335.51, 31.20) {2};

\path[draw=drawColor,line width= 0.4pt,line join=round,line cap=round] (223.90, 51.76) -- (223.90,145.69);

\path[draw=drawColor,line width= 0.4pt,line join=round,line cap=round] (223.90, 51.76) -- (217.90, 51.76);

\path[draw=drawColor,line width= 0.4pt,line join=round,line cap=round] (223.90, 83.07) -- (217.90, 83.07);

\path[draw=drawColor,line width= 0.4pt,line join=round,line cap=round] (223.90,114.38) -- (217.90,114.38);

\path[draw=drawColor,line width= 0.4pt,line join=round,line cap=round] (223.90,145.69) -- (217.90,145.69);

\node[text=drawColor,rotate= 90.00,anchor=base,inner sep=0pt, outer sep=0pt, scale=  0.80] at (214.30, 51.76) {-1};

\node[text=drawColor,rotate= 90.00,anchor=base,inner sep=0pt, outer sep=0pt, scale=  0.80] at (214.30, 83.07) {0};

\node[text=drawColor,rotate= 90.00,anchor=base,inner sep=0pt, outer sep=0pt, scale=  0.80] at (214.30,114.38) {1};

\node[text=drawColor,rotate= 90.00,anchor=base,inner sep=0pt, outer sep=0pt, scale=  0.80] at (214.30,145.69) {2};

\path[draw=drawColor,line width= 0.4pt,line join=round,line cap=round] (223.90, 48.00) --
	(339.80, 48.00) --
	(339.80,149.45) --
	(223.90,149.45) --
	(223.90, 48.00);
\end{scope}
\begin{scope}
\path[clip] (187.90,  0.00) rectangle (375.80,173.45);
\definecolor{drawColor}{RGB}{0,0,0}

\node[text=drawColor,anchor=base,inner sep=0pt, outer sep=0pt, scale=  0.80] at (281.85,159.05) {$L^{\text{zo}}(j,\tau^\dagger)$};
\end{scope}
\end{tikzpicture}
% }
\vspace{-.4in}
\caption{\tzq{Classification using the prediction mapping $\tilde\tau$ (left) or $\tau^\dag$ (right), defined in Proposition~\ref{pro:hinge-regret-bd}, with $\tau=(\tau_1,\tau_2)\in \bbR^2$
for $m=3$.} Each region separated by solid lines from others is classified by the index of a maximum component of $\tilde\tau$ or $\tau^\dag$.}\label{fig:class-boundary}
\end{figure}

There is an interesting feature in the regret bound (\ref{eq:hinge-regret-bd}) for $L^{\text{cw3}}(j,\tau)$, compared with the regret bound
(\ref{eq:hinge-regret-bd2}) for $L^{\text{zo4}}(j,\tau)$.
The prediction mapping associated with the loss $L^{\text{cw3}}(j,\tau)$ with $\tau\in\bbR^{m-1}$ is
$ \tau^\dag = (\tau_1, \ldots, \tau_{m-1}, 1-\sum_{k=1}^{m-1} \tau_{k+})^\T$, whose components may sum to less than one,
instead of $\tilde \tau = (\tau_1, \ldots, \tau_{m-1}, 1-\sum_{k=1}^{m-1} \tau_k)^\T$, whose components necessarily sum to one.
Although this difference warrants further study, using $\tau^\dag$ instead of $\tilde\tau$ for classification ensures that
the predicted value for $m$th class, $1-\sum_{k=1}^{m-1} \tau_{k+}$, is not affected by any negative components among $(\tau_1,\ldots,\tau_{m-1})$.
For example, if $m=3$ and $(\tau_1,\tau_2) = (.6,-.3)$, then $1 - \sum_{k=1}^2 \tau_{k+} = .4$ but $1 - \sum_{k=1}^2 \tau_k = .7$.
Using $\tau^\dag$ means that class $1$ is predicted, whereas using $\tilde\tau$ means that class $3$ is predicted, which seems to be artificially caused by the negative value of $\tau_2$.
See Figure~\ref{fig:class-boundary} for an illustration and \tzq{Proposition~\ref{pro:general-hinge-bd} for an explanation}.

The regret bounds (\ref{eq:hinge-regret-bd}) and (\ref{eq:hinge-regret-bd2}) for the losses $L^{\text{cw3}}$ and $L^{\text{zo4}}$ are similar to those for
the losses $L^{\text{LLW2}}$ and $L^{\text{DKR2}}$ in \cite{duchi2018}.
In fact, the regret bound for $L^{\text{LLW2}}$ in Duchi et al. can be seen as (\ref{eq:hinge-regret-bd}) with
$L^{\text{cw3}}  (j, \tau)$ and  $L^{\text{cw}}  (j, \tau^\dag)$  replaced by $L^{\text{LLW2}}  (j, \tau)$
and $L^{\text{cw}}  (j, \tilde \tau)$ respectively,
because $L^{\text{LLW2}}$ is $L^{\text{LLW}}$ multiplied by $m$ after a reparametrization noted earlier.
The regret bound for $L^{\text{DKR2}}$ in Duchi et al. can be seen as (\ref{eq:hinge-regret-bd2}) with
$L^{\text{zo4}}  (j, \tau)$ and  $L^{\text{zo}}  (j, \tilde\tau)$  replaced by $L^{\text{DKR2}}  (j, \tau)$
and $L^{\text{zo}}  (j, \tilde \tau)$ respectively.
Further research is desired to compare these hinge-like losses in both theory and empirical evaluation.

\subsection{\tzq{General characterization and regret bounds}} \label{sec:general-hinge-bd}

Our new hinge-like losses are explicitly derived to induce the same generalized entropy as the zero-one or cost-weighted classification loss,
and shown to achieve comparable regret bounds to those for existing hinge-like losses. \tzq{In this section, we provide a general result
indicating that all losses with the same generalized entropy as the zero-one loss achieve a classification regret bound similarly as in Proposition~\ref{pro:hinge-regret-bd}.
This result relies on a general characterization of such losses in terms of the value manifold defined below.}

For a loss $L(j,\gamma)$ with action space $\mathcal A$, the value manifold is defined as $\mathcal S_L = \overline{\mathrm{conv}} (\mathcal R_L)$, where
$\overline{\mathrm{conv}}$ denotes the closure of the convex hull and
\begin{align*}
& \mathcal R_L = \left\{ (L(1,\gamma), \ldots, L(m,\gamma))^\T: \gamma \in \mathcal A \right\}.
\end{align*}
The concept of the set $\mathcal R_L$ and its convex hull, $\mathrm{conv}(\mathcal R_L)$, also plays an important role in \cite{tewari2007},
where the admissibility of $\mathrm{conv}(\mathcal R_L)$ can be equivalently defined as that of $\mathcal S_L$ because
$\mathrm{conv}(\mathcal R_L)$ and $\mathcal S_L$ share the same boundary, denoted as $\partial \mathcal S_L$.
Then the generalized entropy of $L$ can be expressed such that for any $\eta \in \Delta_m$,
\begin{align}
H_L(\eta) = \inf_{\gamma \in \mathcal A} \left\{ \sum_{j=1}^m \eta_j L(j,\gamma) \right\} = \inf_{z \in \mathcal R_L} \eta^\T z = \inf_{z \in \mathcal S_L} \eta^\T z , \label{eq:entropy-manifold}
\end{align}
similarly as in \cite{tewari2007}, Eq.~(7).
%For the cost-weighted classification loss $L^{\text{cw}}$, the value manifold is denoted as
%\begin{align*}
%\mathcal S^{\text{cw}} = \left\{ \textstyle{\sum_{j=1}^m} \lambda_j C_j:  (\lambda_1,\ldots,\lambda_m)^\T \in \Delta_m \right\} .
%\end{align*}
%The set $\mathcal S^{\text{cw}}$ is an $(m-1)$-dimensional polytope in $\bbR^m$, with the vertices corresponding to the columns $C_1,\ldots,C_m$ in the cost matrix.
For the zero-one loss $L^{\text{zo}}$, the value manifold is denoted as
\begin{align*}
\mathcal S^{\text{zo}} = \left\{(z_1, \ldots, z_m)^\T: \textstyle{\sum_{j=1}^m} z_j = m-1 \mbox{ and } 0 \le z_1, \ldots, z_m \le 1 \right\} .
\end{align*}
The set $\mathcal S^{\text{zo}}$ is an $(m-1)$-dimensional polytope in $\bbR^m$, where each vertex is a $m$-dimensional vector with one component 0 and the remaining 1.

\begin{pro} \label{pro:value-manifold}
A loss $L(j,\gamma)$ induces the same generalized entropy as the zero-one loss, i.e., $H_L(\eta) = H^{\text{zo}}(\eta) =1- \max_{k\in[m]} \eta_k$ for $\eta \in \Delta_m$ if and only if
\begin{align}
\mathcal S^{\text{zo}} \subset \mathcal S_L \subset \mathcal S^{\text{zo}*}, \label{eq:val-manifold-inclusion}
\end{align}
where $S^{\text{zo}*} =\{z+b: z \in \mathcal S^{\text{zo}}, b \in \bbR^m_+\}$, also denoted as $\mathcal S^{\text{zo}} + \bbR^m_+$.
\end{pro}

\begin{figure}
	\begin{centering}
		\includegraphics[width=\textwidth]{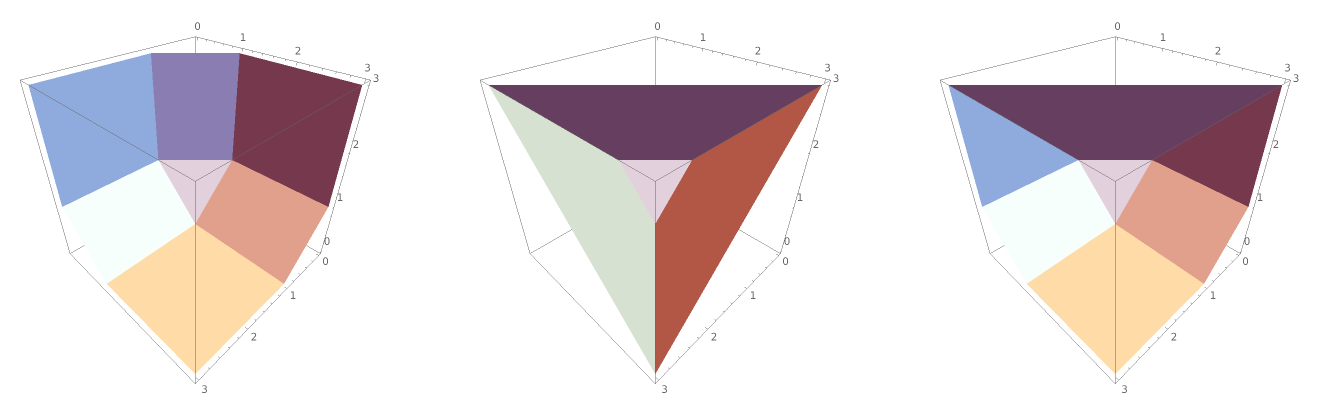}
		\caption{The boundaries of value manifolds for three-class hinge-like losses $L^{\text{DKR2}}$/$L^{\text{zo4}}$ (left), $L^{\text{LLW2}}$ (middle) and $L^{\text{zo3}}$ (right).
The triangle polytope in the center of each plot is the value manifold $\mathcal S^{\text{zo}}$ for the zero-one loss. The boundary of $\mathcal S^{\text{zo}*}$,
defined as $\mathcal S^{\text{zo}} + \bbR^m_+$, is the same as in the left plot.}
		\label{fig:value-manifold}
	\end{centering}
\end{figure}

Figure~\ref{fig:value-manifold} shows, in the three-class setting, the value manifolds for the zero-one and several hinge-like losses
with the same generalized entropy as the zero-one loss. These value manifolds all satisfy the inclusion property as stated in Proposition~\ref{pro:value-manifold}.

The following result establishes a general link from the generalized entropy of the zero-one loss to classification regret bounds.
The link involves a particular prediction mapping $\sigma_L(\gamma)$, defined from the negative values of a given loss $L$.
Such a prediction mapping is also exploited to study classification calibration in \cite{tewari2007}, with an additional assumption that the value manifold $\mathcal S_L$ is symmetric.

\begin{pro} \label{pro:general-hinge-bd}
Suppose that a loss $L(j,\gamma)$ with action space $\mathcal A$ induces the same generalized entropy as the zero-one loss. Then
for $\eta\in\Delta_m$ and $\gamma\in\mathcal A$, $m^{-1} B_{L^{\text{zo}}}  ( \eta, \sigma_L(\gamma) ) \le B_L  ( \eta, \gamma)$, that is,
\begin{align}
\frac{1}{m} \left\{  \sum_{j=1}^m \eta_j L^{\text{zo}}  (j, \sigma_L(\gamma) ) - H^{\text{zo}}( \eta)  \right\} \le
\sum_{j=1}^m \eta_j L  (j, \gamma) - H_L  ( \eta) , \label{eq:general-hinge-bd}
\end{align}
where $\sigma_L(\gamma) = ( -L(1,\gamma), \ldots, -L(m,\gamma))^\T$. Moreover, the loss $L(j,\gamma)$, with
the prediction mapping $\sigma_L$, is classification calibrated for the zero-one loss.
\end{pro}

Proposition~\ref{pro:general-hinge-bd} provides a theoretical support for our approach
in constructing hinge-like losses with the same generalized entropy as the zero-one loss in Section~\ref{sec:hinge-construction}.
The regret bound (\ref{eq:general-hinge-bd}) implies classification calibration similarly as discussed in Section~\ref{sec:hinge-regret-bd}. Hence our result
gives a more concrete sufficient condition for achieving classification calibration than in \cite{tewari2007}, Section 4.

By the nature of the zero-one loss, the regret bound (\ref{eq:general-hinge-bd}) remains valid
with $\sigma_L(\gamma)$ replaced by another prediction mapping $\sigma(\gamma)$ subject to the monotonicity property that the components in $\sigma(\gamma)$
are in the same order as in $\sigma_L(\gamma)$ for any $\gamma\in\mathcal A$.
Then the regret bounds for $L^{\text{zo3}}$, as a special case of $L^{\text{cw3}}$, and $L^{\text{zo4}}$ in Proposition~\ref{pro:hinge-regret-bd} can be deduced from (\ref{eq:general-hinge-bd}),
because the monotonicity property is satisfied by the prediction mappings $\tau^\dag$ and $\tilde\tau$ used in (\ref{eq:hinge-regret-bd}) and (\ref{eq:hinge-regret-bd2}).
See the Supplement for details.
Possible extensions of Proposition~\ref{pro:general-hinge-bd} to cost-weighted classification can be studied in future work.

\section{Regret bounds for proper scoring rules} \label{sec:regret-bd}

We return to proper scoring rules and derive classification regret bounds, which compare the regrets of the losses
with those of the corresponding zero-one and cost-weighted classification losses, similarly as in Proposition~\ref{pro:hinge-regret-bd} for hinge-like losses.
All such bounds are also called surrogate regret bounds, in the sense that the a proper scoring rule or a hinge-like loss
can be considered a surrogate criterion for the zero-one or cost-weighted classification loss.
Similarly as discussed in Section~\ref{sec:hinge-regret-bd}, these results provide
a quantitative guarantee on classification calibration \citep{zhang2004a, tewari2007}.

Compared with hinge-like losses, a potential gain in using proper scoring rules is that classification regret bounds can be obtained
with respect to a range of cost-weighted classification losses with different cost matrices $C$ for a proper scoring rule, defined independently of $C$.
The cost matrix is involved only to convert an action (in the form of a probability vector) from the scoring rule to a prediction for the cost-weighted classification loss.
See Corollary \ref{cor:regret-bd2-proper}.
In contrast, for the regret bound (\ref{eq:hinge-regret-bd}), the hinge-like loss $L^{\text{cw3}}$ depends on
the cost-matrix $C$ used in the classification loss $L^{\text{cw}}$.
A similar observation is made by \citet[Corollary 28]{reid2011} in two-class settings.

A general basis for deriving regret bounds, applicable to not just scoring rules but arbitrary losses $L(j,\gamma)$ with an action space $\mathcal A \subset \bbR^m$, can be seen as
\begin{align}
\psi ( B^{\text{zo}}(\eta, \gamma)) \le  B_L (\eta,\gamma), \label{eq:regret-general}
\end{align}
where $B^{\text{zo}}(\eta, \gamma) = B_{L^{\text{zo}}}(\eta, \gamma)$, the regret of the zero-one loss $L^{\text{zo}}(\eta, \gamma)$, and
$$
\psi(t ) = \inf_{ \eta^\prime\in\Delta_m, \gamma^\prime\in\mathcal A: \, B^{\text{zo}}(\eta^\prime, \gamma^\prime)=t } B_L (\eta^\prime, \gamma^\prime), \quad t\ge 0.
$$
In fact, (\ref{eq:regret-general}) is a tautology from the definition of $\psi$. Various regret bounds can be obtained by identifying
convenient lower bounds of $\psi$.
In the two-class setting, the regret for the zero-one loss is
$B^{\text{zo}}(\eta, \gamma)=  | 2\eta_1 - 1| 1\{ (2\eta_1 -1) (\gamma_1 - \gamma_2) \le 0\}$ for $\eta = (\eta_1, \eta_2)^\T \in \Delta_2$
and $\gamma = (\gamma_1, \gamma_2)^\T$. For $t>0$, $B^{\text{zo}}(\eta, \gamma)=t$ means $\eta_1 = (1\pm t)/2$, and hence $\psi(t)$ can be simplified as
\begin{align*}
\psi^{\text{BJM}} (t) = \min\left\{ \inf_{\gamma^\prime: t(\gamma^\prime_1 - \gamma^\prime_2)\le 0} B_L^1 \left( \frac{1+t}{2}, \gamma^\prime \right),
\inf_{\gamma^\prime: t (\gamma^\prime_1 - \gamma^\prime_2) \ge 0}  B_L^1 \left( \frac{1-t}{2}, \gamma^\prime \right) \right\},
\end{align*}
where $B_L^1 (\eta_1, \gamma)$ denotes $B_L( \eta, \gamma)$ as a function of $(\eta_1,\gamma)$.
Moreover, $\psi^{\text{BJM}}(t)$ at $t=0$ also satisfies $\psi^{\text{BJM}} (0) =0 \le \psi(0)$.
Therefore, (\ref{eq:regret-general}) holds with  $\psi$ replaced by $\psi^{\text{BJM}}$:
\begin{align}
\psi^{\text{BJM}} ( B^{\text{zo}}(\eta, \gamma)) \le  B_L (\eta,\gamma). \label{eq:regret-BJM}
\end{align}
In the multi-class setting, the regret $B^{\text{zo}}(\eta, \gamma)$ does not admit a direct simplification.
Nevertheless, our results below for proper scoring rules can be seen as further manipulation of (\ref{eq:regret-general})
by exploiting the fact that the regret (\ref{eq:bregman}) is a Bregman divergence due to the canonical representation (\ref{eq:savage}) for proper scoring rules.

\begin{rem} \label{rem:BJM}
Replacing $\psi^{\text{BJM}}$ in (\ref{eq:regret-BJM}) by the greatest convex lower bound on $\psi^{\text{BJM}}$
(or the Fenchel biconjugate of $\psi^{\text{BJM}}$) recovers the regret bound in \citet[Theorem 1]{bartlett2006}
in the symmetric case where $L(1,\gamma) = L(2,-\gamma)$.
In general, there is a benefit from such a modification in the setting where covariates are restored, instead of being lifted out in most of our discussion.
For a regret bound in the form $\phi ( B^{\text{zo}}(\eta, \gamma)) \le  B_L (\eta,\gamma)$,
if $\phi$ is convex, then application of Jensen's inequality gives
\begin{align*}
& \phi \big[ E \{ B^{\text{zo}}(\eta(X), \gamma(X)) \}  \big] \le  E \big[\phi\{ B^{\text{zo}}(\eta(X), \gamma(X)) \} \big]\\
& \le  E\{ B_L (\eta(X),\gamma(X)) \},
\end{align*}
where $E \{ B^{\text{zo}}(\eta(X), \gamma(X))$ and $E\{B_L (\eta(X),\gamma(X))\}$ are the average regret over $X$.
\end{rem}

%The regret from the class-weighted classification loss $L^{\text{cw0}}(j,\gamma)$ is
%\begin{align*}
%& B^{\text{cw0}}(\eta, \gamma) = \sum_{j=1}^m \eta_j L^{\text{cw0}}(j, \gamma) - H^{\text{cw0}}(\eta)\\
%& = \max_{k\in[m]} \eta_k c_{k0} - \sum_{j=1}^m \eta_j c_{j0} 1\{j=\argmax_{k \in [m]} \gamma_k\}.
%\end{align*}

\subsection{\tzq{Zero-one classification}} \label{sec:zero-one}

\tzq{Before presenting our general regret bounds for proper scoring rules with respect to cost-weighted classification in Sections~\ref{sec:cost-transform}--\ref{sec:cost-indep},
we demonstrate novel implications of our general results in the simple but important setting of zero-one classification.}

For a proper scoring rule $L(j,q)$,
an application of our regret bound (\ref{eq:cw0-regret-bd-proper}) or (\ref{eq:cw0-regret-bd2-proper}) with respect to the zero-one loss with $C_0 = 1_m$ shows that
for any $\eta,q\in\Delta_m$,
\begin{align}
\underline \psi \left( B^{\text{zo}}(\eta, q) \right) \le  B_L (\eta, q) , \label{eq:zo-regret-bd-proper}
\end{align}
where $\underline \psi(\cdot)$ is defined as
\begin{align}
\underline \psi (t) = \inf_{\substack{\eta^\prime, q^\prime\in \Delta_m: \, \| \eta^\prime -q^\prime \|_{\infty2} = t,\\ \max_{j\in[m]} q^\prime_j \le 1/2 }} B_L(\eta^\prime, q^\prime) , \quad t \ge 0.
\label{eq:underline-psi}
\end{align}
For a vector $b =(b_1,\ldots,b_m)^\T$, $\| b \|_{\infty2}$ denotes $\max_{j\not= k \in [m]}$ $ (|b_j| + |b_k|)$.
\tzq{Inequality (\ref{eq:zo-regret-bd-proper}) can be seen to extend the two-class regret bound (\ref{eq:regret-BJM}) in Bartlett et al.~(2006) to multi-class settings for proper scoring rules.
Unlike the two-class setting, additional effort is needed to find a simple meaningful lower bound of $\underline\psi$ in the multi-class setting.}
Our current approach involves deriving a lower bound on the regret (or Bregman divergence) $B_L$ by the $L_1$ norm, such that for any $\eta, q\in\Delta_m$,
\begin{align}
B_L( \eta, q ) = H_L(q) - H_L(\eta) - (q-\eta)^\T \partial H_L(q) \ge \frac{\kappa_L}{2}  \| \eta - q\|_1 ^2 ,  \label{eq:bregman-bd}
\end{align}
where $\kappa_L >0$ is a constant depending on $L$, and $\|b\|_1= \sum_{j=1}^m |b_j|$ is the $L_1$ norm for any vector $b=(b_1, \ldots, b_m)^\T$.
Hence (\ref{eq:bregman-bd}) can be interpreted as saying that $-H_L$ is strongly convex with respect to the $L_1$ norm with modulus $\kappa_L$.
Because $\|\eta-q\|_1 \ge \|\eta-q\|_{\infty2}$, the regret bound (\ref{eq:zo-regret-bd-proper}) together with (\ref{eq:bregman-bd}) implies that for any $\eta,q\in\Delta_m$,
\begin{align}
\frac{\kappa_L}{2} \left( B^{\text{zo}}(\eta, q) \right)^2 \le  B_L (\eta, q) . \label{eq:regret-bd2-proper}
\end{align}
In general, the preceding discussion shows that a lower bound on the Bregman divergence $B_L(\eta,q)$ by some norm of $\eta-q$ can be translated into a corresponding regret bound.

Our current approach does not exploit the restriction that $\max_{j\in[m]} q_j^\prime \le 1/2$ in the definition of $\underline \psi$.
Hence it is interesting to study how our results here can be improved. On the other hand, such an improvement, even if achieved, may be limited.
See the later discussion on regret bounds for the pairwise exponential loss.

Our approach leads to the following result for two classes of proper scoring rules discussed in Section~\ref{sec:scoring-examples}:
a class of pairwise losses (\ref{eq:pairwise-s}) with $f_0$ associated with a Beta family of weight functions as studied in \cite{buja2005}, and
a class of simultaneous losses (\ref{eq:simul-loss}). In all these cases, inequalities (\ref{eq:bregman-bd}) can be of independent interest.

\begin{pro} \label{pro:bregman-bd}
Inequalities (\ref{eq:bregman-bd}) and (\ref{eq:regret-bd2-proper}) hold for the following proper scoring rules.
\begin{itemize}
\item[(i)] Consider a pairwise loss $L = L_{f_0}^{\text{pw,s}} $ in (\ref{eq:pairwise-s}),  with a univariate function $f_0$ defined such that
(\ref{eq:2class-grad}) holds with a weight function $w(q_1) =2^{2\nu} q_1^{\nu-1} q_2^{\nu-1}$
for $(q_1,q_2)^\T \in \Delta_2$.
If $\nu \le 0$, then (\ref{eq:bregman-bd}) and (\ref{eq:regret-bd2-proper}) are valid with $\kappa_L =2$.

\item[(ii)] Consider a simultaneous loss $L = L_\beta$ in (\ref{eq:simul-loss}).
Then (\ref{eq:bregman-bd}) and (\ref{eq:regret-bd2-proper}) are valid with
\begin{align*}
\kappa_L = \left\{ \begin{array}{cl}
(1-\beta) m^{(1-1/\beta)(2\beta-1)}2^{2-2\beta}, & \mbox{if } \beta \in [1/2,1), \\
(1-\beta)2^{1/\beta-1}, & \mbox{if } \beta \in (0,1/2].
\end{array}\right.
\end{align*}
The bounds from the two segments both give $\kappa_L = 1$ at $\beta=1/2$.
\end{itemize}
\end{pro}

We discuss several specific examples.
The standard likelihood loss $L(j,q) = -\log q_j$ is equivalent to the simultaneous loss $L_\beta$  in the limit of $\beta \to 1$ after properly rescaled.
In this case, Pinsker's inequality states that (\ref{eq:bregman-bd}) holds with $\kappa_L=1$ \citep[Lemma 12.6.1]{cover1991}:
\begin{align}
\sum_{j=1}^m \eta_j \log (\eta_j/ q_j) \ge \frac{1}{2} \left( \sum_{j=1}^m |\eta_j-q_j| \right)^2. \label{eq:bregman-bd-KL}
\end{align}
The resulting regret bound (\ref{eq:regret-bd2-proper}) for the standard likelihood loss $L$ then gives
\begin{align}
\frac{1}{2} \left( B^{\text{zo}}(\eta, q) \right)^2 \le  B_L (\eta, q) .  \label{eq:regret-likelihood}
\end{align}
This surrogate regret bound for the multinomial likelihood appears new, even though the Bregman divergence bound (\ref{eq:bregman-bd-KL}) is known.
\tzq{In the Supplement, we verify that (\ref{eq:bregman-bd-KL}) can be recovered from (\ref{eq:bregman-bd}) using Proposition~\ref{pro:bregman-bd}(ii).}

\tzq{The pairwise exponential loss associated with multi-class boosting} is defined equivalently as
$L_e^{\text{pw,s}} (j, q) = 2 (L_{1/2} -1) = 2 \sum_{k\in[m], k\not=j} \sqrt{q_k/q_j}$ in Section~\ref{sec:scoring}.
The two inequalities (\ref{eq:bregman-bd}) obtained from Proposition~\ref{pro:bregman-bd}, part (i) with $\nu=-1/2$ and part (ii) with $\beta = 1/2$,
are equivalent to each other and both lead to
\begin{align}
H_{L_{1/2}} (q) - H_{L_{1/2}}  (\eta) - (q-\eta)^\T \partial H_{L_{1/2}}  (q) \ge \frac{1}{2} \| \eta - q\|_1 ^2 , \label{eq:bregman-bd-exp}
\end{align}
where $H_{1/2}(q) = \|q\|_{1/2}$ and $L_{1/2}(j,q) = (\| q \|_{1/2} / q_j)^{1/2} = \sum_{k=1}^m \sqrt{q_k/q_j}$.
The resulting regret bound (\ref{eq:regret-bd2-proper}) for the  rescaled pairwise exponential loss $L_{1/2}$ gives
\begin{align}
\frac{1}{2} \left( B^{\text{zo}}(\eta, q) \right)^2 \le  B_{L_{1/2}} (\eta, q) .   \label{eq:regret-exponential}
\end{align}
The two bounds (\ref{eq:regret-likelihood}) and (\ref{eq:regret-exponential}) for the likelihood and rescaled pairwise exponential losses
happen to be of the same form, due to the scaling used.
For the two-class exponential loss defined as $L_e = L_{1/2}-1$, the existing regret bound (\ref{eq:regret-BJM}), corresponding
to an exact calculation of $\underline\psi$ by the proof of (\ref{eq:cw0-regret-bd-2class}) later, is
\begin{align*}
 1-\sqrt{1-(B^{\text{zo}}(\eta, q) )^2} \le B_{L_{1/2}} (\eta, q),
\end{align*}
which is slightly stronger than (\ref{eq:regret-exponential}) because $1-\sqrt{1-\delta^2} \ge \delta^2/2$ for $\delta\in [0,1]$,
but $(1-\sqrt{1-\delta^2})/ (\delta^2/2) \to 1$ as $\delta\to 0$. Therefore, our result (\ref{eq:regret-exponential}) provides a reasonable extension of
existing regret bounds to multi-class pairwise exponential losses.

A notable proper scoring rule which is not informed by Proposition~\ref{pro:bregman-bd} for $m \ge 3$ is the simultaneous exponential loss $L^{\text{r}}_0$ as used in \cite{zou2008},
even though the loss $L^{\text{r}}_0$ is equivalent to the exponential loss for $m=2$. \tzq{See the Supplement for details.}

\begin{rem} \label{rem:bregman-f}
\tzq{Inequality (\ref{eq:bregman-bd}) on the Bregman divergence in general differs from generalized Pinsker inequalities relating (two-distribution) $f$-divergences
to the total variation studied in \cite{reid2011}, Section 7.2, for binary experiments.}
For the pairwise exponential loss, the Bregman divergence on the left-hand side of (\ref{eq:bregman-bd-exp}) can be calculated as
$(\sum_{j\in[m]} \sqrt{q_j} )(\sum_{j\in[m]} \eta_j/\sqrt{q_j}) - (\sum_{j\in[m]} \sqrt{\eta_j} )^2$,
which is apparently not any $f$-divergence between probability vectors $\eta$ and $q$.
An exception is the classical Pinsker inequality (\ref{eq:bregman-bd-KL}): the Kullback--Liebler divergence on the left-hand side of (\ref{eq:bregman-bd-KL})
is both an $f$-divergence with $f(t)=t\log t$ and a Bregman divergence with $H_L(\eta) = -\sum_{j \in[m]} \eta_j \log \eta_j$.
\end{rem}

\begin{rem} \label{rem:AUC}
\tzq{In the two-class setting, a scoring rule satisfying inequality (\ref{eq:bregman-bd}) is called a strongly proper loss,
and surrogate regret bounds are obtained for strongly proper losses with respect to the area under the curve (AUC) in \cite{agarwal2013}.}
It is interesting to investigate possible extensions of such results to the multi-class setting.
\end{rem}

\subsection{Cost-transformed losses} \label{sec:cost-transform}

We study two types of classification regret bounds with respect to a general cost-weighted classification loss as defined in Section~\ref{sec:misclass}.
This subsection deals with the first type where a classification regret bound is derived for a loss, allowed to depend on a pre-specified cost matrix $C$,
similarly as the hinge-like loss $L^{\text{cw3}}$ in (\ref{eq:cw-loss3}).
An action of the loss is directly taken as a prediction for the cost-weighted classification loss.
See the next subsection on the second type of classification regret bounds.

For a general loss $L(j,\gamma)$ (not just scoring rules), define a cost-transformed loss, depending on a cost matrix $C$, as
\begin{align}
\tilde L (j,\gamma) = c_{jM} L (j,\gamma) + \sum_{k\in[m], k\not=j} (c_{jM} - c_{jk}) \{L(k,\gamma)-1\},  \label{eq:transform-loss}
\end{align}
where $c_{jM} = \max_{k\in[m]} c_{jk}$.
In the special case where $C = 1_m 1_m^\T - I_m$ for the zero-one loss, the transformed loss $\tilde L(j,\gamma)$ reduces to the original loss $L(j,\gamma)$.
A motivation for this construction is that the cost-weighted classification loss can also be obtained in this way from the zero-one loss:
$ L^{\text{cw}}(j, \gamma) = \tilde L^{\text{zo}}(j,\gamma)$.
In general, the risk and regret of the transformed loss can be related to those of the original loss as follows.

\begin{lem} \label{lem:scaling}
The risks of the losses $\tilde L(j,\gamma)$ and $L(j,\gamma)$ satisfy
\begin{align*}
R_{\tilde L} (\eta, \gamma) = (1_m^\T \tilde \eta) \, R_L ( \tilde{\tilde \eta}, \gamma) -D(\eta),
\end{align*}
where $D(\eta) =  \sum_{j\in[m]} \sum_{k\in[m], k\not=j} \eta_j (c_{jM} - c_{jk})$, $\tilde {\tilde \eta} = \tilde \eta / (1_m^\T \tilde\eta) \in \Delta_m $, and
$\tilde \eta = (\tilde\eta_1,\ldots, \tilde\eta_m)^\T$ $ \in \bbR^m_+$ with
\begin{align*}
\tilde \eta_j = c_{jM} \eta_j + \sum_{k\in[m], k\not=j} (c_{kM} - c_{kj}) \eta_k.
\end{align*}
Moreover, the regrets of $\tilde L$ and $L$ satisfy $B_{\tilde L} ( \eta, \gamma) =  (1_m^\T \tilde \eta) \, B_L (\tilde{\tilde\eta}, \gamma)$.
\end{lem}

For a scoring rule $L(j,q)$ with actions defined as probability vectors $q\in\Delta_m$, there is a simple upper bound on the regret of the associated zero-one loss $L^{\text{zo}}(j, q)$,
which is instrumental to our derivation of classification regret bounds.

\begin{lem} \label{lem:misclass}
For any $\eta, q\in \Delta_m$, it holds that
\begin{align*}
B^{\text{zo}}(\eta, q) \le \| \eta- q \|_{\infty2},
\end{align*}
where $\| b \|_{\infty2} = \max_{j\not= k \in [m]}$ $ (|b_j| + |b_k|)$ for any vector $b =(b_1,\ldots,b_m)^\T$.
The bound is tight for any $m\ge 2$ in that there exist $\eta,q \in \Delta_m$ for which the bound becomes exact.
\end{lem}

Combining the preceding two lemmas and invoking a similar argument as indicated by (\ref{eq:regret-general}) leads to the following regret bound, depending on the action $q$.

\begin{pro} \label{pro:regret-bd}
For a scoring rule $L(j,q)$, define a nondecreasing function $\psi_q$:
\begin{align*}
\psi_q(t) = \inf_{\eta^\prime \in\Delta_m: \|\eta^\prime -q \|_{\infty2}\ge t} B_L(\eta^\prime, q) , \quad t \ge 0.
\end{align*}
Then the regrets of the cost-weighted classification loss $L^{\text{cw}}(j,q)$ and the cost-transformed scoring rule $\tilde L(j,q)$ satisfy
\begin{align}
\psi_q \left( \frac{B^{\text{cw}}(\eta, q)}{1_m^\T \tilde\eta } \right) \le \frac{B_{\tilde L} (\eta, q)}{1_m^\T \tilde\eta} , \label{eq:cw-regret-bd}
\end{align}
where $B^{\text{cw}}= B_{L^{\text{cw}}}$, and $\tilde\eta $ is defined, depending on $\eta$ and $C$,  as in Lemma~\ref{lem:scaling}.
\end{pro}

A cost-transformed loss (\ref{eq:transform-loss}) from a proper scoring rule can be easily shown to remain a proper scoring rule.
In this case, a uniform regret bound can be obtained from (\ref{eq:cw-regret-bd}), by taking an infimum over $q$ and incorporating simplification
due to the representation of the regret (\ref{eq:bregman}) as a Bregman divergence for a proper scoring rule.

\begin{cor} \label{cor:regret-bd-proper}
For a proper scoring rule $L(j,q)$, the regrets of $L^{\text{cw}}(j, q)$ and $\tilde L(j,q)$ satisfy
\begin{align}
\underline \psi \left( \frac{B^{\text{cw}}(\eta, q)}{1_m^\T \tilde\eta } \right) \le \frac{B_{\tilde L} (\eta, q)}{1_m^\T \tilde\eta} , \label{eq:cw-regret-bd-proper}
\end{align}
where $\underline \psi $ is defined in (\ref{eq:underline-psi}), and
$\tilde\eta $ is defined, depending on $\eta$ and $C$, in Lemma~\ref{lem:scaling}.
\end{cor}

It is instructive to examine the regret bound (\ref{eq:cw-regret-bd-proper}) in the special case of class-weighted costs, where $C = C_0 1_m^\T - \diag(C_0)$ with $C_0=(c_{10}, \ldots, c_{m0})^\T$.
The cost-transformed loss $\tilde L$ reduces to $\tilde L(j,q) = c_{j0} L(j,q)$.
The regret bound (\ref{eq:cw-regret-bd-proper}) becomes
\begin{align}
\underline\psi \left( \frac{B^{\text{cw0}}(\eta, q)}{C_0^\T\eta} \right) \le \frac{B_{\tilde L} (\eta, q)}{C_0^\T\eta} , \label{eq:cw0-regret-bd-proper}
\end{align}
where $B^{\text{cw0}}= B_{L^{\text{cw0}}}$. We defer a discussion of these results until after Corollary~\ref{cor:regret-bd2-proper}.

\subsection{Cost-independent losses} \label{sec:cost-indep}

We derive a different type of classification regret bounds than in the preceding subsection.
Here a loss used for training is defined independently of any cost matrix, but an action of the loss can be converted after training to a prediction, depending on the
cost matrix $C$, for the cost-weighted classification loss.
For scoring rules, our derivation relies on the following extension of Lemma~\ref{lem:misclass} on the regret of the cost-weighted classification loss,
where a prediction is linearly converted from a probability vector.

\begin{lem} \label{lem:misclass2}
For any $\eta, q\in \Delta_m$, it holds that
\begin{align*}
B^{\text{cw}}(\eta, \overline C^\T q) \le \| \overline C^\T (\eta- q) \|_{\infty2},
\end{align*}
where $\overline C = C_M 1_m^\T - C$ and $C_M=(c_{1M}, \ldots, c_{mM})^\T$ with $c_{jM} = \max_{k\in[m]} c_{jk}$ for $j\in[m]$ as
defined in the transformed loss (\ref{eq:transform-loss}).
\end{lem}

By a similar argument as indicated by (\ref{eq:regret-general}), we obtain a regret bound which compares the regret of
a scoring rule $L(j,q)$ with that of the cost-weighted classification loss with a prediction depending on both $q$ and $C$ as in Lemma~\ref{lem:misclass2}.

\begin{pro} \label{pro:regret-bd2}
For a scoring rule $L(j,q)$, define a nondecreasing function $\psi_q^C $:
\begin{align*}
\psi_q^C (t) = \inf_{\eta^\prime \in\Delta_m: \| \overline C^\T (\eta^\prime -q) \|_{\infty2}\ge t} B_L(\eta^\prime, q) , \quad t \ge 0.
\end{align*}
Then  the regrets of the cost-weighted classification loss $L^{\text{cw}}(j, \overline C^\T q)$ and the scoring rule $L(j,q)$ satisfy
\begin{align}
\psi_q^C \left( B^{\text{cw}}(\eta, \overline C^\T q) \right) \le B_L (\eta, q) . \label{eq:cw-regret-bd2}
\end{align}
\end{pro}

For a {\it proper} scoring rule $L(j,q)$, the regret bound (\ref{eq:cw-regret-bd2}) can be strengthened (see the Supplement for a proof) such that for each $w \in \mathcal W_{\eta,q}$,
\begin{align}
\psi^C_{q^w}  \left( B^{\text{cw}}(\eta, \overline C^\T q) \right) \le  B_L (\eta, q), \label{eq:cw-regret-bd3}
\end{align}
where $q^w=(1-w)\eta+wq$ and
\begin{align*}
\mathcal W_{\eta,q} = \left\{w \in [0,1]: \overline C^\T_k q^w = \max_j (\overline C^\T_j q^w) \mbox{ for } k=\argmax_j (\overline C^\T_j q) \right\} \ni 1.
\end{align*}
By definition, $w\in \mathcal W_{\eta,q}$ means that using $q^w$ yields the same classification as using $q$.
Moreover, a uniform regret bound can be obtained from (\ref{eq:cw-regret-bd3}) by minimizing over $q^w$ with $w \in \mathcal W_{\eta,q}$ such that
$\max_{j \in [m]} \overline C^\T_j q^w \le 1_m^\T \overline C^\T  q^w /2$.

\begin{cor} \label{cor:regret-bd2-proper}
For a proper scoring rule $L(j,q)$, define
\begin{align*}
\underline \psi^C (t) = \inf_{\substack{\eta^\prime, q^\prime\in \Delta_m:  \, \|\overline C^\T (\eta^\prime -q^\prime) \|_{\infty2} = t,\\
\max_{j\in[m]} (\overline C_j^\T q^\prime) \le 1_m^\T \overline C^\T q^\prime /2}}
B_L(\eta^\prime, q^\prime) , \quad t \ge 0.
\end{align*}
Then  the regrets of $L^{\text{cw}}(j, \overline C^\T q)$ and $L(j,q)$ satisfy
\begin{align}
\underline \psi^C \left( B^{\text{cw}}(\eta, \overline C^\T q) \right) \le  B_L (\eta, q) . \label{eq:cw-regret-bd2-proper}
\end{align}
\end{cor}

In the special case of class-weighted costs, corresponding to $C = C_0 1_m^\T - \diag(C_0)$ with $C_0=(c_{10}, \ldots, c_{m0})^\T$, define
\begin{align*}
\underline \psi^{C_0} (t) = \inf_{\substack{\eta^\prime, q^\prime\in \Delta_m:  \, \|C_0\circ(\eta^\prime -q^\prime) \|_{\infty2} = t,\\
\max_{j\in[m]} (c_{j0}q^\prime_j) \le C_0^\T q^\prime/2 }}
B_L(\eta^\prime, q^\prime) , \quad t \ge 0,
\end{align*}
where $\circ$ denotes the component-wise product between two vectors.
The regret bound (\ref{eq:cw-regret-bd2-proper}) for proper scoring rules reduces to
\begin{align}
\underline \psi^{C_0} \left( B^{\text{cw0}}(\eta, C_0\circ q) \right) \le  B_L (\eta, q) . \label{eq:cw0-regret-bd2-proper}
\end{align}
It is interesting to compare the two regret bounds (\ref{eq:cw0-regret-bd-proper}) and (\ref{eq:cw0-regret-bd2-proper}).
On one hand, for the zero-one loss with $C_0= 1_m$, both of these bounds lead to the regret bound (\ref{eq:zo-regret-bd-proper}) discussed in Section~\ref{sec:zero-one}.
On the other hand, the two bounds (\ref{eq:cw0-regret-bd-proper}) and (\ref{eq:cw0-regret-bd2-proper}) in general serve different purposes.
%let alone the fact that the function $\underline \psi$ in (\ref{eq:cw0-regret-bd-proper}) is independent of $C_0$, whereas
%$\underline \psi^{C_0}$ in (\ref{eq:cw0-regret-bd2-proper}) depends on $C_0$.
The bound (\ref{eq:cw0-regret-bd-proper}) compares the regrets of the transformed scoring rule $\tilde L$ depending on $C_0$ and
the classification loss $L^{\text{cw0}}$ with the prediction always set to $q$.
To use $\tilde L$, a different round of training is required for a different choice of $C_0$.
The bound (\ref{eq:cw0-regret-bd2-proper}) relates the regrets of the original scoring rule $L$, independent of $C_0$, and
the classification loss $L^{\text{cw0}}$ with the prediction defined as $C_0 \circ q$.
Only one round of training is needed to determine $q$ when using $L$, and then the prediction can be adjusted from $q$ according to
the choice of $C_0$. Hence the bound (\ref{eq:cw0-regret-bd2-proper}) can be potentially more useful than (\ref{eq:cw0-regret-bd-proper}).

For binary classification with $m=2$, the regret bound (\ref{eq:cw0-regret-bd2-proper}) for proper scoring rules can be shown to
recover Theorem 25 in \cite{reid2011}. For a {\it proper} scoring rule $L(j,q)$ and any $\eta,q\in \Delta_2$, it holds that
\begin{align}
\min \left\{ \psi^{\text{RW}}  (\delta), \psi^{\text{RW}} (-\delta) \right\} \le B_L (\eta, q),\label{eq:cw0-regret-bd-2class}
\end{align}
where $\delta = B^{\text{cw0}}(\eta, C_0\circ q) $, $\psi^{\text{RW}}  (\delta) = B_L^1 ((c_{20}+\delta)/(c_{10}+c_{20}), c_{20}/(c_{10}+c_{20}))$,
and $B_L^1(\eta_1,q_1) = B_L(\eta,q)$ with $\eta=(\eta_1,\eta_2)^\T$ and $q=(q_1,q_2)^\T$,
that is, $B_L^1(\eta_1,q_1)$ is $B_L(\eta, q)$  treated as a function of $(\eta_1,q_1)$ only.
See the Supplement for details.

\section{Conclusion}

In this article, we are mainly concerned with constructing losses and establishing corresponding regret bounds in multi-class settings.
Various topics are of interest for further research.
Large sample theory can be studied regarding estimation and approximation errors, similarly as in \cite{zhang2004b} and \cite{bartlett2006}, by taking advantage of our multi-class regret bounds.
Computational algorithms need to be developed for implementing our new hinge-like losses and, in connection with boosting algorithms, for implementing composite losses based on new proper scoring rules.
Numerical experiments are also desired to evaluate empirical performance of new methods.

\vskip 0.2in

\bibliographystyle{ztan}  % already in .sty file

\bibpunct{(}{)}{;}{a}{,}{,}
\bibliography{loss-rev}

%%%%%%%%%%%%%%%%%%%%%%%%%%%%%%%%%%%%%%%%%%%%%%%%%%%%%%%%%%%%%%%%%%%%%%%%%%%%%%%%%%%%%%%%%%

\clearpage

\setcounter{page}{1}

\setcounter{section}{0}
\setcounter{equation}{0}

\setcounter{figure}{0}
\setcounter{table}{0}

\renewcommand{\theequation}{S\arabic{equation}}
\renewcommand{\thesection}{\Roman{section}}

\renewcommand\thefigure{S\arabic{figure}}
\renewcommand\thetable{S\arabic{table}}

\setcounter{lem}{0}
\renewcommand{\thelem}{S\arabic{lem}}

\begin{center}
{\Large Supplementary Material for}

{\Large ``On Loss Functions and Regret Bounds for Multi-category Classification"}

\vspace{.1in} {\large Zhiqiang Tan and Xinwei Zhang}
\end{center}

\section{Technical details}\label{app:proof}
%In this appendix, we provide remaining proofs for lemmas and propositions in the main text. For completeness, we also provide some background facts which are used in the proof.

\subsection{Preparation}

For a convex function $\psi$ defined on a convex domain $\Omega$, the Bregman divergence is defined as
\begin{align*}
B_\psi(x,y ) = \psi(x) - \psi(y) - (x-y) ^\T \partial \psi (y) ,
\end{align*}
where $\partial \psi$ is a sub-gradient of $\psi$.
The symmetrized Bregman divergence is
\begin{align*}
B_\psi(x,y ) + B_\psi(y,x) = (y-x)^\T \{\partial \psi (y)  - \partial \psi(x)\}.
\end{align*}
The following lemma shows that the Bregman divergence is nondecreasing as the first (or second) argument, $x$ or $y$, moves away from the other argument, $y$ or $x$, along a straight line, while  the second (or respectively first) argument remains fixed.

\begin{lem} \label{lem:breg-mon}
For any $x,y\in \Omega$ and $w \in [0,1]$, we have
\begin{align}
B_\psi(x,y ) \ge B_\psi(x^w, y), \label{eq:breg-mon1}\\
B_\psi(x,y) \ge B_\psi(x, x^w). \label{eq:breg-mon2}
\end{align}
where $x^w = (1-w)x+wy$.
\end{lem}

\begin{prf}
If $w=0$ or 1, then (\ref{eq:breg-mon1}) and (\ref{eq:breg-mon2}) hold trivially. In the following, assume $w \in (0,1)$.
To show (\ref{eq:breg-mon1}), direct calculation yields
\begin{align*}
& B_\psi(x,y ) - B_\psi(x^w, y) = \psi(x) - \psi(x^w) - (x-x^w)^\T \partial \psi (y)\\
& = B_\psi(x, x^w) + (x-x^w)^\T \{ \partial \psi(x^w) - \partial \psi (y) \} \\
& = B_\psi(x, x^w) + \frac{w}{1-w} (x^w -y)^\T \{ \partial \psi(x^w) - \partial \psi (y) \} .
\end{align*}
Hence (\ref{eq:breg-mon1}) follows because $(x^w -y)^\T \{ \partial \psi(x^w) - \partial \psi (y) \}$ is the symmetrized Bregman divergence between $x^w$ and $y$.
From the preceding equations, we see
\begin{align*}
& B_\psi(x,y ) - B_\psi(x, x^w) %= \psi(x^w) - \psi(y) - (x-y)^\T \partial \psi (y) + (x-x^w)^\T \partial \psi(x^w)\\
%& = B_\psi(x^w, y) + (x-x^w)^\T \{ \partial \psi(x^w) - \partial \psi (y) \} \\
 = B_\psi(x^w,y) + \frac{w}{1-w} (x^w -y)^\T \{ \partial \psi(x^w) - \partial \psi (y) \} .
\end{align*}
Hence (\ref{eq:breg-mon2}) follows because $(x^w -y)^\T \{ \partial \psi(x^w) - \partial \psi (y) \} \ge 0$ again.
\end{prf}

\subsection{Proofs of results in Section~\ref{sec:construct-loss}}

\noindent\textbf{Proof of Proposition~\ref{pro:f-loss2}.}\;
Denote by $\partial^\dag f$ the set of all sub-gradients of $f$. For any $u \in \overline\bbR_+^{m-1}$ and
$s = \partial f(u) \in \partial^\dag f(u)$, Fenchel's conjugacy property implies that
$f^*(s) = u^\T s -f(u)$ and hence $s \in \dom(f^*)$. Moreover, we have
\begin{align*}
& \sum_{j=1}^m \eta_j L_{f2}(j,u) =\sum_{j=1}^{m-1}  \eta_j(-\partial_j f(u)) + \eta_m (u^\T \partial f(u) - f(u))\\
& \quad= \sum_{j=1}^{m-1} (-\eta_j s_j) + \eta_m f^*(s) = \sum_{j=1}^m \eta_j L_f(j, s) .
\end{align*}
Therefore,
\begin{align*}
\inf_{u \in \overline\bbR_+^{m-1}} \left\{ \sum_{j=1}^m \eta_j L_{f2}(j,u) \right\}
\ge \inf_{s \in \dom(f^*)}  \left\{ \sum_{j=1}^m \eta_j L_f(j, s) \right\} = H_f(\eta).
\end{align*}
Next, we show the reverse inequality.
For any $\eta \in \Delta_m$, denote $u^\eta = (\eta_1 / \eta_m, \ldots, \eta_{m-1} / \eta_m)^\T$ and $s^\eta = \partial f(u^\eta) \in\partial^\dag f(u^\eta)$.
Then
\begin{align*}
& H_f (\eta) = - \eta_m f( u^\eta) = -\eta_m \left\{ \sum_{j=1}^{m-1} u_j^\eta s_j^\eta - f^*(s^\eta) \right\}\\
& =\sum_{j=1}^{m-1} (-\eta_j s_j^\eta) + \eta_m f^*(s^\eta) =\sum_{j=1}^{m-1}  \eta_j(-\partial_j f(u^\eta)) + \eta_m (u^{\eta\T} \partial f(u^\eta) - f(u^\eta))  ,
\end{align*}
where Fenchel's conjugacy property, $u^{\eta\T} s^\eta = f(u^\eta) + f^*(s^\eta)$, is used in the last equalities on the first and second lines. Hence
$ H_f(\eta) \ge \inf_{u \in \overline\bbR_+^{m-1}} \{ \sum_{j=1}^m \eta_j L_{f2}(j,u) \}$. \hfill$\blacksquare$\\[2mm]

\noindent\textbf{Proof of equation~(\ref{eq:savage2}).}\;
By definition, $H_f(q) = - q_m f( q_1/q_m, \ldots, q_{m-1}/q_m)$.
The sub-gradient of $-H_f$ can be calculated as
\begin{align*}
-\partial_j H_f(q_1, \ldots, q_m) = \left\{
\begin{array}{cl}
\partial_j f(u^q), & \mbox{if } j \in [m-1], \\
f(u^q) - \sum_{j=1}^{m-1} \frac{q_j}{q_m}  \partial_j f(u^q) , & \mbox{if } j=m,
\end{array} \right.
\end{align*}
where $u^q = ( q_1/q_m, \ldots, q_{m-1}/q_m)^\T$. Substituting these expressions into $H_f( q) - \sum_{j=1}^m  (q_j - \eta_j) \partial_j H_f( q)$ yields the second equality in Eq.~(\ref{eq:savage2}):
\begin{align*}
& H_f( q) - \sum_{j=1}^m  (q_j - \eta_j) \partial_j H_f( q) \\
& = -q_m f(u^q) + \sum_{j=1}^{m-1} (q_j - \eta_j) \partial_j f(u^q) + (q_m-\eta_m) \left\{ f(u^q) - \sum_{j=1}^{m-1} \frac{q_j}{q_m}  \partial_j f(u^q) \right\}\\
& = -\sum_{j=1}^{m-1}\eta_j \partial_j f(u^q) + \eta_m \left\{ -f(u^q) + \sum_{j=1}^{m-1} \frac{q_j}{q_m}  \partial_j f(u^q) \right\} .
\end{align*}
\hfill$\blacksquare$\\[2mm]

\subsection{Proofs of results in Section~\ref{sec:scoring-examples}}

\noindent\textbf{Proof of equation~(\ref{eq:pairwise-s}).}\; By manipulating the summation, we have
\begin{align*}
 & L_{f_0}^{\text{pw,s}}(j, q) = \sum_{l,k \in [m], k\not= l} \left[- \one_k(j) \partial  f_0( \frac{q_k}{q_l} ) +
 \one_l(j) \left\{ \frac{q_k}{q_l} \partial f_0( \frac{q_k}{q_l}) - f_0( \frac{q_k}{q_l}) \right\} \right] \\
 & = \sum_{l\in [m]} \sum_{k\in [m], k\not=l} \left\{ - \one_k(j) \partial  f_0( \frac{q_k}{q_l} ) \right\} +
  \sum_{k\in [m]} \sum_{l\in [m], l \not=k } \one_l(j) \left\{ \frac{q_k}{q_l} \partial f_0( \frac{q_k}{q_l}) - f_0( \frac{q_k}{q_l}) \right\} \\
 & = \sum_{l\in [m], j\not=l} \left\{ - \partial  f_0( \frac{q_j}{q_l} ) \right\} +
  \sum_{k\in [m], j \not=k} \left\{ \frac{q_k}{q_j} \partial f_0( \frac{q_k}{q_j}) - f_0( \frac{q_k}{q_j}) \right\},
\end{align*}
which yields the desired result.
\hfill$\blacksquare$\\[2mm]

\noindent\textbf{Convexity of two-class composite losses.}\;
Consider a logistic link $q^{h_0} =(q_1^{h_0}, q_2^{h_0})^\T$, where $q_1^{h_0} = \{1+\exp(-h_0)\}^{-1}$ or equivalently $q_1^{h_0}/q_2^{h_0} = \exp(h_0)$.
Then it can be easily shown that the three composite losses, $L_\ell(j, q^{h_0})$, $L_e (j, q^{h_0})$, and $L_c (j, q^{h_0})$, are convex in $h_0$, with the following gradients:
\begin{align*}
& \frac{\dif}{\dif h_0} L_\ell( j, q^{h_0}) = -\left\{\one_1(j) - q_1^{h_0} \right\}  ,\\
& \frac{\dif}{\dif h_0} L_e( j, q^{h_0}) =  -\left\{\one_1(j) -q_1^{h_0} \right\} (q_2^{h_0} q_1^{h_0}) ^{-1/2}  ,\\
& \frac{\dif}{\dif h_0} L_c( j, q^{h_0}) = -\left\{ \one_1(j) / q_1^{h_0} -1 \right\} /2 .
\end{align*}
The corresponding weight functions are in the Beta family, $w(q_1) =2^{\nu_1+\nu_2} q_1^{\nu_1-1} q_2^{\nu_2-1}$, with $(\nu_1,\nu_2)= (0,0), (-1/2,-1/2)$ and $(-1,0)$ respectively.
As discussed in \citet[Section 15]{buja2005}, a proper scoring rule with a logistic link and $w(q_1)$ in the Beta family is convex in $h_0$ if and only if $\nu_1,\nu_2 \in [-1,0]$.
Hence the likelihood and calibration losses are at the boundary of achieving convexity in $h_0$.
\hfill$\blacksquare$\\[2mm]

\noindent\textbf{Proof of Proposition~\ref{pro:limit}.}\;
The scoring rules are obtained directly from Proposition 3. First, we show the three limits of $H_\beta$ for $\beta=0,1,\infty$.

(\textit{i}) Rewrite $H_\beta(q)$ as
\begin{align*}
  H_\beta(q) &=  \frac{ \exp\{\frac{1}{\beta}\log(1+\frac{\sum_{j=1}^m (q_j^\beta-1)}{m})\} - m^{-\frac{1}{\beta}}}{m^{-1}-m^{-\frac{1}{\beta}}}.
\end{align*}
Using $\log(1+x) /x \to 1$ as $x\rightarrow 0$, we have
\begin{align*}
\lim_{\beta\to 0+}  H_\beta(q) &= \lim_{\beta\to 0+}\frac{ \exp\{\frac{\sum_{j=1}^m (q_j^\beta-1)}{\beta m}\} - m^{-\frac{1}{\beta}}}{m^{-1}-m^{-\frac{1}{\beta}}}
=m \left(\prod_{j=1}^m q_j\right)^{\frac{1}{m}},
\end{align*}
where the last step holds because  $\lim_{\beta\rightarrow 0+} (q_j^\beta -1)/\beta = \log q_j$ by L'Hopital's rule.

(\textit{iii}) Rewrite $H_\beta(q)$ as
\begin{align*}
    H_\beta(q) &=  \frac{ \exp\{\frac{1}{\beta}\log(\sum_{j=1}^m q_j^\beta)\} - 1}{\exp\{(\frac{1}{\beta}-1)\log m\}-1}.
\end{align*}
Using $(\me^x -1)/x \to 1$ as $x\rightarrow 0$, we obtain
\begin{align*}
   \lim_{\beta\to 1} H_\beta(q) &=  \lim_{\beta\to 1} \frac{ \log(\sum_{j=1}^m q_j^\beta) }{(1-\beta)\log m}.
\end{align*}
Applying L'Hopital's rule yields
\begin{align*}
   \lim_{\beta\rightarrow 1} H_\beta(q) =  \lim_{\beta\rightarrow 1} \frac{-\sum_{j=1}^mq_j^\beta \log q_j}{(\log m)(q_1^\beta+\dots+q_m^\beta  )} = \frac{ -\sum_{j=1}^m q_j\log q_j}{\log m}.
\end{align*}

(\textit{iv}) The result follows from the standard limit of $L^p$-norm, $\lim_{p\rightarrow \infty} \|x\|_p = \|x\|_\infty$,
where $\|x\|_p= (\sum_{j=1}^m |x_j|^p)^{1/p}$ and $\|x\|_\infty=\max_{ j\in [m]}|x_j|$ for $x\in\bbR^{m}$.

Finally, we show that the composite loss $L^{\text{r}}_\beta(j, q^h)$ is convex in $h$ for $\beta \in [0,1]$. The case $\beta=0$ or 1
can be verified directly, corresponding to the simultaneous exponential or likelihood composite loss.
For $\beta\in (0,1)$, the unscaled composite loss $L_\beta(j, q^h)$ is
\begin{align*}
& L_\beta(j, q^h) = \left\{1 + \sum_{i\neq j} \exp(\beta(h_i-h_j)) \right\}^{\frac{1}{\beta}-1}.
\end{align*}
It suffices to show that for $\beta\in(0,1)$, the function
\begin{align*}
g(x) = \left\{1 + \sum_{i=1}^{m-1} \exp(\beta x_i) \right\}^{\frac{1}{\beta}-1}
\end{align*}
is convex in $x \in\mathbb{R}^{m-1}$. Rewrite $g(x)$ as
\begin{align*}
g(x) = \exp \left[ \left(\frac{1}{\beta}-1 \right) \log \left\{1 + \sum_{i=1}^{m-1} \exp(\beta x_i ) \right\} \right].
\end{align*}
Note that $\log \{1 + \sum_{i=1}^{m-1} \exp(\beta x_i )\}$ is convex in $x$ \citep[Example~3.14]{boyd2004}.
The convexity of $g(x)$ follows by the scalar composition rule in \citet[Section~3.2.4]{boyd2004}.
\hfill$\blacksquare$\\[2mm]

\subsection{Proofs of results in Sections~\ref{sec:hinge-construction}--\ref{sec:hinge-regret-bd}}

\noindent\textbf{Proof of Lemma~\ref{lem:cw-conj}.}\;
Note that $f^{\text{cw}}(t) = \max_{k \in[m]} \, (-C_k^\T \tilde t)$, that is, the maximum of $m$ functions $-C_1^\T \tilde t, ..., -C_m^\T \tilde t$.
By a direct extension of Eq.~(1) in \cite{bot2008} to allow multiple functions, we have
\begin{align*}
f^{\text{cw}*} (s)  = \min_{\lambda \in \Delta_m} f_\lambda^*(s) ,
\end{align*}
where $f_\lambda = -(C\lambda)^\T \tilde t$. For each $\lambda\in\Delta_m$, direct calculation yields
\begin{align*}
f_\lambda^*(s) = \sup_{t \in \bbR_+^{m-1}} \left\{ st + (C\lambda)^\T \tilde t \right\} = \left\{
\begin{array}{cl}
(C\lambda)_m , & \mbox{if } s_j \le - (C\lambda)_j, j\in [m-1],\\
\infty, & \mbox{otherwise}.
\end{array} \right.
\end{align*}
The desired result then follows. \hfill$\blacksquare$\\[2mm]

\noindent\textbf{Proof of Lemma~\ref{lem:cw-conj2}.}\;
We need to show that for $\eta\in \Delta_m$,
\begin{align*}
H^{\text{cw}}  (\eta) = \inf_{\lambda \in \Delta_m } \left\{\sum_{j=1}^{m-1} \eta_j (C \lambda)_j + \eta_m (C\lambda)_m  \right\}.
\end{align*}
Although this can be directly established, we give a proof based on Proposition~\ref{pro:f-loss}. In fact, applying Proposition~\ref{pro:f-loss} with $f= f^{\text{cw}} $ yields
\begin{align*}
H^{\text{cw}}  (\eta) = \inf_{s\in \dom(f^{\text{cw}*} ) } \left\{ \sum_{j=1}^{m-1} \eta_j (-s_j) + \eta_m f^{\text{cw}*} (s) \right\}
\end{align*}
For each $s \in \dom(f^{\text{cw}*}) $, there exists some $\lambda^s\in \Delta_m$ such that $s_j \le - (C\lambda^s)_j, j\in [m-1] $ and hence by Lemma~\ref{lem:cw-conj},
\begin{align*}
& \sum_{j=1}^{m-1} \eta_j (-s_j) + \eta_m f^{\text{cw}*}(s)
 \ge \sum_{j=1}^{m-1} \eta_j (C\lambda^s)_j  + \eta_m (C\lambda^s)_m.
\end{align*}
Therefore,
\begin{align*}
H^{\text{cw}}  (\eta) \ge \inf_{\lambda \in \Delta_m } \left\{\sum_{j=1}^{m-1} \eta_j (C \lambda)_j + \eta_m (C\lambda)_m  \right\}.
\end{align*}
The reverse inequality can be obtained by using the fact that
for each $\lambda \in \Delta_m$,  the vector $s^\lambda$ is contained in $\dom(f^{\text{cw}*} ) $ with $s^\lambda_j = -(C\lambda)_j$.
\hfill$\blacksquare$\\[2mm]

\noindent\textbf{Proof of Proposition~\ref{pro:cw-conj3}.}
We need to show that for $\eta\in\Delta_m$,
\begin{align}
H^{\text{cw}} (\eta) = \inf_{\tau \in \bbR^{m-1}} \left\{\sum_{j=1}^m \eta_j L^{\text{cw3}} (j,\tau) \right\}. \label{eq:prf-cw-H}
\end{align}
In fact, Lemma~\ref{lem:cw-conj2} implies that for $\eta\in\Delta_m$,
\begin{align*}
H^{\text{cw}}(\eta) = \inf_{\lambda \in \Delta_m } \left\{\sum_{j=1}^m \eta_j L^{\text{cw2}} (j,\lambda) \right\},
\end{align*}
where by definition
\begin{align}
L^{\text{cw2}} (j,\lambda) =\left\{
\begin{array}{cl}
(C\lambda)_j = c_{jm} \lambda_m + \sum_{k\in[m-1],k\not=j} c_{jk} \lambda_k ,  & \mbox{if } j \in [m-1],\\
(C\lambda)_m = \sum_{k\in[m-1]} c_{mk} \lambda_k, & \mbox{if } j=m,
\end{array} \right. \label{eq:prf-cw2-loss}
\end{align}
It suffices to show that
\begin{itemize}
\item[(i)] $L^{\text{cw3}}$ is an extension of $L^{\text{cw2}}$ from $\Delta_m$ to $\bbR^{m-1}$, and

\item[(ii)] the minimum in (\ref{eq:prf-cw-H}) is achieved at $\tau \in \bbR^{m-1}$ such that $\tilde\tau\in\Delta_m$, where $\tilde\tau= (\tau_1, \ldots, \tau_{m-1}, 1-\sum_{k=1}^{m-1} \tau_k)^\T$.
\end{itemize}
For the extension in (i), $L^{\text{cw2}}(j,\lambda)$ is considered a function of $j$ and $(\lambda_1,\ldots,\lambda_{m-1})^\T$, with $\lambda_m =1- \sum_{j=1}^{m-1} \lambda_j$,
such that $\lambda\in \Delta_m$.

Result (i) is immediate by comparison of (\ref{eq:cw-loss3}) with (\ref{eq:prf-cw2-loss}). For any $\tau \in \bbR^{m-1}$ such that $\tilde\tau\in\Delta_m$, we have
$\tau_{k+} = \tau_k$ for $k \in [m-1]$, $\tau^{(j)}_{m+} =  \tau^{(j)}_m = 1- \sum_{k\in[m-1]} \tau_k$ for any $j \in[m-1]$,
and hence $L^{\text{cw3}} (j, \tau) = L^{\text{cw2}}(j, \tilde\tau)$ for $j \in [m-1]$ or $j=m$.

For result (ii), we distinguish two cases. First, we show that for any $\tau\in \bbR^{m-1}$ with one or more negative components and $j \in [m]$,
\begin{align}
L^{\text{cw3}} (j,\tau^\prime) \le  L^{\text{cw3}} (j,\tau) , \label{eq:prf-cw-conj3-ineq1}
\end{align}
where $\tau^\prime$ is obtained from $\tau$ by setting all negative components of $\tau$ to 0.
In fact, by examining (\ref{eq:cw-loss3}), we have
$L^{\text{cw3}} (m,\tau^\prime)= L^{\text{cw3}} (m,\tau)$, because $\tau^\prime_{k+} = \tau_{k+}$ for each $k\in[m-1]$.
Moreover, $L^{\text{cw3}} (j,\tau^\prime) \le  L^{\text{cw3}} (j,\tau)$, by noting that $\tau^\prime_j \ge \tau_j$ and $(\tau^\prime)^{(j)}_m \le \tau^{(j)}_m $ for $j\in[m-1]$,
where $(\tau^\prime)^{(j)}_m $ is defined by (\ref{eq:tau-m}) with $\tau$ replaced by $\tau^\prime$.

Second, we show that for any $\tau\in \bbR_+^{m-1}$ (i.e., all components of $\tau$ are nonnegative) with $\sum_{k=1}^{m-1} \tau_k >1$ and $j \in [m]$,
\begin{align}
L^{\text{cw3}} (j, \tau^\dprime) \le  L^{\text{cw3}} (j,\tau) , \label{eq:prf-cw-conj3-ineq2}
\end{align}
where $\tau^\dprime = (\tau^\dprime_1, \ldots, \tau^\dprime_{m-1})^\T \in \bbR_+^{m-1}$ with $\tau^\dprime_k = (\tau_k - b)_+$
and $b>0$ chosen such that $\sum_{k=1}^{m-1} \tau^\dprime_k = 1$.
This choice of $b$ exists, because $\sum_{k=1}^{m-1}(\tau_k - b)_+$ is continuous in $b$, attaining a value $>1$ at $b=0$ but a value $<1$ at a sufficiently large $b$.
By examining (\ref{eq:cw-loss3}), we have
$L^{\text{cw3}} (m,\tau^\dprime) \le L^{\text{cw3}} (m,\tau)$ because $\tau^\dprime_k \le \tau_k$ for each $k\in[m-1]$.
Moreover,  $L^{\text{cw3}} (j,\tau^\dprime) \le L^{\text{cw3}} (j,\tau)$ for $j \in [m-1]$, by noting that
$(\tau^\dprime)^{(j)}_m %= 1- \tau^\dprime_j - \sum_{k\in[m-1],k\not=j} \tau^\dprime_{k+}
= 1-\sum_{k\in[m-1]} \tau^\dprime_k=0$,
$\tau^{(j)}_m %= 1- \tau_j - \sum_{k\in[m-1],k\not=j} \tau_{k+}
= 1-\sum_{k\in[m-1]} \tau_k <0$, and hence $(\tau^\dprime)^{(j)}_{m+} = \tau^{(j)}_{m+}=0$.

By combining the preceding two steps, the minimum in (\ref{eq:prf-cw-H}) is achieved at some $\tau \in \bbR_+^{m-1}$ with $\sum_{k=1}^{m-1} \tau_k \le 1$, that is, satisfying $\tilde\tau \in\Delta_m$.
\hfill$\blacksquare$\\[2mm]

\noindent\textbf{Proof of Proposition~\ref{pro:hinge-regret-bd}(i).}\;
Note that $H_{L^{\text{cw}}} (\eta) = H_{L^{\text{cw3}}} (\eta)$ by Proposition~\ref{pro:cw-conj3}. Then inequality (\ref{eq:hinge-regret-bd}) is equivalent to
\begin{align}
\frac{1}{m} R_{L^{\text{cw}}} (\eta,\tau^\dag) + \frac{m-1}{m} H_{L^{\text{cw}}} (\eta) \le R_{L^{\text{cw3}}} (\eta,\tau) ,  \label{eq:prf-hinge-bd}
\end{align}
where $\tau^\dag = (\tau_1,\ldots, \tau_{m-1}, 1-\sum_{k=1}^{m-1} \tau_{k+})^\T$. We distinguish three cases.

In the first case, suppose that  $\tau \in \bbR^{m-1}$ with one or more negative components.
We show that for any $\eta \in \Delta_m$,
\begin{align*}
R_{L^{\text{cw}}} (\eta, {\tau^\prime}^\dag) = R_{L^{\text{cw}}} (\eta, \tau^\dag), \quad R_{L^{\text{cw3}}} (\eta, \tau^\prime)  \le R_{L^{\text{cw3}}} (\eta,\tau) ,
\end{align*}
where $\tau^\prime$ is obtained from $\tau$ by setting all negative components of $\tau$ to 0.
The second inequality follows from (\ref{eq:prf-cw-conj3-ineq1}) directly.
To see the first equality, note that a maximum component among $\tau^\dag= (\tau_1, \ldots, \tau_{m-1}, \tau^\dag_m)^\T$ must be positive;
otherwise, $\tau_j \le 0$ for each $j\in[m-1]$ and hence $\tau^\dag_m =1$, a contradiction.
A maximum component among ${\tau^\prime}^\dag= (\tau_1^\prime, \ldots, \tau_{m-1}^\prime, {\tau^\prime}^\dag_m)^\T$ must also be positive.
But for $j\in[m-1]$, we have $\tau_j^\prime = \tau_j$ whenever $\tau_j^\prime$ or $\tau_j$ is positive.
Moreover, we have ${\tau^\prime}^\dag_m = \tau^\dag_m$, regardless of the signs of ${\tau^\prime}^\dag_m$ and $\tau^\dag_m$, because $\tau^\prime_{k+} = \tau_{k+}$ for each $k\in[m-1]$.
Therefore, $\argmax_{j\in [m]} (\tau^\prime)^\dag_j$ and  $\argmax_{j\in [m]} \tau^\dag_j$ can be set to be same, and the first equality above holds.

In the second case, suppose that $\tau\in \bbR_+^{m-1}$ (i.e., all components of $\tau$ are nonnegative) with $\sum_{k=1}^{m-1} \tau_k >1$.
We show that for any $\eta \in \Delta_m$,
\begin{align*}
R_{L^{\text{cw}}} (\eta, {\tau^\dprime}^\dag) = R_{L^{\text{cw}}} (\eta, \tau^\dag), \quad
R_{L^{\text{cw3}}} (\eta, \tau^\dprime)  \le R_{L^{\text{cw3}}} (\eta,\tau) ,
\end{align*}
where $\tau^\dprime = (\tau^\dprime_1, \ldots, \tau^\dprime_{m-1})^\T \in \bbR_+^{m-1}$ are defined as in Proof of Proposition~\ref{pro:cw-conj3}.
The second inequality follows from (\ref{eq:prf-cw-conj3-ineq2}) directly.
To see the first equality, note that  $\argmax_{j\in [m]} (\tau^\dprime)^\dag_j$ and  $\argmax_{j\in [m]} \tau^\dag_j$ must lie in the set $[m-1]$ because
$(\tau^\dprime)^\dag_m= 1-\sum_{k\in[m-1]} \tau^\dprime_k=0$ and
$\tau^\dag_m= 1-\sum_{k\in[m-1]} \tau_k <0$.
But the first $m-1$ components of $\tau^\dprime$, $(\tau^\dprime)^\dag_j = \tau^\dprime_j = (\tau_j-b)_+$ for $j\in[m-1]$, are ordered in the same way as
those of $\tau$.
Hence $\argmax_{j\in [m]} (\tau^\dprime)^\dag_j$ and  $\argmax_{j\in [m]} \tau^\dag_j$ can be set to be same, and the desired equality holds.

From the preceding discussion, it suffices to show (\ref{eq:prf-hinge-bd}) in the third case where $\tau \in \bbR_+^{m-1}$ with $\sum_{k=1}^{m-1} \tau^\dag_k \le 1$,
and hence $\tau^\dag =(\tau^\dag_1, \ldots, \tau^\dag_{m-1}, 1- \sum_{j \in [m-1]} \tau^\dag_j)^\T \in \Delta_m$.
Let $k = \argmin_{j\in [m]} \eta^\T C_j$ and $l = \argmax_{j \in [m]} \tau^\dag_j$.  Then $\tau^\dag_l \ge m^{-1}$ and
\begin{align*}
& R_{L^{\text{cw}}} (\eta,\tau^\dag) = \eta^\T C_l, \quad H_{L^{\text{cw}}} (\eta) = \eta^\T C_k.
\end{align*}
Moreover, direct calculation yields
\begin{align*}
& R_{L^{\text{cw3}}} (\eta,\tau) = \sum_{j=1}^m \eta^\T C_j \tau^\dag_j
\ge \tau^\dag_l \eta^\T C_l + (1-\tau^\dag_l) \eta^\T C_k .
\end{align*}
The right-hand side above is non-decreasing in $\tau^\dag_l$ because $\eta^\T C_l \ge \eta^\T C_k$,
and hence is no smaller than its value at $\tau^\dag_l=m^{-1}$, that is, the left-hand side of (\ref{eq:prf-hinge-bd}).
\hfill$\blacksquare$\\[2mm]

\noindent\textbf{Proof of equivalence between $L^{\text{LLW}}$ and $L^{\text{LLW2}}$.}\;
Suppose that $\tau_k = (1+\gamma_k)/m$ for $k\in [m-1]$. Then it is immediate $L^{\text{LLW2}} (m,\tau) = L^{\text{LLW}} (m,\gamma)/m$.
Moreover, because $0=\sum_{k=1}^m \gamma_k  = \gamma_m + \sum_{k=1}^{m-1} (m\tau_k-1)$, we have
\begin{align*}
1 + \gamma_m = m - m \sum_{k=1}^{m-1} \tau_k.
\end{align*}
Substituting this into the definition of $L^{\text{LLW}}$ and using $1+\gamma_k = m \tau_k$ for $k\in[m-1]$
yields $L^{\text{LLW2}} (j,\tau) = L^{\text{LLW}} (j,\gamma)/m$ for $j \in [m-1]$.
\hfill$\blacksquare$\\[2mm]

\noindent\textbf{Comparison between $L^{\text{zo3}}$ and $L^{\text{LLW2}}$.}\;
On one hand, the two losses $L^{\text{zo3}}$ and $L^{\text{LLW2}}$ share some similar properties.
It can be verified that, similarly to $L^{\text{zo3}}$,
$L^{\text{LLW2}}$ is a convex extension of
$L^{\text{zo2}}$ in (\ref{eq:zo-loss2}), considered a function of $j$ and $(\lambda_1,\ldots,\lambda_{m-1})^\T$ with $\lambda_m=1-\sum_{k=1}^m \lambda_k$.
Moreover, by Proposition~\ref{pro:cw-conj3} and \citet[Example~5]{duchi2018}, the losses $L^{\text{zo3}}$ and
$L^{\text{LLW2}}$ lead to the same generalized entropy $H^{\text{zo}}$.
Our result, Proposition \ref{pro:hinge-regret-bd}, also yields a classification regret bound for $L^{\text{zo3}}$, similar to that for $L^{\text{LLW2}}$ in \cite{duchi2018}.
On the other hand, there are interesting differences between $L^{\text{zo3}}$ and $L^{\text{LLW2}}$.
While $L^{\text{zo3}}(j,\tau)$ and $L^{\text{LLW2}}(j,\tau)$ are aligned with $L^{\text{zo2}} (j, \tilde \tau)$
for $\tilde\tau = (\tau_1,\ldots,\tau_{m-1}, 1-\sum_{k=1}^{m-1} \tau_k)^\T \in \Delta_m$,
the loss $L^{\text{zo3}}$ stays uniformly lower than $L^{\text{LLW2}}$,
\begin{align*}
0 \le L^{\text{zo3}}(j,\tau) \le L^{\text{LLW2}}(j,\tau), \quad j \in [m], \tau \in \bbR^{m-1},
\end{align*}
because $L^{\text{zo3}}(j,\tau)$ can be written as $\sum_{k\in[m-1],k\not=j} \tau_{k+} + ( 1-\tau_j - \sum_{k\in[m-1], k\not=j} \tau_{k+} )_+ $ for $j \in [m-1]$.
Hence the loss $L^{\text{zo3}}$ is a tighter convex extension than $L^{\text{LLW2}}$.
Another remarkable difference is that $L^{\text{zo3}}(j,\tau)$ appears to be geometrically simpler with fewer non-differentiable ridges than $L^{\text{LLW2}}(j,\tau)$ for $j \in [m-1]$.
See Figure~\ref{fig:LLW2} for an illustration in the three-class setting.
Further research is needed on whether the aforementioned differences can be translated into advantages in classification performance.
\hfill$\blacksquare$\\[2mm]

\noindent\textbf{Proof of Proposition~\ref{pro:cw-conj4}.}\;
We need to show that for $\eta\in\Delta_m$,
\begin{align}
H^{\text{zo}}(\eta) = \inf_{\tau \in \bbR^{m-1}} \left\{\sum_{j=1}^m \eta_j L^{\text{zo4}} (j,\tau) \right\}. \label{eq:prf-cw4-H}
\end{align}
Similarly as in the proof of Proposition~\ref{pro:cw-conj3}, it suffices to show that
\begin{itemize}
\item[(i)] $L^{\text{zo4}}$ is an extension of $L^{\text{zo2}}$ from $\Delta_m$ to $\bbR^{m-1}$, and

\item[(ii)] the minimum in (\ref{eq:prf-cw4-H}) is achieved at $\tau \in \bbR^{m-1}$ such that $\tilde\tau\in\Delta_m$, where $\tilde\tau= (\tau_1, \ldots, \tau_{m-1}, 1-\sum_{k=1}^{m-1} \tau_k)^\T$.
\end{itemize}
We use the following equivalent expressions for $S^{(j)}_\tau$:
\begin{align}
S_\tau^{(j)} =  \max \left\{0, \tilde \tau_j-1, \frac{-\tilde \tau_{j(m-1)} }{m-1}, \frac{-\tilde \tau_{j(m-1)} -\tilde \tau_{j(m-2)} }{m-2}, \label{eq:Stau-expression1}
\ldots, \frac{-\tilde \tau_{j(m-1)} - \cdots - \tilde\tau_{j(2)} }{2} \right\} .
\end{align}
and, if $m_\tau^{(j)} \ge 1$,
\begin{align}
S_\tau^{(j)}  = \max \left\{\tilde \tau_j-1, \max_{ m_\tau^{(j)} \le l \le m-2 } \frac{ %-\sum_{k\in[m]: k\not=j, \tilde\tau_k \le 0}  \tilde \tau_k
- \sum_{k=1}^{l} \tilde \tau_{j(m-k)} } { m-l } \right\}, \label{eq:Stau-expression2}
\end{align}
where $m_\tau^{(j)} =\# \{k\in[m]: k\not=j, \tilde\tau_k \le 0\}$.
The first expression is immediate because $\sum_{k=1}^m \tilde\tau_k =1$. The second expression follows
because $\{\tilde\tau_{j(m-k)}: 1\le k \le m_\tau^{(j)}\}$, the smallest $m_\tau^{(j)}$ components among $\tilde\tau$ excluding $\tilde\tau_j$,
are $\{\tilde \tau_k: \tilde \tau_k\le 0, k\not=j, k \in [m]\}$, and  $- \sum_{k=1}^{l} \tilde \tau_{j(m-k)}  /(m-l)$ is nonnegative and nondecreasing in $1 \le l\le m_\tau^{(j)}  $.

Result (i) can be directly verified. For any $\tau\in\bbR^{m-1}$ such that $\tilde\tau \in \Delta_m$,
we have $S_\tau^{(j)} =0$ for $j\in [m]$ by examining the expression (\ref{eq:Stau-expression1}), and
hence $L^{\text{zo4}} (j, \tau) = L^{\text{zo2}}(j, \tilde\tau) = 1 - \tilde\tau_j$ for $j \in [m]$
by the definitions (\ref{eq:zo-loss4}) and (\ref{eq:zo-loss2}).

For result (ii), we show that for any $\tau\in\bbR^{m-1}$ with one or more negative components in $\tilde\tau$, there exists $\tau^\prime=(\tau^\prime_1,\ldots,\tau^\prime_{m-1})^\T\in \bbR^{m-1}$
such that $\tilde\tau^\prime=(\tau^\prime_1,\ldots,\tau^\prime_{m-1}, 1-\sum_{k=1}^{m-1} \tau^\prime_k)^\T \in \Delta_m$ and for $j\in[m]$,
\begin{align}
L^{\text{zo4}} (j,\tau^\prime) \le  L^{\text{zo4}} (j,\tau) . \label{eq:prf-cw-conj4-ineq}
\end{align}
Then the minimum in (\ref{eq:prf-cw4-H}) is achieved at some $\tau\in\bbR^{m-1}$ with $\tilde\tau \in\Delta_m$.

First, let $\tau^\dprime = (\tilde\tau^\dprime_1, \ldots, \tilde\tau^\dprime_{m-1})^\T$ and $\tilde\tau^\dprime = (\tilde\tau^\dprime_1,\ldots,\tilde\tau^\dprime_m)^\T$ with
\begin{align*}
\tilde \tau^\dprime_j = \left\{ \begin{array}{cl}
\tilde \tau_j -b , & \mbox{if } \tilde\tau_j \ge  0, \\
\tilde \tau_j + \frac{m_\tau^+}{m_\tau^-} b, & \mbox{if } \tilde\tau_j < 0,
\end{array} \right.
\end{align*}
for $j \in [m]$, where $m_\tau^- = \# \{k\in[m]: \tilde\tau_k < 0\} \ge 1$, $m_\tau^+ = \# \{k\in[m]: \tilde\tau_k \ge 0\} = m - m_\tau^-$, and $b >0$ is determined such that
$\max\{ \tilde\tau_k + (m_\tau^+/m_\tau^-)b : k\in[m], \tilde\tau_k <0\}$ equals
$\min(0, \min\{ \tilde\tau_k -b : k\in[m], \tilde\tau_k  \ge 0\})$. Then the following properties hold:
\begin{itemize}
\item[(a)] $\sum_{k=1}^m \tilde\tau^\dprime_k = \sum_{k=1}^m \tilde \tau_k =1$.
\item[(b)] The ordering among components of $\tilde\tau^\dprime$ remains the same as that among $\tilde\tau$.
\item[(c)] If $\tilde\tau_k \le 0$ then $\tilde\tau_k^\dprime \le 0$ for $k\in[m]$.
\end{itemize}
It can be shown that $L^{\text{zo4}} (j,\tau^\dprime) \le L^{\text{zo4}} (j,\tau)$ for $j \in[m]$, depending on the sign of $\tilde\tau_j$.
\begin{itemize}
\item Suppose $\tilde\tau_j \ge 0$. Then $m_\tau^- \le m_{\tau}^{(j)}$ by definition, and $ m_{\tau}^{(j)} \le m_{\tau^\dprime}^{(j)} $
by property (c). For $m_\tau^- \le l \le m-2$, by property (b),
\begin{align*}
& 1 - \tilde\tau^\dprime_j - \frac{\sum_{k=1}^{l} \tilde \tau^\dprime_{j(m-k)} } { m-l } -
\left\{ 1 - \tilde\tau_j - \frac{\sum_{k=1}^{l} \tilde \tau_{j(m-k)} } { m-l }\right\} \nonumber \\
& = -(\tilde\tau^\dprime_j - \tilde\tau_j) -
\frac{\sum_{k=1}^{m_\tau^-} (\tilde \tau^\dprime_{j(m-k)}-\tilde \tau_{j(m-k)}) + \sum_{k=1+m_\tau^-}^l (\tilde \tau^\dprime_{j(m-k)}-\tilde \tau_{j(m-k)})   } { m-l } \\
& = b - \frac{m_\tau^+ b - (l - m_\tau^-) b  } { m-l } = 0.
\end{align*}
By combining the preceding properties with (\ref{eq:Stau-expression2}),
\begin{align}
& L^{\text{zo4}} (j,\tau)  = \max \left\{0, 1-\tilde\tau_j+ \max_{ m_\tau^{(j)} \le l \le m-2 } \frac{
- \sum_{k=1}^{l} \tilde \tau_{j(m-k)} } { m-l } \right\}, \label{eq:new-expression1} \\
& L^{\text{zo4}} (j,\tau^\dprime)  = \max \left\{0, 1-\tilde\tau_j^\dprime + \max_{ m_{\tau^\dprime}^{(j)} \le l \le m-2 } \frac{
- \sum_{k=1}^{l} \tilde \tau^\dprime_{j(m-k)} } { m-l } \right\}. \label{eq:new-expression2}
\end{align}
we see that $L^{\text{zo4}} (j,\tau^\dprime) = L^{\text{zo4}} (j,\tau)$.

\item Suppose $\tilde\tau_j <0$ and $m_\tau^- \ge 2$.
Then $1\le m_\tau^- -1 \le m_{\tau}^{(j)}$ by definition, and $ m_{\tau}^{(j)} \le m_{\tau^\dprime}^{(j)} $ by property (c).
For $m_\tau^- -1 \le l \le m-2$, by property (b),
\begin{align*}
& 1 - \tilde\tau^\dprime_j - \frac{\sum_{k=1}^{l} \tilde \tau^\dprime_{j(m-k)} } { m-l } -
\left\{ 1 - \tilde\tau_j - \frac{\sum_{k=1}^{l} \tilde \tau_{j(m-k)} } { m-l }\right\} \nonumber \\
& = -(\tilde\tau^\dprime_j - \tilde\tau_j) -
\frac{\sum_{k=1}^{m_\tau^- -1} (\tilde \tau^\dprime_{j(m-k)}-\tilde \tau_{j(m-k)}) + \sum_{k=m_\tau^-}^l (\tilde \tau^\dprime_{j(m-k)}-\tilde \tau_{j(m-k)})   } { m-l } \\
& = -\frac{m_\tau^+}{m_\tau^-} b - \frac{(m_\tau^+ b -\frac{m_\tau^+}{m_\tau^-}b) - (l - m_\tau^- +1) b  } { m-l } =
-\left(\frac{m_\tau^+}{m_\tau^-} +1 \right)b + \frac{ (\frac{m_\tau^+}{m_\tau^-} +1 )b}{m-l} < 0 .
\end{align*}
Hence $L^{\text{zo4}} (j,\tau^\dprime) \le L^{\text{zo4}} (j,\tau)$ by the expressions (\ref{eq:new-expression1})--(\ref{eq:new-expression2}).

\item Suppose $\tilde\tau_j <0$ and $m_\tau^- =1$. Then $\tilde\tau_k \ge 0$ for $k \in [m]$ and $k \not= j$, and hence $L^{\text{zo4}} (j,\tau) = 1 -\tilde\tau_j$
by (\ref{eq:Stau-expression1}). Moreover, $\tilde\tau^\dprime_j \le 0$ and $\tilde\tau^\dprime_k \ge -b$ for $k\in [m]$ and $k \not =j$, and hence
by (\ref{eq:Stau-expression1}) applied to $\tau^\dprime$,
\begin{align*}
& L^{\text{zo4}} (j,\tau^\dprime) \le 1 - \tilde\tau^\dprime_j + (m-2) b/2 \\
& = 1-\tilde\tau_j - (m-1)b + (m-2)b/2 \le L^{\text{zo4}} (j,\tau).
\end{align*}
\end{itemize}

If $\tilde\tau^\dprime$ has no negative components, then $\tilde\tau^\dprime \in \Delta_m$ and (\ref{eq:prf-cw-conj4-ineq}) holds with $\tau^\prime = \tau^\dprime$.
Otherwise, the preceding mapping from $\tilde\tau$ to $\tilde\tau^\dprime$, denoted as $\mathcal F(\cdot)$, can be iteratively applied.
Let $\tilde\tau^{(0)} = \tilde\tau$ and for $i=1,2,\ldots$, if $\tilde\tau^{(i-1)}$ has one or more negative components, then
let $\tilde\tau^{(i)} = \mathcal F( \tilde \tau^{(i-1)})$. It suffices to show that this process necessarily terminates after finite steps.
The final iteration $\tilde\tau^{(i)}$ has no negative components and hence $\tilde\tau^{(i)} \in \Delta_m$.
The first $m-1$ components of $\tilde\tau^{(i)}$ can be taken as the desired $\tau^\prime$ in (\ref{eq:prf-cw-conj4-ineq}).

Denote the set of $m_\tau^-$ negative components of $\tilde\tau$ (or equivalently the $m_\tau^-$ smallest components of $\tilde\tau$)
as $0 > \tilde\tau_{j1} \ge \ldots \ge \tilde\tau_{j,m_\tau^-}$.
By property (b), $\tilde\tau^{(i)}_{j1} \ge \ldots \ge \tilde\tau^{(i)}_{j,m_\tau^-}$ remain
the smallest $m_\tau^-$ components of $ \tilde\tau^{(i)}$ for each $i \ge 1$.
It suffices to show that $ \tilde\tau^{(i)}_{j1}$ becomes 0 for a certain finite $i \ge 1$. Then
the number of negative components of $ \tilde\tau^{(i)}$ decreases to $m_\tau^- -1$ or smaller.
Applying this argument repeatedly shows that $\tilde\tau^{(i)}_{j,m_\tau^-}$ necessarily becomes 0
(or equivalently the number of negative components of $ \tilde\tau^{(i)}$ decreases to 0)
for a certain finite $i$, hence proving the finite-termination of the iterations.

Return to the mapping from $\tilde\tau$ to $\tilde\tau^\dprime = \tilde\tau^{(1)}$.
By the choice of $b$, $\tilde\tau^\dprime_{j1}$ either equals 0 or $\tilde\tau^\dprime_k $ for some $k\in[m]$ such that $\tilde\tau_k \ge 0$ but  $\tilde\tau^\dprime_k < 0$.
In the latter case, the number of negative components of $ \tilde\tau^\dprime$ increases to at least $m_\tau^- +1$.
Applying this argument repeatedly shows that $\tilde\tau^{(i)}_{j1}$ necessarily equals 0 for some $i \le m-m_\tau^-$.
Otherwise, the number of negative components of $\tilde\tau^{(i)}_{j1}$ would be $m$, which contradicts the fact that
all the components of $\tilde\tau^{(i)}$ sum up to 1, by property (a).
\hfill$\blacksquare$\\[2mm]

\noindent\textbf{Proof of Proposition~\ref{pro:hinge-regret-bd}(ii).}\;
Because $L^{\text{zo4}}$ induces the same generalized entropy as the zero-one loss by Proposition~\ref{pro:cw-conj4},
the result can be obtained from Proposition~\ref{pro:general-hinge-bd}. Alternatively, the following gives a direct proof, building on
the proof of Proposition~\ref{pro:cw-conj4}.

The main steps of the proof are similar as in the proof of Proposition~\ref{pro:hinge-regret-bd}(i).
First, note that $H_{L^{\text{zo}}} (\eta) = H_{L^{\text{zo4}}} (\eta)$ by Proposition~\ref{pro:cw-conj4}. Then inequality (\ref{eq:hinge-regret-bd2}) is equivalent to
\begin{align}
\frac{1}{m} R_{L^{\text{zo}}} (\eta,\tilde \tau) + \frac{m-1}{m} H_{L^{\text{zo}}} (\eta) \le R_{L^{\text{zo4}}} (\eta,\tau) . \label{eq:prf-hinge-bd2}
\end{align}

Second, for any $\tau\in\bbR^{m-1}$ with one or more negative components in $\tilde\tau$, there exists $\tau^\prime=(\tau^\prime_1,\ldots,\tau^\prime_{m-1})^\T\in \bbR^{m-1}$
such that $\tilde\tau^\prime=(\tau^\prime_1,\ldots,\tau^\prime_{m-1}, 1-\sum_{k=1}^{m-1} \tau^\prime_k)^\T \in \Delta_m$ and
for any $\eta \in \Delta_m$,
\begin{align*}
R_{L^{\text{zo}}} (\eta, \tilde\tau^\prime) = R_{L^{\text{zo}}} (\eta, \tilde\tau), \quad
R_{L^{\text{zo4}}} (\eta, \tau^\prime)  \le R_{L^{\text{zo4}}} (\eta,\tau) .
\end{align*}
The second equality follows from (\ref{eq:prf-cw-conj4-ineq}) directly. Moreover, in the proof of (\ref{eq:prf-cw-conj4-ineq}), $\tilde\tau^\prime$ is obtained
from $\tilde\tau$ by iteratively applying the mapping $\mathcal F(\cdot)$ from $\tilde\tau$ to $\tilde\tau^\dprime$. By property (b),
the ordering among components of $\tilde\tau$ is preserved (although not strictly preserved) under the mapping.
Hence $\argmax_{j\in[m]} \tilde\tau^\dprime_j$ and, through iterations, $\argmax_{j\in[m]} \tilde\tau^\prime_j$ can all be set to be same as $\argmax_{j\in[m]} \tilde\tau_j$.
The first equality holds.

Finally, it suffices to show (\ref{eq:prf-hinge-bd2}) for $\tau \in \bbR^{m-1}$ with $\tilde\tau\in\Delta_m$.
Let $k = \argmax_{j\in [m]} \eta_j$ and $l = \argmax_{j \in [m]} \tilde\tau_j$.  Then $\tilde\tau_l \ge m^{-1}$. Direct calculation yields
\begin{align*}
& R_{L^{\text{zo4}}} (\eta,\tau) = \sum_{j\in[m]} \eta_j (1-\tau_j) = 1- \sum_{j\in[m]} \eta_j  \tilde\tau_j ,
\end{align*}
and
\begin{align*}
& \frac{1}{m} R_{L^{\text{zo}}} (\eta,\tilde\tau) + \frac{m-1}{m} H_{L^{\text{zo}}} (\eta) \\
& =\frac{1}{m} (1-\eta_l ) + \frac{m-1}{m} (1-\eta_k)
= 1- \left( \frac{1}{m} \eta_l + \frac{m-1}{m} \eta_k \right).
\end{align*}
Inequality (\ref{eq:prf-hinge-bd2}) can be obtained by comparing the above two expressions:
$\sum_{j\in[m]} \eta_j  \tilde\tau_j$ is upper-bounded by $\eta_l \tilde\tau_l + \eta_k (1-\tilde\tau_l)$, which is nonincreasing in $\tilde\tau_l$
with $\eta_l \le \eta_k$, and hence is no greater than its value at $\tilde\tau_l=m^{-1}$.
\hfill$\blacksquare$\\[2mm]

\noindent\textbf{Comparison between $L^{\text{zo4}}$ and $L^{\text{DKR2}}$.}\;
The comparison is similar to that between $L^{\text{zo3}} $ and $L^{\text{LLW2}}$.
On one hand, it can be verified that $L^{\text{DKR2}}$ is a convex extension of $L^{\text{zo2}}$ similarly to $L^{\text{zo4}}$, and
by Proposition~\ref{pro:cw-conj4} and \citet[Example~3]{duchi2018}, both $L^{\text{zo3}}$ and $L^{\text{DKR2}}$ lead to the same
generalized entropy $H^{\text{zo}}$.
Our result, Proposition~\ref{pro:hinge-regret-bd}, also gives a classification regret bound for $L^{\text{zo4}}$, similar to that for $L^{\text{DKR2}}$ in \cite{duchi2018}.
On the other hand, there are interesting differences between $L^{\text{zo4}}$ and $L^{\text{DKR2}}$.
While $L^{\text{zo4}} (j, \tau)$ and $L^{\text{DKR2}}(j, \tau)$ coincide with  $L^{\text{zo2}} (j, \tilde\tau)$
provided $\tilde \tau \in \Delta_m$, the loss $L^{\text{zo4}}$ gives a tighter convex extension than $L^{\text{DKR2}}$:
\begin{align*}
0 \le L^{\text{zo4}} (j, \tau) \le L^{\text{DKR2}}(j, \tau), \quad j \in [m], \tau \in \bbR^{m-1},
\end{align*}
because $S^{(j)}_\tau\le S_{\tilde\tau}$ for $j\in[m]$, with $S^{(j)}_\tau$ being the maximum of $m$ numbers which are respectively no greater than those in the definition of $S_{\tilde\tau}$.
Moreover, $L^{\text{zo4}}(j,\tau)$ appears to be geometrically simpler with fewer non-differentiable ridges than $L^{\text{DKR2}}(j, \tau)$ for $j \in [m]$.
See Figure~\ref{fig:DKR2} for an illustration in the three-class setting.
\hfill$\blacksquare$\\[2mm]

\subsection{Proofs of results in Section~\ref{sec:general-hinge-bd}}

\noindent\textbf{Proof of Proposition~\ref{pro:value-manifold}.} \;
Denote by $v_1, \ldots, v_m$ the vertices of $\mathcal S^{\text{zo}}$, where $v_j \in \bbR^m$ has $j$th component 0 and the remaining components 1.
For two vectors $x,y \in \bbR^m$, write $x \preceq z$ if $x_j \le y_j$ for $j \in [m]$.

First, suppose that the inclusion property (\ref{eq:val-manifold-inclusion}) holds. Then by (\ref{eq:entropy-manifold}), we have for $\eta\in \Delta_m$,
\begin{align*}
& H_L(\eta)  = \inf_{z \in \mathcal S_L} \eta^\T z = \inf_{z \in \mathcal S^{\text{zo}}} \eta^\T z \\
& = \inf_{\lambda\in \Delta_m} \sum_{j\in[m]} \lambda_j \eta^\T v_j = \inf_{\lambda\in \Delta_m} \sum_{j\in[m]} \lambda_j (1-\eta_j) = 1- \max_{j \in [m]} \eta_j .
\end{align*}
The equality $\inf_{z \in \mathcal S_L} \eta^\T z = \inf_{z \in \mathcal S^{\text{zo}}} \eta^\T z $ on the first line can be shown as follows.
On one hand,
$\inf_{z \in \mathcal S_L} \eta^\T z  \le  \inf_{z \in \mathcal S^{\text{zo}}} \eta^\T z $ because $\mathcal S_L \supset \mathcal S^{\text{zo}} $.
On the other hand,
$\inf_{z \in \mathcal S_L} \eta^\T z  \ge  \inf_{z \in \mathcal S^{\text{zo}*}} \eta^\T z  \ge \inf_{z \in \mathcal S^{\text{zo}}} \eta^\T z $ because
$\mathcal S_L \subset \mathcal S^{\text{zo}*}$ and  for any $z \in \mathcal S^{\text{zo}*}$ there exists $\tilde z \in \mathcal S^{\text{zo}}$ such that
$\tilde z \preceq z$ by the definition of $\mathcal S^{\text{zo}*}$.

Next, suppose that $H_L(\eta)  =   1- \max_{j \in [m]} \eta_j$ for $\eta \in \Delta_m$. We show that the inclusion property (\ref{eq:val-manifold-inclusion}) holds.
To argue by contradiction, there are two possible alternative cases if (\ref{eq:val-manifold-inclusion}) does not hold.

In the first case, $\mathcal S^{\text{zo}} \subset \mathcal S_L$ but $\mathcal S_L \not\subset \mathcal S^{\text{zo}*}$.
Then there exists a point $x \in \mathcal S_L$ but $x \not\in \mathcal S^{\text{zo}*}$.
The set $\mathcal S^{\text{zo}*}$ is easily seen to be closed and convex.
By the support hyperplane theorem, there exists a hyperplane which strictly separates $x$ and $\mathcal S^{\text{zo}*}$, that is, there exists some $\eta \in \bbR^m$ and $b \in \bbR$
such that $\eta^\T x < b$, but $\eta^\T z > b$ for all $z \in \mathcal S^{\text{zo}*}$.
The coefficient vector $\eta$ must be nonzero, $\eta\not=0$, and have all nonnegative components, $\eta \in \bbR^m_+$. Otherwise, suppose that, for example, $\eta_1 <0$ and
fix some point $\tilde z\in \mathcal S^{\text{zo}}$.
Define $\tilde z_k = \tilde z + k e_1$, where $e_1 = (1,0,\ldots,0)^\T$. Then
$\tilde z_k \in \mathcal S^{\text{zo}*}$ for all $k$, but $\eta^\T \tilde z_k = \eta^\T \tilde z + k \eta_1 \to -\infty$ as $k\to \infty$,
which contradicts the fact that $\eta^\T z > b$ for all $z \in \mathcal S^{\text{zo}*}$.
Hence $\eta$ can be normalized such that $\eta \in \Delta_m$. But then
\begin{align*}
\inf_{z \in \mathcal S_L} \eta^\T z \le \eta^\T x < \inf_{z \in \mathcal S^{\text{zo}*}} \eta^\T z = \inf_{z \in \mathcal S^{\text{zo}}} \eta^\T z = 1 - \max_{k\in [m]} \eta_k,
\end{align*}
a contradiction to the assumption that $H_L(\eta)  =   1- \max_{j \in [m]} \eta_j$.

In the second case, $\mathcal S^{\text{zo}} \not\subset \mathcal S_L$.
Then there exists a vertex of $\mathcal S^{\text{zo}}$ which is not contained in $\mathcal S_L$; otherwise $\mathcal S^{\text{zo}} \subset S_L$ by the convexity of $\mathcal S_L$.
Without loss of generality, assume that $v_1 \not\in \mathcal S_L$. Then
$v_1 \not\in \mathcal S_L + \bbR^m_+$. Otherwise, there exist some $x \in \mathcal S_L$ and $y \,(\not= 0) \in \bbR^m_+$ such that $v_1 = x + y$.
Then $\sum_{j\in [m]} (x_j + y_j) = m-1$, which contradicts the fact that $\sum_{j\in [m]} (x_j + y_j) > \sum_{j \in [m]} x_j \ge  m-1$.
The second equality holds because $\inf_{z \in \mathcal S_L} \sum_{j\in[m]} z_j = m-1$
by the assumption that $H_L(1_m/m) = \inf_{z \in \mathcal S_L} (1_m^\T/m) z = 1- 1/m$ for $1_m /m= (1/m,\ldots,1/m)^\T\in \Delta_m$.
The set $\mathcal S_L + \bbR^m_+$ is closed and convex. By the support hyperplane theorem,
there exists a hyperplane which strictly separates $v_1$ and $\mathcal S_L + \bbR^m_+$, that is, there exists some $\eta \in \bbR^m$ and $b \in \bbR$
such that $\eta^\T v_1 < b$, but $\eta^\T z > b$ for all $z \in \mathcal S_L + \bbR^m_+$.
Similarly as in the first case, $\eta$ must be nonzero, $\eta\not=0$, and have all nonnegative components, $\eta \in \bbR^m_+$.
Hence $\eta$ can be normalized such that $\eta \in \Delta_m$. But then
\begin{align*}
\inf_{z \in \mathcal S_L} \eta^\T z \ge \inf_{z \in \mathcal S_L + \bbR^m_+} \eta^\T z > \eta^\T v_1 \ge \inf_{z \in \mathcal S^{\text{zo}}} \eta^\T z = 1 - \max_{k\in [m]} \eta_k,
\end{align*}
again a contradiction to the assumption that $H_L(\eta)  =   1- \max_{j \in [m]} \eta_j$.

Combining the preceding two cases shows that (\ref{eq:val-manifold-inclusion}) holds as desired.
\hfill$\blacksquare$\\[2mm]

\noindent\textbf{Proof of Proposition~\ref{pro:general-hinge-bd}.} \;
Note that $H_L  ( \eta) = H_{L^{\text{zo}}}(\eta) $ by assumption. Inequality (\ref{eq:general-hinge-bd}) reduces to
\begin{align*}
\frac{1}{m}   R_{L^{\text{zo}}} (\eta,\sigma_L(\gamma)) + \frac{m-1}{m} H_{L^{\text{zo}}} (\eta) \le R_L(\eta,\gamma).
\end{align*}
By definition, $\sigma_L(\gamma) =  (-L(1,\gamma), \ldots, -L(m,\gamma))^\T $.
The preceding inequality can be stated such that for $\eta\in\Delta_m$ and $z = (L(1,\gamma), \ldots, L(m,\gamma))^\T \in \mathcal R_L$,
\begin{align}
\frac{1}{m} (1- \eta_l) + \frac{m-1}{m} (1- \eta_k ) \le \eta^\T z , \label{eq:prf-general-hinge-bd}
\end{align}
where  $l = \argmin_{j\in [m]} z_j$, and $k = \argmax_{j\in[m]} \eta_j$.
In the following, we show that (\ref{eq:prf-general-hinge-bd}) holds for $\eta\in\Delta_m$ and $z \in \mathcal S_L$. The notation $\preceq$ is used as in the proof of Proposition~\ref{pro:value-manifold}.

First, we show that for any $ z\in\mathcal S_L$, there exists some $\tilde z \in \mathcal S^{\text{zo}}$ such that
\begin{align}
\tilde z \preceq z, \quad \argmin_{j\in [m]} \tilde z_j = \argmin_{j \in [m]} z_j, \label{eq:prf2-general-hinge-bd}
\end{align}
which means that $\argmin_{j\in [m]} \tilde z_j$ can be set to be same as $\argmin_{j \in [m]} z_j$. Because $\mathcal S_L \subset \mathcal S^{\text{zo}*}$ by Proposition~\ref{pro:value-manifold},
it suffices to show that for any $z \in \mathcal S^{\text{zo}*}$, there exists $\tilde z \in \mathcal S^{\text{zo}}$ such that (\ref{eq:prf2-general-hinge-bd}) holds.
Without loss of generality, assume that $z_1 \ge z_2 \ge \cdots \ge z_m$.
Let $b = \sup\,\{ b^\prime \ge 0: z- b^\prime \, e_m \in \mathcal S^{\text{zo}*}\}$, such that
$z^\prime = z- b \, e_m \in \partial \mathcal S^{\text{zo}*}$.
Then $z^\prime \preceq z$ and the $m$th component of $z$, $z^\prime_m$, remains a minimum component of $z^\prime$.
By the definition of $\mathcal S^{\text{zo}*}$, there exists some $\tilde z \in \mathcal S^{\text{zo}}$ satisfying $\tilde z \preceq z^\prime$.
For any such point $\tilde z$, we have
\begin{itemize}
\item[(i)] $\tilde z_m = z^\prime_m$ and
\item[(ii)] $\tilde z_j \ge z^\prime_m$ for $j \in [m-1]$,
\end{itemize}
which then imply that (\ref{eq:prf2-general-hinge-bd}) is satisfied.
Property (i) follows because if $\tilde z_m < z^\prime_m$, then by the definition of $\mathcal S^{\text{zo}*}$,
$z^\prime - (z^\prime_m - \tilde z_m) e_m = \tilde z + (z^\prime-\tilde z -(z^\prime_m - \tilde z_m) e_m) \in \mathcal S^{\text{zo}*}$,
but this contradicts the definition of $b$.
To show property (ii), suppose that there exists $\tilde z \in \mathcal S^{\text{zo}}$ such that
$\tilde z \preceq z$ and $\tilde z_j < z^\prime_m$ for some $j\in[m-1]$.
Let $\tilde{\tilde z} = (\tilde z_1, \ldots, \tilde z_{j-1}, \tilde z_m, \tilde z_{j+1}, \ldots, \tilde z_j)^\T$, by exchanging the $j$th and $m$th components of $\tilde z$.
Then $\tilde{\tilde z} \in \mathcal S^{\text{zo}}$ by symmetry of $\mathcal S^{\text{zo}}$. Moreover,
$\tilde{\tilde z} \preceq z^\prime$, because $\tilde{\tilde z}_j = \tilde z_m = z^\prime_m \le z^\prime_j$ by property (i) and $\tilde{\tilde z}_m = \tilde z_j < z^\prime_m$ for $j \in [m-1]$.
Then $\tilde{\tilde z}$ must also satisfy property (i), i.e., $\tilde{\tilde z}_m = z^\prime_m$, a contradiction.

By the preceding result, it suffices to show that (\ref{eq:prf-general-hinge-bd}) holds for any $\eta \in \Delta_m$ and $\tilde z \in \mathcal S^{\text{zo}}$.
This can be obtained as follows:
\begin{align*}
& \eta^\T \tilde z = 1- \sum_{j \in [m]} \eta_j (1-\tilde z_j) =  1- \eta_l (1-\tilde z_l)  - \sum_{j\not=l} \eta_j (1-\tilde z_j) \\
& \ge  1- \eta_l (1-\tilde z_l) - \sum_{j\not=l} \eta_k (1-\tilde z_j)  = 1- \eta_l + (\eta_l -\eta_k) \tilde z_l ,\\
& \ge 1- \eta_l + \frac{m-1}{m} ( \eta_l -\eta_k) .
\end{align*}
The second line above uses the fact that $\eta_k = \max_{j\in[m]} \eta_j$, $(0\le)\,\tilde z_j \le 1$ for $j\in[m]$, and $\sum_{j\in[m]} \tilde z_j= m-1$.
The last line holds because $1- \eta_l + (\eta_l - \eta_k) \tilde z_l$ is non-increasing in $\tilde z_l$ with $\eta_l \le \eta_k$, and hence is no smaller than its value at $\tilde z_l = \frac{m-1}{m}$,
where $\tilde z_l =\min_{j\in[m]} \tilde z_j \le  \frac{m-1}{m}$ with $\sum_{j\in[m]} \tilde z_j = m-1$.
\hfill$\blacksquare$\\[2mm]

\noindent\textbf{Simplification of prediction mapping $\sigma_L$.} \;
We show that for each of the four losses, $L^{\text{LLW2}}$, $L^{\text{DKR2}}$, $L^{\text{zo3}}$, and $L^{\text{zo4}}$,
the prediction mapping $\sigma_L$ in Proposition~\ref{pro:general-hinge-bd}
is monotonically related to that in the corresponding regret bound discussed in Section~\ref{sec:hinge-regret-bd}.

The loss $L^{\text{LLW2}}$ can be written as
\begin{align*}
 L^{\text{LLW2}} (j, \tau)=
\sum_{k\in[m], k\not=j} \tilde \tau_{k+} = -\tilde \tau_j + \sum_{k\in[m]} \tilde \tau_{k+}, \quad j \in [m],
\end{align*}
where $ \tilde\tau = (\tau_1, \ldots, \tau_{m-1}, 1- \sum_{k=1}^{m-1} \tau_k )^\T$.
Hence if $\tilde \tau_j \le \tilde \tau_k$, then $L^{\text{LLW2}} (j, \tau) \ge L^{\text{LLW2}} (k, \tau)$.
Similarly, it is easily seen that
if $\tilde \tau_j \le \tilde \tau_k$, then $L^{\text{DKR2}} (j, \tau) \ge L^{\text{DKR2}} (k, \tau)$.

The loss $L^{\text{zo3}}$ can be written as
\begin{align*}
 L^{\text{zo3}} (j,\tau)= \left\{
\begin{array}{cl}
 \max\left(1-\tau_j, \,1- \tau^\dag_m - \tau_{j+} \right), & \mbox{if } j \in [m-1],\\
 1- \tau^\dag_m , & \mbox{if } j=m,
\end{array} \right.
\end{align*}
where $\tau^\dag = (\tau_1, \ldots, \tau_{m-1}, 1- \sum_{k=1}^{m-1} \tau_{k+} )^\T$.
For $j,k\in[m-1]$, if $ \tau_j \le  \tau_k$, then $L^{\text{zo3}} (j, \tau) \ge L^{\text{zo3}} (k, \tau)$ trivially.
For $j\in [m-1]$, if $\tau_j \le \tau^\dag_m$, then $L^{\text{zo3}} (j, \tau) \ge 1-\tau_j \ge  L^{\text{zo3}} (m, \tau)$, and
if $\tau_j \ge \tau^\dag_m$, then $L^{\text{zo3}} (j, \tau) \le 1-\tau^\dag_m =  L^{\text{zo3}} (m, \tau)$.

The loss $L^{\text{zo4}}$ can be written as
\begin{align*}
  L^{\text{zo4}} (j,\tau) =\max &  \left\{0, \frac{1+ \tilde\tau_{j(1)} - \tilde\tau_j }{2},  \ldots,
 \frac{m-2 +(\tilde\tau_{j(1)} -\tilde\tau_j) + \cdots+(\tilde\tau_{j(m-2)} - \tilde\tau_j) }{m-1} ,  \right. \\
 & \qquad  \left. \frac{m-1 +(\tilde\tau_{j(1)} -\tilde\tau_j) + \cdots+(\tilde\tau_{j(m-1)} - \tilde\tau_j) }{m} \right\}, \quad j \in [m],
\end{align*}
where $\tilde \tau_{j(1)} \ge \cdots \ge \tilde\tau_{j(m-1)}$ are the sorted components of $\tilde\tau$ excluding $\tilde\tau_j$.
Without loss of generality, assume that $\tilde\tau_1 \ge \ldots \ge \tilde \tau_m$. Then
\begin{align*}
  L^{\text{zo4}} (j,\tau) =\max &  \left\{0, \frac{1+ \tilde\tau_{1} - \tilde\tau_j }{2},  \ldots,
 \frac{j-1 +(\tilde\tau_{1} -\tilde\tau_j) + \cdots+(\tilde\tau_{j-1} - \tilde\tau_j) }{j} ,  \right. \\
 & \qquad  \left. \ldots, \frac{j +(\tilde\tau_{1} -\tilde\tau_j) + \cdots+(\tilde\tau_{j-1} - \tilde\tau_j) +(\tilde\tau_{j+1} - \tilde\tau_j) }{j+1} , \right. \\
 & \qquad  \left. \ldots, \frac{m-1 +(\tilde\tau_{1} -\tilde\tau_j) + \cdots+(\tilde\tau_{j-1} -\tilde\tau_j) +(\tilde\tau_{j+1} - \tilde\tau_j) +(\tilde\tau_{m} - \tilde\tau_j) }{m} \right\} .
\end{align*}
Denote the $i$th term in the curly brackets above as $\ell_i(j,\tau)$ for $i=1,\ldots,m$, that is, $\ell_1(j,\tau)=0$ and for $i=2,\ldots,m$,
\begin{align*}
\ell_i(j,\tau) = \left\{ \begin{array}{cl}
\frac{i-1 + \sum_{h=1}^{i-1} (\tilde \tau_h - \tilde \tau_j) }{i}, & \mbox{if } i <j,\\
\frac{i-1 + \sum_{h=1}^{i} (\tilde \tau_h - \tilde \tau_j) }{i}, & \mbox{if } i \ge j .
\end{array} \right.
\end{align*}
For $j>k$ with $\tilde\tau_j \le \tilde \tau_k$, if $ i<k$ or $i \ge j$, then $\ell_i(j,\tau) \ge \ell_i(k,\tau)$, and
if $ k \le i <j$, then
\begin{align*}
& \ell_i(j,\tau) = \frac{i-1 + \sum_{h=1}^{i-1} (\tilde \tau_h - \tilde \tau_j) }{i} \\
& \ge \frac{i-1 + \sum_{h=1}^{i-1} (\tilde \tau_h - \tilde \tau_k) }{i}
\ge \frac{i-1 + \sum_{h=1}^{i} (\tilde \tau_h - \tilde \tau_k) }{i} = \ell_i(k,\tau),
\end{align*}
where the second inequality follows because $\tilde \tau_i \le \tilde \tau_k$.
In summary, if $\tilde\tau_j \le \tilde \tau_k$, then $\ell_i(j,\tau) \ge \ell_i(k,\tau)$ for $i=1,\ldots,m$, and hence $L^{\text{zo4}} (j, \tau) \ge L^{\text{zo4}} (k, \tau)$.
\hfill$\blacksquare$\\[2mm]

\subsection{Proofs of results in Section \ref{sec:zero-one}}

\noindent\textbf{Proof of Proposition~\ref{pro:bregman-bd}.}\;
We use the fact that for a twice differentiable generalized entropy $H_L$, the Bregman divergence can be written as a double integral of a quadratic form
based on the Hessian matrix of $\nabla^2 H_L$:
\begin{align}\label{eq:bd_expansion}
    B_L(\eta,q) = -\int_0^1\int_0^1(\eta-q)^\T \nabla^2 H_L(q+ts(\eta-q))(\eta-q) t \, \dif s \dif t, \quad \eta,q\in\Delta_m.
\end{align}
Identity (\ref{eq:bd_expansion}) follows from a second-order Taylor expansion with an integral remainder for
the univariate function $H_L(q + t(\eta-q))$ with $t \in [0,1]$.

(i) By definition (\ref{eq:LtoH}),
the generalized entropy corresponding to the pairwise (symmetrized) loss $L$ in (\ref{eq:pairwise-s}) is $H_L(q) = - \sum_{i=1}^m\sum_{j\neq i} q_i f_0(q_j/q_i)$. See also Supplement Table~S1.
The first-order and second-order derivatives of $H(q)$ are
\begin{align*}
    \frac{\partial H}{\partial q_i} &= -\sum_{j\neq i}\left\{f_0(\frac{q_j}{q_i})-\frac{q_j}{q_i}f_0'(\frac{q_j}{q_i}) + f_0'(\frac{q_i}{q_j})\right\},\\
    \frac{\partial^2 H}{\partial q_i^2} &= -\sum_{j\neq i}\left\{\frac{q_j^2}{q_i^3}f_0^\dprime(\frac{q_j}{q_i})+\frac{1}{q_j}f_0^\dprime(\frac{q_i}{q_j})\right\},\\
    \frac{\partial^2 H}{\partial q_j \partial q_i} &= \frac{q_j}{q_i^2}f_0^\dprime(\frac{q_j}{q_i})+\frac{q_i}{q_j^2}f_0^\dprime(\frac{q_i}{q_j}),\quad j\neq i.
\end{align*}
By the relationship $w(q_1) =f_0^\dprime(u^q) / q_2^3$, we obtain $f_0^\dprime(u^q) = 2^{2\nu} q_1^{\nu-1} q_2^{\nu+2}$ from $w(q_1) = 2^{2\nu} q_1^{\nu-1} q_2^{\nu-1}$.
Then the quadratic form $-x^\T\nabla^2 H_L(\tilde\eta)x$ with $x=\eta-q$ and $\tilde\eta = q+ts(\eta-q)\in \Delta_m$ can be written as
\begin{align*}
     -x^\T\nabla^2H_L(\tilde\eta)x  &= \sum_{i=1}^m \sum_{j=1}^m \left[ \left \{ \frac{\tilde\eta_j^2}{\tilde\eta_i^3} f_0^\dprime(\frac{\tilde\eta_j}{\tilde\eta_i}) + \frac{1}{\tilde\eta_j} f_0^\dprime(\frac{\tilde\eta_i}{\tilde\eta_j})  \right\} x_i^2 - \left\{ \frac{\tilde\eta_j}{\tilde\eta_i^2} f_0^\dprime(\frac{\tilde\eta_j}{\tilde\eta_i})  +\frac{\tilde\eta_i}{\tilde\eta_j^2} f_0^\dprime(\frac{\tilde\eta_i}{\tilde\eta_j}) \right\}x_ix_j\right]\\
    &= 2^{2\nu} \sum_{i=1}^m \sum_{j=1}^m \frac{(\tilde\eta_i\tilde\eta_j)^\nu}{(\tilde\eta_i+\tilde\eta_j)^{2\nu+1}}(\tilde\eta_j^{\frac{1}{2}}\tilde\eta_i^{-\frac{1}{2}}x_i-\tilde\eta_i^{\frac{1}{2}}\tilde\eta_j^{-\frac{1}{2}}x_j)^2.
\end{align*}
With $\nu\leq 0$, note that $(\tilde\eta_i \tilde\eta_j)^\nu/(\tilde\eta_i+\tilde\eta_j)^{2\nu+1}\geq 2^{-2\nu}$ because $\tilde\eta_i\tilde\eta_j\leq 2^{-2}(\tilde\eta_i+\tilde\eta_j)^2$
and $\eta_i+\eta_j\leq 1$ for each pair $(i,j)$. Hence we have
\begin{align*}
-x^\T\nabla^2H_L(\tilde\eta)x &\geq  \sum_{i=1}^m \sum_{j=1}^m(\tilde\eta_j^{\frac{1}{2}}\tilde\eta_i^{-\frac{1}{2}}x_i-\tilde\eta_i^{\frac{1}{2}}\tilde\eta_j^{-\frac{1}{2}}x_j)^2\\
    &=2 \left\{ (\sum_{i=1}^m \tilde\eta_i^{-1}x_i^2)(\sum_{j=1}^m\tilde\eta_j) -(\sum_{i=1}^m x_i)^2\right\}
    \geq 2 \|x\|_1^2,
\end{align*}
where the last inequality follows because $(\sum_{i=1}^m \tilde\eta_i^{-1}x_i^2)(\sum_{j=1}^m\tilde\eta_j) \ge (\sum_{i=1}^m |x_i|)^2$ by the Cauchy--Schwartz inequality
and $\sum_i x_i = \sum_i p_i -\sum_i \eta_i = 0$. Combining this lower bound with (\ref{eq:bd_expansion}) and integrating over $s$ and $t$ yield $\kappa_L = 2$.

(ii\,a) Suppose $\beta \in [1/2,1]$. The generalized entropy corresponding to the simultaneous loss $L$ in (\ref{eq:simul-loss}) is $H(q)=\|q\|_\beta$. The first-order and second-order derivatives are
\begin{align*}
    \frac{\partial H}{\partial q_i} &= q_i^{\beta-1}\|q\|_\beta^{1-\beta},\\
    \frac{\partial^2 H}{\partial q_i^2}& = -(1-\beta)(\sum_{j\neq i}q_j^{\beta}q_i^{\beta-2})\|q\|_\beta^{1-2\beta},\\
    \frac{\partial^2 H}{\partial q_i \partial q_j}& = (1-\beta)(q_jq_i)^{\beta-1}\|q\|_\beta^{1-2\beta}, \quad j\neq i.
\end{align*}
The quadratic form $-x^\T\nabla^2H_L(\tilde\eta)x$ with $x=\eta-q$ and $\tilde\eta = q+ts(\eta-q)\in \Delta_m$ can be written as
\begin{align*}
   -x^\T\nabla^2H_L(\tilde\eta)x  =  \frac{1-\beta}{2}\left\{\sum_{i=1}^m\sum_{j=1}^m \|\tilde\eta\|_\beta^{1-2\beta} (\tilde\eta_i\tilde\eta_j)^{\beta-1} (\tilde\eta_j^{\frac{1}{2}}\tilde\eta_i^{-\frac{1}{2}}x_i-\tilde\eta_i^{\frac{1}{2}}\tilde\eta_j^{-\frac{1}{2}}x_j)^2\right\}.
\end{align*}
With $\beta \in [1/2,1)$ and hence $\beta-1<0$, it holds that $(\tilde\eta_i\tilde\eta_j)^{\beta-1}\geq 2^{2-2\beta}$ by inverting the inequality
$\tilde\eta_i\tilde\eta_j \leq 2^{-2}(\tilde\eta_i+\tilde\eta_j)^2\leq 2^{-2}$. In addition, $\|\tilde\eta\|_\beta$ is concave and attains the maximum
$m^{1-1/\beta}$ over $\Delta_m$ when $\tilde\eta_i=1/m$ for $i\in[m]$. Because $1-2\beta \le 0$, it follows that the minimum of $\|\tilde\eta\|_\beta^{1-2\beta}$ over $\Delta_m$
is $m^{(1-1/\beta)(2\beta-1)}$. Then the quadratic form is lower bounded by
\begin{align*}
   -x^\T\nabla^2H_L(\tilde\eta)x  &\geq \frac{1-\beta}{2} m^{\frac{(\beta-1)(2\beta-1)}{\beta}} 2^{2-2\beta}
   \left\{\sum_{i=1}^m\sum_{j=1}^m (\tilde\eta_j^{\frac{1}{2}}\tilde\eta_i^{-\frac{1}{2}}x_i-\tilde\eta_i^{\frac{1}{2}}\tilde\eta_j^{-\frac{1}{2}}x_j)^2\right\}\\
   & =  (1-\beta)m^{\frac{(\beta-1)(2\beta-1)}{\beta}} 2^{2-2\beta}
   \left\{(\sum_{i=1}^m \tilde\eta_i^{-1}x_i^2)(\sum_{j=1}^m\tilde\eta_j)-(\sum_{i=1}^m x_i)^2\right\}\\
   &\geq(1-\beta)m^{\frac{(\beta-1)(2\beta-1)}{\beta}} 2^{2-2\beta}\| x \|_1^2,
\end{align*}
where the last inequality follows similarly as in the proof of (i), by the Cauchy–-Schwartz inequality and $\sum_i x_i =0$.
Integration of (\ref{eq:bd_expansion}) with the preceding lower bound yields $\kappa_L = (1-\beta)m^{(1-1/\beta)(2\beta-1)} 2^{2-2\beta}$.

(ii\,b) Suppose $\beta \in (0,1/2]$. The generalized entropy, derivatives and quadratic form remain the same as in (iia). With $\beta \in (0,1/2]$ and hence $1-2\beta>0$, we have
\begin{align*}
    \|\tilde\eta\|_\beta^{1-2\beta}(\tilde\eta_i\tilde\eta_j)^{\beta-1}\geq (\tilde\eta_i^\beta+\tilde\eta_j^\beta)^{\frac{1-2\beta}{\beta}}(\tilde\eta_i\tilde\eta_j)^{\beta-1}
    = \left\{ \frac{(\tilde\eta_i\tilde\eta_j)^\beta}{(\tilde\eta_i^\beta + \tilde\eta_j^\beta)^{2}}\right\}^{1-\frac{1}{\beta}}
     (\tilde\eta_i^\beta + \tilde\eta_j^\beta)^{-\frac{1}{\beta}} \geq 2^{\frac{1}{\beta}-1}.
\end{align*}
The first inequality holds trivially. The second inequality holds because $1-1/\beta<0$, $(\eta_i\eta_j)^\beta \leq 2^{-2}(\eta_i^\beta+\eta_j^\beta)^2$,
and $(\tilde\eta_i^\beta+\tilde\eta_j^\beta)^{-1/\beta}$ is lower bounded by $2^{1-1/\beta}$. Similarly as in (ii\,a),  the quadratic form is lower bounded by
\begin{align*}
     -x^\T\nabla^2H_L(\tilde\eta)x  &\geq (1-\beta)  2^{\frac{1}{\beta}-1}
    \left\{(\sum_{i=1}^m \tilde\eta_i^{-1}x_i^2)(\sum_{j=1}^m\tilde\eta_j)-(\sum_{i=1}^m x_i)^2\right\} \geq (1-\beta) 2^{\frac{1}{\beta}-1}\|x\|_1^2,
\end{align*}
where the last inequality follows from the Cauchy–-Schwartz inequality and $\sum_i x_i =0$.
Integration of (\ref{eq:bd_expansion}) with the preceding lower bound yields $\kappa_L = (1-\beta) 2^{1/\beta-1}$.
\hfill$\blacksquare$\\[2mm]

\noindent\textbf{Discussion on the multinomial likelihood loss.}\;
The standard likelihood loss $L(j,q) = -\log q_j$ is equivalent to (\ref{eq:simul-loss}) in the limit of $\beta \to 1$ after properly rescaled.
By Proposition~\ref{pro:bregman-bd}(ii), inequality (\ref{eq:bregman-bd}) can be shown to hold for the rescaled entropy $H^{\text{r}}_\beta$ in (\ref{eq:H-beta2}) with
$\kappa_L = \kappa_\beta / ( m^{1/\beta-1} -1) $,
where $\kappa_\beta = (1-\beta) m^{(1-1/\beta)(2\beta-1)}2^{2-2\beta}$ if $\beta \in [1/2,1)$.
Then (\ref{eq:bregman-bd-KL}) can be recovered from (\ref{eq:bregman-bd}) as $\beta \to 1$, because
\begin{align*}
\lim_{\beta \to 1} \frac{\kappa_\beta } { m^{1/\beta-1} -1} = \lim_{\beta \to 1} \frac{1-\beta} { m^{1/\beta-1} -1} = (\log m)^{-1} .
\end{align*}
\hfill$\blacksquare$\\[2mm]

\noindent\textbf{Discussion on simultaneous exponential loss.}\;
The simultaneous exponential loss $L^{\text{r}}_0$ as used in \cite{zou2008}
can be obtained from (\ref{eq:simul-loss}) in the limit of $\beta\to 0+$ after properly rescaled, by Proposition~\ref{pro:limit}(i).
However, for $m\ge 3$, the corresponding modulus $\kappa_L$ from Proposition~\ref{pro:bregman-bd}(ii) for $L^{\text{r}}_\beta$ as $\beta \to 0+$ gives 0:
\begin{align*}
\lim_{\beta \to 0+} \frac{ (1-\beta)2^{1/\beta-1} } {m^{1/\beta-1}-1}  = 0 .
\end{align*}
The limit above gives 1 for  $m=2$, in agreement with the relationship $L^{\text{r}}_0 = L_{1/2}-1$ with $m=2$.
Our further calculation (not shown) suggests that a uniform bound in the form of (\ref{eq:bregman-bd}) might not be feasible on the associated Bregman divergence.
Hence an alternative approach would be needed to analyze $\underline\psi$ and deduce a concrete meaningful implication from regret bound (\ref{eq:zo-regret-bd-proper})
for the simultaneous exponential loss $L^{\text{r}}_0$.
\hfill$\blacksquare$\\[2mm]

\subsection{Proofs of results in Sections~\ref{sec:cost-transform}--\ref{sec:cost-indep}}

\noindent\textbf{Proof of Lemma~\ref{lem:scaling}.}\;
By definition,
\begin{align*}
& R_{\tilde L}(\eta,\gamma)= \sum_{j\in[m]} \eta_j c_{jM} L (j,\gamma) + \sum_{j\in[m]} \eta_j \sum_{k\in[m], k\not=j} (c_{jM} - c_{jk}) \{L(k,\gamma)-1\} \\
& = \sum_{j\in[m]} \eta_j c_{jM} L (j,\gamma) + \sum_{j\in[m]} \eta_j \sum_{k\in[m], k\not=j} (c_{jM} - c_{jk}) L(k,\gamma) - D(\eta),
\end{align*}
where $D(\eta) =  \sum_{j\in[m]} \sum_{k\in[m], k\not=j} \eta_j (c_{jM} - c_{jk})$.
By an exchange of indices $j$ and $k$, the second term above is
$\sum_{k\in[m]} \eta_k \sum_{j\in[m], j\not=k} (c_{kM} - c_{kj}) L(j,\gamma)$.
Substituting this into the preceding expression for $R_{\tilde L}(\eta,\gamma)$ yields
\begin{align*}
& R_{\tilde L}(\eta,\gamma)=  \sum_{j\in[m]} L(j,\gamma) \tilde \eta_j - D(\eta) = (1_m^\T \tilde \eta) R_L (\tilde{\tilde \eta},\gamma) - D(\eta) .
\end{align*}
The generalized entropy from $\tilde L$ is
\begin{align*}
H_{\tilde L} (\eta) = \inf_{\gamma}R_{\tilde L} (\eta,\gamma)  = (1_m^\T \tilde \eta) \inf_{\gamma}  R_L (\tilde{\tilde \eta},\gamma) - D(\eta) =  (1_m^\T \tilde \eta) H_L ( \tilde{\tilde \eta}) -D(\eta).
\end{align*}
The desired result on $B_{\tilde L}(\eta,\gamma)$ then follows.\hfill$\blacksquare$\\[2mm]

\noindent\textbf{Proof of Lemmas~\ref{lem:misclass} and \ref{lem:misclass2}.}\;
The bound in Lemma~\ref{lem:misclass} is a special case of Lemma~\ref{lem:misclass2} with $C = 1_m 1_m^\T - I_m$.
If $\eta=(\eta_1, \eta_2, 0, \ldots, 0)^\T$ and $q = (1/2, 1/2, 0, \ldots, 0)^\T$, then the bound becomes exact:
$ B^{\text{zo}} (\eta, q) = |2\eta_1 - 1|  = |\eta_1 - q_1| + |\eta_2-q_2|$ with $\eta_1+ \eta_2=1$.

For Lemma~\ref{lem:misclass2}, let $l= \argmax_{j \in[m]} (\overline C^\T_j \eta )$ and $k = \argmax_{j\in[m]} (\overline C^\T_j q)$,
where $\overline C= (\overline C_1, \ldots, \overline C_m)$ is a column representation of $\overline C$.
By definition, $\overline C_j = C_j - C_M$ and $\overline C^\T_j \eta = C_M^\T \eta - C_j^\T \eta$ for $j\in [m]$. Direct calculation yields
\begin{align*}
& B^{\text{cw}} (\eta, \overline C^\T q)  = R_{L^{\text{cw}}} (\eta, \overline C^\T q) - H^{\text{cw}} (\eta) \\
% & = \sum_{j=1}^m \eta_j c_{jl} - C_l^\T \eta \\
& = C_k^\T \eta - C_l ^\T \eta = \overline C_l^\T \eta - \overline C_k^\T \eta.
\end{align*}
Then $B^{\text{cw}} (\eta, \overline C^\T q)  = \overline C^\T_l \eta - \overline C^\T_k \eta \le \overline C^\T_l \eta - \overline C^\T_l q
+ \overline C^\T_k q - \overline C^\T_k \eta$ because $\overline C^\T_l q \le \overline C^\T_k q$ by definition.
Hence $B^{\text{cw}} (\eta, \overline C^\T q) \le | \overline C^\T_l \eta - \overline C^\T_l q| + |\overline C^\T_k q -\overline C^\T_k \eta| \le \|\overline C^\T ( \eta- q)\|_{\infty 2}$.
\hfill$\blacksquare$\\[2mm]

\noindent\textbf{Proof of Proposition~\ref{pro:regret-bd}.}\;
Note that $L^{\text{cw}}(j, \gamma) = \tilde L^{\text{zo}}(j,\gamma)$ by direct calculation.
Applying Lemma~\ref{lem:scaling} to $L$ and $L^{\text{zo}}$ shows that for any $\eta, q\in\Delta_m$,
\begin{align}
& B_{\tilde L} (\eta, q) = (1_m^\T \tilde \eta) B_L ( \tilde{\tilde \eta}, q), \label{eq:L-relation}\\
& B^{\text{cw}}(\eta, q) =(1_m^\T \tilde \eta) B^{\text{zo}}(\tilde{\tilde\eta}, q)  , \label{eq:zo-relation}
\end{align}
where $\tilde\eta$ and $\tilde{\tilde\eta}$ are defined as in Lemma~\ref{lem:scaling}.
The desired result then follows because
$B^{\text{zo}}(\tilde{\tilde\eta}, q) \le  \| \tilde{\tilde \eta}- q \|_{\infty2}$ by Lemma \ref{lem:misclass},
$\psi_q (  \| \tilde{\tilde \eta}- q \|_{\infty2} ) \le B_L(\tilde{\tilde \eta}, q)$ by definition, and $\psi_q(\cdot)$ is nondecreasing.
\hfill$\blacksquare$\\[2mm]

\noindent\textbf{Proof of Corollary~\ref{cor:regret-bd-proper}.}\;
Applying (\ref{eq:cw0-regret-bd2-proper}) with $C_0 = 1_m$ and $\eta$ replaced by $\tilde{\tilde\eta}$ yields
\begin{align*}
\underline \psi \left( B^{\text{zo}}(\tilde{\tilde\eta}, q) \right) \le  B_L (\tilde{\tilde\eta}, q).
\end{align*}
Combining this with (\ref{eq:L-relation}) and (\ref{eq:zo-relation}) gives the desired result.
\hfill$\blacksquare$\\[2mm]

\noindent\textbf{Proof of Proposition~\ref{pro:regret-bd2}.}\;
The desired result is obtained by combining the following observations: $B^{\text{cw}}(\eta, \overline C^\T q) \le \| \overline C^\T (\eta- q) \|_{\infty2}$ by Lemma~\ref{lem:misclass2},
$\psi_q^C ( \| \overline C^\T ( \eta- q) \|_{\infty2} ) \le B_L(\eta, q)$ by definition, and $\psi_q^C(\cdot)$ is nondecreasing.
\hfill$\blacksquare$\\[2mm]

\noindent\textbf{Proof of inequality (\ref{eq:cw-regret-bd3}).}\;
For any $w \in \mathcal W_{\eta,q}$, $\argmax_{j\in[m]} (\overline C^\T_j q^w)$ can be set to be same as $\argmax_{j\in[m]} (\overline C^\T_j q)$ and
hence $B^{\text{cw}}(\eta, \overline C^\T q^w) =  B^{\text{cw}}(\eta, \overline C^\T q)$.
Then inequality (\ref{eq:cw-regret-bd2}) with $q$ replaced by $q^w$ shows that
$\psi_{q^w}^C ( B^{\text{cw}}(\eta, \overline C^\T q) ) \le  B_L (\eta, q^w)$. The desired result then follows because  $B_L (\eta, q^w)\le B_L(\eta, q)$ by
the representation of $B_L(\eta,q)$ as  the Bregman divergence (\ref{eq:bregman}) and inequality (\ref{eq:breg-mon2}) in Lemma~\ref{lem:breg-mon}.
\hfill$\blacksquare$\\[2mm]

\noindent\textbf{Proof of Corollary~\ref{cor:regret-bd2-proper}.}\;
By the representation of $B_L(\eta,q)$ as  the Bregman divergence (\ref{eq:bregman})
and inequality (\ref{eq:breg-mon1}) in Lemma~\ref{lem:breg-mon},
$\psi_q^C(t)$ in Proposition~\ref{pro:regret-bd2} can be equivalently defined with $\eta^\prime \in \Delta_m$
restricted such that $\|\overline C^\T (\eta^\prime -q) \|_{\infty2}= t$.

We distinguish three cases. Let $k=\argmax_{j\in[m]} (\overline C^\T_j q)$.
First, if $\overline C^\T_k q > 1_m^\T \overline C^\T q /2$ and $\overline C^\T_k \eta > 1_m^\T \overline C^\T \eta /2$,
then $k=\argmax_{j\in[m]} (\overline C^\T_j \eta)$ and hence $B^{\text{cw}}(\eta, \overline C^\T q)=0$ and (\ref{eq:cw-regret-bd2-proper}) holds trivially.
Second, if $\overline C^\T_k q > 1_m^\T \overline C^\T q /2$ and $\overline C^\T_k \eta \le 1_m^\T \overline C^\T \eta /2$, then
(\ref{eq:cw-regret-bd3}) holds with some
$w \in \mathcal W_{\eta,q}$ such that $\overline C^\T_k q^w = 1_m^\T \overline C^\T  q^w /2$ and hence $\max_{j \in [m]} (\overline C^\T_j q^w) = 1_m^\T \overline C^\T q^w /2$,
because $\overline C^\T_k q^w  /(1_m^\T \overline C^\T  q^w)$ is continuous in $w \in [0,1]$, while taking a value $\le 1/2$ at $w=0$ and $>1/2$ at $w=1$ by assumption.
Third, if $\overline C^\T_k q \le 1_m^\T \overline C^\T q /2$, then  (\ref{eq:cw-regret-bd3}) holds with $w = 1 \in \mathcal W_{\eta,q}$.
In the latter two cases, inequality (\ref{eq:cw-regret-bd2-proper}) can be shown as follows:
\begin{align*}
& \underline \psi^C \left( B^{\text{cw}}(\eta, \overline C^\T q) \right)
\le \psi^C_{q^w} \left( B^{\text{cw}}(\eta, \overline C^\T q) \right) \le  B_L (\eta, q) .
\end{align*}
The first inequality holds because $\underline \psi^C (t)\le \psi_{q^w}^C(t)$ with $q^w$ satisfying $\max_{j \in [m]} (\overline C^\T_j q^w) \le 1_m^\T \overline C^\T  q^w /2$.
The second inequality holds by (\ref{eq:cw-regret-bd3}) with $w \in \mathcal W_{\eta,q}$.
\hfill$\blacksquare$\\[2mm]

\noindent\textbf{Proof of inequality (\ref{eq:cw0-regret-bd-2class}).}\;
For $\eta,q\in\Delta_2$, we have $\| C_0\circ(\eta - q)\|_{\infty2} =(c_{10}+c_{20}) |\eta_1 - q_1|$, where $\eta=(\eta_1,\eta_2)^\T$ and $q=(q_1,q_2)^\T$.
Moreover, $\max_{j=1,2} (c_{j0}q_j) \le C_0^\T q /2$ for $q\in \Delta_2$ leads to a single probability vector $q= (c_{20},c_{10})^\T /(c_{10}+c_{20})$.
From these expressions, $\underline \psi^{C_0}(t)$
can be simplified as $\underline \psi^{C_0}(t) = \min\{ \psi^{\text{RW}} (-t), \psi^{\text{RW}} (t)\}$.
\hfill$\blacksquare$\\[2mm]

\newpage

\begin{sidewaystable}[h]
\centering
    \renewcommand{\thetable}{S\arabic{table}(a)}

\begin{threeparttable}
\resizebox{\columnwidth}{!}{
\begin{tabular}{@{}llccc@{}}\toprule
         &   \multicolumn{1}{c}{Name}  & Loss $L(j,q)$ & Dissimilarity function $f(t)$ & Generalized Entropy $H(\eta)$ \\ \hline
\multicolumn{4}{l}{\small \underline{\textsc{Two Class Loss}}} \\
&Likelihood  & $-\mathbbm{1}_{\{j=1\}}\log q_1 - \mathbbm{1}_{\{j=2\}}\log q_2$        &       $t\log t -(t+1)\log(t+1)$  &   $-\eta_1\log \eta_1 - \eta_2\log\eta_2$       \\
&Exponential & $\mathbbm{1}_{\{j=1\}}\sqrt{\frac{q_2}{q_1}} + \mathbbm{1}_{\{j=2\}}\sqrt{\frac{q_1}{q_2}}$        &     $(\sqrt{t}-1)^2$    &    $-(\sqrt{\eta_1}-\sqrt{\eta_2})^2$      \\
&Calibration$_a$   &    $ \mathbbm{1}_{\{j=1\}}\frac{q_2}{2q_1}  + \mathbbm{1}_{\{j=2\}}\frac{1}{2}(\log\frac{q_1}{q_2} - 1)$   & $-\frac{1}{2}\log t$        & $\frac{\eta_2}{2}\log \frac{\eta_1}{\eta_2}$        \\
&Calibration$_s$   & $ \mathbbm{1}_{\{j=1\}} \frac{1}{2}(\log \frac{q_2}{q_1}+\frac{q_2}{q_1} -1)  + \mathbbm{1}_{\{j=2\}}\frac{1}{2}(\log\frac{q_1}{q_2} + \frac{q_1}{q_2} - 1)$        & $\frac{1}{2}(t\log t - \log t)$ & $\frac{1}{2}(\eta_1\log \frac{\eta_2}{\eta_1}+ \eta_2\log \frac{\eta_1}{\eta_2})$   \\

\multicolumn{4}{l}{\small \underline{\textsc{Multi-class Pairwise Asymmetric}}} \\
&Likelihood  & $\mathbbm{1}_{\{j\in[m-1]\}}\log(1+\frac{q_m}{q_j}) + \mathbbm{1}_{\{j=m\}}\sum_{i=1}^{m-1}\log(1+\frac{q_i}{q_m})$ &             $\sum_{i=1}^{m-1}\{t_i\log t_i -(1+t_i)\log(1+t_i)\}$ &   $-\sum_{i=1}^{m-1}(\eta_m\log \frac{\eta_m}{\eta_i+\eta_m} + \eta_i\log \frac{\eta_i}{\eta_i+\eta_m})$       \\
&Exponential & $\mathbbm{1}_{\{j\in[m-1]\}}(\sqrt{\frac{q_m}{q_j}}-1) +\mathbbm{1}_{\{j=m\}}\sum_{i=1}^{m-1}(\sqrt{\frac{q_i}{q_m}}-1) $       &  $\sum_{i=1}^{m-1}(\sqrt{t_i}-1)^2$    &       $-\sum_{i=1}^{m-1}(\sqrt{\eta_i}-\sqrt{\eta_m})^2$   \\
&Calibration   &    $\mathbbm{1}_{\{j\in[m-1]\}}\frac{q_m}{2q_j}+\mathbbm{1}_{\{j=m\}}\sum_{i=1}^{m-1}\frac{1}{2}(\log\frac{q_i}{q_m}-1)$     &     $-\sum_{i=1}^{m-1}\frac{1}{2} \log t_i$    &   $\sum_{i=1}^{m-1}\frac{\eta_m}{2}\log\frac{\eta_i}{\eta_m}$      \\

\multicolumn{4}{l}{\small \underline{\textsc{Multi-class Pairwise Symmetric}}} \\

&Likelihood  & $\sum_{i\neq j}2\log(1+\frac{q_i}{q_j})$ &             $-\sum_{i=1}^m\sum_{j\neq i}2t_i\log(1+\frac{t_j}{t_i})$ &   $\sum_{i=1}^m\sum_{j\neq i}2\eta_i\log(1+\frac{\eta_j}{\eta_i})$      \\
&Exponential & $\sum_{i\neq j}2(\sqrt{\frac{q_i}{q_j}}-1)$       &  $\sum_{i=1}^m\sum_{j\neq i}(\sqrt{t_i}-\sqrt{t_j})^2$    &       $-\sum_{i=1}^m\sum_{j\neq i}(\sqrt{\eta_i}-\sqrt{\eta_j})^2$   \\
&Calibration   &    $\sum_{i\neq j}\frac{1}{2}(\log \frac{q_i}{q_j}+\frac{q_i}{q_j}-1)$     &        $-\sum_{i=1}^m\sum_{j\neq i}\frac{t_i}{2} \log\frac{t_j}{t_i}$  &   $\sum_{i=1}^m\sum_{j\neq i}\frac{\eta_i}{2} \log\frac{\eta_j}{\eta_i}$      \\

\multicolumn{4}{l}{\small \underline{\textsc{Multi-class Simultaneous}}} \\

&$L_\beta$ Family & $({m^{\frac{1}{\beta}-1}-1})^{-1}[\{(1 + \sum_{i\neq j} (\frac{q_i}{q_j})^{\beta}\}^{\frac{1}{\beta}-1}-1]$ &  $-({m^{\frac{1}{\beta}-1}-1})^{-1}\{1+\sum_{i=1}^{m-1}t_i^\beta)^{\frac{1}{\beta}} -t_{\bullet}\}$ &   $({m^{\frac{1}{\beta}-1}-1})^{-1}\{(\sum_{i=1}^m\eta_i^\beta)^{\frac{1}{\beta}}-1\}$      \\
&Pairwise Exp($\beta=0$) & $(\prod_{i\neq j} \frac{q_i}{q_j})^{\frac{1}{m}}$       &  $-m(\prod_{i=1}^{m-1} t_i)^{\frac{1}{m}}$    &       $m(\prod_{i=1}^{m} \eta_i)^{\frac{1}{m}}$   \\
&Simulteneous Exp($\beta=\frac{1}{2}$)   &   $(m-1)^{-1}\sum_{i\neq j}\sqrt{\frac{q_i}{q_j}}  $   &        $-(m-1)^{-1}\{(1+\sum_{i=1}^{m-1}\sqrt{t_i})^2-t_{\bullet}\}$  &   $(m-1)^{-1}\{(\sum_{i=1}^{m}\sqrt{\eta_i})^2-1\}$       \\
&{Multinomial Lik($\beta=1$)}   & $-(\log m)^{-1}\log(q_j)$ &             $(\log m)^{-1}\sum_{i=1}^m t_i\log \frac{t_i}{t_{\bullet}}$ &   $-(\log m)^{-1}\sum_{i=1}^m \eta_i \log \eta_i$     \\
\bottomrule
\end{tabular} 
}

\begin{tablenotes}[para]
  {\scriptsize Note: $t_m=1$ and $t_{\bullet}=\sum_{i=1}^m t_i$. Calibration$_a$ and Calibration$_s$ are the asymmetric and symmetric versions.}
  \end{tablenotes}
  \end{threeparttable}
   \caption{Examples of losses, dissimilarity functions, and generalized entropies} \label{tab:1}
\end{sidewaystable}

\begin{table}[h]  \centering
 \addtocounter{table}{-1}
  \renewcommand{\thetable}{S\arabic{table}(b)}
    \begin{tabular}{@{}llc@{}}\toprule
\rule{0pt}{5pt}         &   \multicolumn{1}{c}{Name}  & Gradients $\frac{\partial}{\partial h_l} L(j, q^h)$  \\ \hline
\multicolumn{3}{l}{\small \underline{\textsc{Two Class Loss}}} \\
&Likelihood  & $(\mathbbm{1}_{\{j=1\}}q_2 - \mathbbm{1}_{\{j=2\}}q_1)(-1)^{\mathbbm{1}_{\{l=1\}}}$        \\
&Exponential & $\frac{1}{2}\big(\mathbbm{1}_{\{j=1\}}\sqrt{\frac{q_2}{q_1}} + \mathbbm{1}_{\{j=2\}}\sqrt{\frac{q_1}{q_2}}\big)(-1)^{\mathbbm{1}_{\{l=1\}}}$       \\
&Calibration$_a$   &    $ \frac{1}{2}\big(\mathbbm{1}_{\{j=1\}}\frac{q_2}{q_1}  + \mathbbm{1}_{\{j=2\}}\cdot 1\big)(-1)^{\mathbbm{1}_{\{l=1\}}}$         \\
&Calibration$_s$   & $\frac{1}{2}\big\{ \mathbbm{1}_{\{j=1\}} \big(1 + \frac{q_2}{q_1}\big)  + \mathbbm{1}_{\{j=2\}}\big(1+\frac{q_1}{q_2}\big) \big\}(-1)^{\mathbbm{1}_{\{l=1\}}}$        \\

\multicolumn{3}{l}{\small \underline{\textsc{Multi-class Pairwise Asymmetric}}} \\

&Likelihood & $\begin{aligned}\begin{cases}-\mathbbm{1}_{\{j=l\}}\frac{q_m}{q_l+q_m}+\mathbbm{1}_{\{j=m\}}\frac{q_l}{q_l+q_m}, & l\in[m-1] \\
\mathbbm{1}_{\{j\neq l\}}\frac{q_m}{q_j+q_m}+\mathbbm{1}_{\{j=m\}}\sum_{i=1}^{m-1}\frac{q_i}{q_i+q_m},& l = m
\end{cases}\end{aligned}$          \\

&Exponential & $\begin{aligned}\begin{cases}\frac{1}{2}\big(-\mathbbm{1}_{\{j=l\}}\sqrt{\frac{q_m}{q_l}}+\mathbbm{1}_{\{j=m\}}\sqrt{\frac{q_l}{q_m}}\big), & l\in[m-1] \\
\frac{1}{2}\big(\mathbbm{1}_{\{j\neq l\}}\sqrt{\frac{q_m}{q_j}}-\mathbbm{1}_{\{j=m\}}\sum_{i=1}^{m-1}\sqrt{\frac{q_i}{q_m}}\big), & l = m
\end{cases}\end{aligned}$ \\
&Calibration   &  $\begin{aligned}\begin{cases}\frac{1}{2}\big(-\mathbbm{1}_{\{j=l\}}\frac{q_m}{q_l}+\mathbbm{1}_{\{j=m\}}\cdot 1\big),& l\in[m-1] \\
\frac{1}{2}\big(\mathbbm{1}_{\{j\neq l\}}\frac{q_m}{q_j}-\mathbbm{1}_{\{j=m\}}\cdot(m-1)\big), & l = m
\end{cases}\end{aligned}$   \\

\multicolumn{3}{l}{\small \underline{\textsc{Multi-class Pairwise Symmetric}}} \\

&Likelihood  & $2(\mathbbm{1}_{\{j\neq l\}} \frac{q_l}{q_l + q_j} - \mathbbm{1}_{\{j=l\}}\sum_{i\neq l}\frac{q_i}{q_i+q_l})$              \\
&Exponential &  $2 \big(\mathbbm{1}_{\{j\neq l\}}\sqrt{\frac{q_l}{q_j}} - \mathbbm{1}_{\{j=l\}}\sum_{i\neq l}\sqrt{\frac{q_i}{q_l}}\big)$      \\
&Calibration   &    $\frac{1}{2}\big(\mathbbm{1}_{\{j\neq l\}} \big(\frac{q_l}{q_j}+1\big) - \mathbbm{1}_{\{j=l\}}\sum_{i\neq l}\big({\frac{q_i}{q_l}}+1\big)\big)$      \\

\multicolumn{3}{l}{\small \underline{\textsc{Multi-class Simultaneous}}} \\

&$L_\beta$ Family & $\frac{1-\beta}{m^{\frac{1}{\beta}-1}-1}(\sum_{i=1}^mq_i^\beta)^{\frac{1}{\beta}-2}\big(
\mathbbm{1}_{\{j\neq l\}}q_j^{\beta-1}q_l^{\beta} -\mathbbm{1}_{\{j=l\}}\sum_{i\neq l} q_{l}^{\beta -1}q_i^\beta\big)$ \\
&Pairwise Exp($\beta=0$) & $\frac{1}{m} \mathbbm{1}_{\{j\neq l\}}(\prod_{i\neq j} \frac{q_i}{q_j})^{\frac{1}{m}} - \frac{m-1}{m}\mathbbm{1}_{\{j=l\}}(\prod_{i\neq l} \frac{q_i}{q_l})^{\frac{1}{m}} $      \\
&Simultaneous Exp($\beta=\frac{1}{2}$)   &   $\frac{1}{2(m-1)}\big(\mathbbm{1}_{\{j\neq l\}} \sqrt{\frac{q_l}{q_j}} - \mathbbm{1}_{\{j=l\}}\sum_{i\neq l}\sqrt{\frac{q_i}{q_l}}\big)$  \\
&{Multinomial Lik($\beta=1$)}   & $\frac{1}{\log m}\left(\mathbbm{1}_{\{j\neq l\}}q_l + \mathbbm{1}_{\{j=l\}}(q_l-1)\right)$             \\
\bottomrule
\end{tabular} 
    \caption{Examples of losses and gradients}\label{tab:2}
\end{table}

\end{document}